\RequirePackage{fix-cm}
\documentclass[mnsc,nonblindrev]{informs3_hide}
\OneAndAHalfSpacedXI
\usepackage{amssymb}
\usepackage{amsfonts}
\usepackage{float}
\usepackage{placeins}
\usepackage{algorithm,algpseudocode}
\floatstyle{boxed}
\restylefloat{algorithm}
\usepackage{setspace}
\usepackage{lmodern}
\usepackage{microtype}
\usepackage{graphicx}
\usepackage[breakable,skins]{tcolorbox}
\newtcolorbox{promptbox}{
  enhanced, breakable,
  colback=black!2, colframe=black!30,
  boxrule=0.4pt, arc=0pt,
  left=8pt, right=8pt, top=7pt, bottom=7pt,
  before skip=6pt, after skip=10pt,
  fontupper=\small,
  before upper={\setlength{\parindent}{0pt}\setlength{\parskip}{3pt}%
    \renewcommand{\baselinestretch}{1.05}\selectfont\raggedright}
}
\usepackage{natbib}
\usepackage{bm,mathtools}
\bibpunct[, ]{(}{)}{,}{a}{}{,}
\def\bibfont{\small}

\TheoremsNumberedThrough
\ECRepeatTheorems
\theoremstyle{TH}

\EquationsNumberedThrough
\usepackage[colorlinks=true,linkcolor=blue,urlcolor=blue,citecolor=blue]{hyperref}

\let\InformsProof\proof
\let\EndInformsProof\endproof
\renewenvironment{proof}[1][Proof]
  {\InformsProof{#1.}}
  {\unskip\nobreak\hfill\Halmos\EndInformsProof}

\newcommand{\E}{\mathbb E}
\newcommand{\Prb}{\mathbb P}

\newcommand{\ind}{\mathbf 1}

\newcommand{\SetCompanionAnchors}{%
  \setcounter{observation}{0}%
  \renewcommand{\theobservation}{EC.\arabic{observation}}%
  \def\theHsection{EC.\arabic{section}}%
  \def\theHsubsection{\theHsection.\arabic{subsection}}%
  \def\theHsubsubsection{\theHsubsection.\arabic{subsubsection}}%
  \def\theHequation{EC.\arabic{equation}}%
  \def\theHfigure{EC.\arabic{figure}}%
  \def\theHtable{EC.\arabic{table}}%
  \def\theHtheorem{EC.\arabic{theorem}}%
  \def\theHlemma{EC.\arabic{lemma}}%
  \def\theHobservation{EC.\arabic{observation}}%
}
\MANUSCRIPTNO{}

\begin{document}
\RUNAUTHOR{}
\RUNTITLE{Agentic Discovery of Simple Symbolic Inventory Policies}
\TITLE{Agentic Discovery of Simple Symbolic Inventory Policies}
\ARTICLEAUTHORS{%
 \AUTHOR{Xi Chen\footnotemark[1]}
 \AFF{Leonard N. Stern School of Business, New York University, New York, NY 10012, USA, \EMAIL{xc13@stern.nyu.edu}}
 \AUTHOR{Yuze Chen\footnotemark[1]}
 \AFF{Qiuzhen College, Tsinghua University, Beijing 100084, China, \EMAIL{yz-chen21@mails.tsinghua.edu.cn}}
\AUTHOR{Baichuan Mo\footnotemark[1]}
\AFF{Department of Civil Engineering, Tsinghua University, Beijing 100084, China, \EMAIL{bmo@tsinghua.edu.cn}}
\AUTHOR{Jie Ying\footnotemark[1]}
\AFF{Qiuzhen College, Tsinghua University, Beijing 100084, China, \EMAIL{yingj24@mails.tsinghua.edu.cn}}
 \AUTHOR{Yuan Zhou\footnotemark[1]}
 \AFF{Yau Mathematical Sciences Center \& Department of Mathematical  Sciences, Tsinghua University, Beijing 100084, China; Beijing Institute of Mathematical Sciences and Applications, Beijing 101408, China \EMAIL{yuan-zhou@tsinghua.edu.cn}}
}
\renewcommand{\thefootnote}{\fnsymbol{footnote}}
\footnotetext[1]{Author names listed in alphabetical order.}

\RRHFirstLine{{\it\theRUNTITLE}}
\LRHFirstLine{{\it\theRUNTITLE}}
\ECAUpunct{}
\ABSTRACT{Designing simple inventory policies becomes increasingly difficult as operational features and system configurations expand. We develop an Agentic Inventory-Policy Discovery System based on explicitly defined simple operational formula (SOF) policy and template classes. Large language models (LLMs) revise symbolic structures, while a numerical engine fits parameters separately for each instance and evaluates the induced policies. Multicriteria parent selection, state-specific operational diagnostics, and reusable symbolic motifs guide structural evolution. A three-stage progression develops shared templates across instances and configurations with varying state dimensions. Experiments on 15,600 instances spanning 50 demand distributions cover lead times, perishability, their combination, and lead times with limited customer patience. The discovered templates improve aggregate empirical performance over evaluated SOF benchmarks and are competitive with policies requiring more elaborate online computations. For a lost-sales system with a three-period lead time, we develop an average-cost decomposition and prove that, under specified demand and cost conditions, a capped dual-shortfall specialization of a discovered policy strictly outperforms every policy in the considered capped base-stock benchmark class. This analysis links the performance advantage to the discovered replenishment mechanism. Together, these results highlight the potential of agentic discovery to generate simple inventory policies with strong empirical performance and identify operational mechanisms amenable to rigorous theoretical analysis.
}

\maketitle

\section{Introduction} \label{sec:intro}

Online decision-making in operations requires translating an evolving system state into actions whose consequences unfold over time. Inventory replenishment illustrates this challenge: an ordering decision must account for available stock, outstanding supply, and future demand, while balancing holding costs against shortages and potential expiration. Dynamic programming provides the standard framework for determining an optimal policy, but its computational burden grows rapidly as the state records more pipeline positions, inventory-age cohorts, and outstanding customer requests. Simple symbolic policies offer a practical alternative. Constant-order, base-stock, and capped base-stock policies show how simple and concise rules based on replenishment rates, inventory-position targets, and ordering caps can achieve strong performance and good interpretability \citep{goldberg2016constant,huh2009orderupto,xin2021capped}. Their success motivates the search for compact formulas that retain the state information most useful for controlling inventory costs.

Constructing such formulas manually becomes increasingly difficult as operational models incorporate additional features. Lead times, perishability, and customer patience introduce interactions among supply arrival, product expiration, and the timing of demand fulfillment. An effective compact policy must capture these interactions through a small number of state summaries and ordering adjustments. Moreover, favorable asymptotic behavior does not by itself establish strong empirical performance at moderate parameter scale. Designing policies that perform well across these ranges often involves substantial manual trial and error: researchers propose candidate formulas, tune their parameters, evaluate their performance, and revise their structures. This process becomes more demanding as the range of demand distributions, cost vectors, and system configurations expands.

Automatic approaches provide an alternative to manual policy construction. Earlier work uses genetic programming to evolve inventory-control policies~\citep{kleinau2004genetic} and interpretable symbolic decision rules~\citep{hein2018interpretable}. Recent preprints have explored the use of large language models (LLMs) to discover inventory policies. In concurrent work, \citet{yang2026aips} develop AIPS, which combines evolutionary proposals of policy structures with numerical parameter optimization, while \citet{baek2026llms} investigates autonomous algorithm design through tool-assisted development and testing. InvEvolve likewise combines policy generation with iterative evaluation and refinement~\citep{huang2026invevolve}. These studies demonstrate the potential of LLMs to construct effective policies using computational feedback, but leave scope for a more systematic connection between evaluation experience and improvements to specific replenishment mechanisms. Overall inventory cost aggregates the consequences of ordering
decisions across states that differ in stock availability, supply arrival timing, and risks of expiration or customer abandonment. However, aggregate cost alone may not provide the LLM with sufficient guidance to identify which state-specific ordering decisions warrant revision or which policy components encode useful operational relationships. A systematic evolutionary procedure may therefore connect evaluation evidence to revisions of specific ordering responses and preserve useful operational constructions across iterations. This process also calls for an explicitly defined symbolic policy class that specifies the available components and their admissible combinations, providing a common representation for guided structural edits and the reuse of discovered subexpressions.

\subsection{Our contributions}

We develop an \emph{Agentic Inventory-Policy Discovery System} for a unified inventory model incorporating replenishment lead times, perishability, and limited customer patience. The system combines explicit symbolic policy classes with an evolutionary procedure that uses evaluation evidence and reusable subexpressions to guide discovery. We evaluate the resulting policies across a broad collection of inventory instances and theoretically investigate one of their replenishment mechanisms. Specifically, we make the following four contributions.

\begin{enumerate}
\item \textbf{Explicit SOF policy and template classes.}
We define the \emph{simple operational formula} (SOF) policy class through a symbolic grammar specifying admissible inputs, elementary operators, and their compositions. We further define SOF templates that share symbolic structure across \emph{system configurations}, each specified by a replenishment lead time, product shelf life, and customer patience horizon. A template combines an ordering expression with state summaries whose symbolic weight functions aggregate state variables across pipeline positions, inventory-age cohorts, and backlog cohorts, accommodating changes in state dimension.  Each template uses a fixed number of parameter symbols across configurations, while their numerical values are fitted separately for each demand--cost instance. These classes provide a common representation for structural edits, motif reuse, and complexity assessment within the evolutionary discovery procedure introduced below.

\item \textbf{Agentic evolution with multicriteria selection, guidance,  motif reuse, and progressive generalization.}  We develop an evolutionary procedure in which an LLM proposes structural edits and a numerical engine fits parameters and evaluates candidate policies. Our evolutionary procedure uses multicriteria selection to assemble a parent set that balances policy performance, simplicity, behavioral diversity, and the potential to generate useful improvements, preserving alternative replenishment mechanisms for refinement and recombination. The procedure also constructs state-specific guidance from state-specific dynamic-programming decision regrets and disagreement among parent actions. These diagnostics highlight costly decisions and contrasting replenishment responses to the same inventory, pipeline, and backlog conditions, helping the agent identify concrete targets for structural revision. Furthermore, a reflection agent extracts symbolic subexpressions from evaluated candidates and adds validated motifs to a library, preserving promising operational constructions that extend standard inventory-position rules for reuse and recombination in subsequent proposals. Finally, a three-stage design progressively broadens discovery from individual instances to multiple demand distributions and cost vectors at a fixed configuration, and then to multiple configurations. The procedure uses promising structures from smaller evaluation sets to initialize broader searches, and jointly revises ordering expressions and state-summary definitions in the final stage.

\item \textbf{Compact policy templates with strong empirical performance.} We apply the system to four settings: lead times, perishability, joint lead times and perishability, and lead times with limited customer patience. Experiments on 15,600 instances spanning 50 demand distributions identify compact templates that achieve better aggregate performance than benchmark policies within the SOF class, including recent LLM-discovered rules, and remain competitive with policies outside this class that require more elaborate state-dependent calculations. The discovered formulas provide interpretable refinements to familiar ordering rules. For example, in the lead-time setting, the discovered \emph{scaled capped dual-shortfall} (SCDS) policy and its variant supplement the inventory-position shortfall with an additional ordering limit based on the most recently placed order or a weighted summary of later-arriving supply. These limits can reduce inventory accumulation from consecutive orders. In the perishable setting, remaining-life summaries distinguish how long different inventory cohorts remain available to satisfy demand.  These findings demonstrate that our agentic discovery system can produce compact, interpretable policy templates that improve the empirical trade-off between inventory cost and policy complexity across diverse operational settings.

\item \textbf{Theoretical insights into a discovered replenishment mechanism.}
We examine the additional replenishment limit that the discovered SCDS policy imposes based on the most recently placed order. For a lost-sales system with a three-period lead time, we develop an average-cost decomposition linking individual ordering changes to long-run holding and lost-sales costs. Using this decomposition, we prove that, under specified demand and cost conditions, a capped dual-shortfall specialization achieves a strictly lower average cost than every policy in an extended CBS benchmark class. Under our integer-order formulation, this class includes the standard fixed-cap CBS policies studied in the literature~\citep{xin2021capped}. The analysis also explains how delaying replenishment after a preceding order can yield holding-cost savings that outweigh the additional lost-sales costs. Although the analysis focuses on a simplified setting, it provides rigorous support for the discovered ordering mechanism and highlights the potential for establishing broader theoretical guarantees for the discovered policies in future work.
\end{enumerate}

\section{Related work}\label{sec:relatedliterature}

\noindent\underline{\bf Lost-sales inventory systems with positive lead times.} For lost-sales inventory systems with positive lead times, early work establishes cost bounds and develops near-myopic policies~\citep{morton1969bounds,morton1971near}, while subsequent comparisons and reviews examine a broader range of
approximations~\citep{zipkin2008old,bijvank2011lostsales}.
One line of research establishes performance guarantees for simple ordering rules in the classical nonperishable model. Base-stock policies are asymptotically optimal as the lost-sales penalty grows relative to the holding cost~\citep{huh2009orderupto}, whereas constant-order policies are asymptotically optimal as the lead time grows~\citep{goldberg2016constant}, with an exponentially decaying optimality gap~\citep{xin2016exponential}. Capped base-stock (CBS) policies combine these structures through an inventory-position target and an order cap. \citet{xin2021capped} establishes their long-lead-time asymptotic optimality, and \citet{xin2026approximation} provides a $7/3$ approximation guarantee for the best CBS policy. A complementary line constructs ordering rules through state-dependent expectations and probabilities. Dual-balancing policies balance expected marginal holding and lost-sales costs and achieve a two-approximation guarantee~\citep{levi2008lostsales}. Projected-inventory-level (PIL) policies target expected inventory upon order arrival;
\citet{vanjaarsveld2024projected} establish their asymptotic
optimality in a high-penalty regime and, for exponential demand,
a long-lead-time regime. Fixed non-stockout-probability policies
target a prescribed probability of avoiding a stockout in the
order-arrival period~\citep{dekok2024fixed}. In their standard formulations, these expectation- and probability-based policies fall outside our SOF class because their state-dependent calculations extend beyond its admissible elementary operations. Our discovered templates are only modestly more complex than traditional constant-order, base-stock, and capped base-stock rules, yet achieve substantially better aggregate empirical performance than these baselines and remain competitive with policies requiring more elaborate online computations.

\medskip
\noindent\underline{\bf Perishable inventory systems.} For perishable inventory systems, early dynamic-programming analyses characterize optimal ordering for products with fixed lifetimes~\citep{fries1975optimal,nahmias1975optimal}, followed by myopic approximations~\citep{nahmias1976myopic}. Subsequent work develops simple ordering rules that account for inventory levels and stock-age information. Stock-level-dependent policies combine
base-stock and constant-order behavior~\citep{haijema2013new,haijema2016stock}, while age weights and expected-waste adjustments further refine replenishment decisions~\citep{haijema2019improved}. For a zero-lead-time FIFO system, \citet{bu2023basestock} prove asymptotic optimality of base-stock policies in several regimes, including increasing product lifetime and demand population size. Balancing policies provide worst-case guarantees~\citep{chao2015approximation}, with truncated balancing addressing high shortage penalties~\citep{zhang2023truncated}. Positive lead times introduce additional interactions between supply arrival and product expiration. In this setting, \citet{broekmeulen2009heuristic} develop an expected-waste-based heuristic for batch ordering and time-varying demand, \citet{chao2018capacitated} establish approximation guarantees under finite ordering capacity, and \citet{bu2026perishable} analyze inventory-position and projected-inventory-level policies. We apply our agentic policy discovery system to perishable systems both without and with replenishment lead times. Our discovered SOF templates account for how long existing stock remains usable and when outstanding supply arrives, refining replenishment beyond aggregate inventory-position targets. The resulting rules achieve empirical performance gains over existing policies, including more computationally involved alternatives, without requiring online expected-waste or inventory-projection calculations.

\medskip
\noindent\underline{\bf Inventory systems with limited customer patience.} Inventory systems with limited customer patience combine backordering and lost sales: customers may initially wait for service but abandon before their demand is fulfilled. \citet{lim2024partial} study a partial-backorder model with geometrically distributed patience. When demand and patience distributions are known and backlog is observable, they establish base-stock optimality at zero lead time and uniform asymptotic optimality as the sales price grows at positive lead times. They also develop $(s,q)$ policies for unobservable backlog and an online learning algorithm with regret guarantees when demand and patience parameters are unknown. \citet{bu2026patience} consider random lead times and waiting-time-dependent customer retention, establishing optimality and asymptotic optimality results for base-stock policies under penalty scaling. We apply our discovery system to settings with observable backlog and fixed lead times and patience horizons. Our discovered SOF templates temper the replenishment response to older backlog and account for supply arrival timing, recognizing that some waiting demand may abandon before newly ordered inventory arrives. These families also recover base-stock behavior so their best achievable average costs inherit the asymptotic optimality guarantees of \citet{bu2026patience}. Their additional structural flexibility yields substantial aggregate empirical improvements over benchmark policies across different lead times and patience horizons.

\medskip
\noindent\underline{\bf Automatic policy discovery.} Genetic programming has been used to discover inventory algorithms
and replenishment rules~\citep{kleinau2004genetic,moore2005reordering},
including the joint evolution of expressions for $(r,Q)$
parameters~\citep{lopes2020cooperative} and replenishment and
transshipment decisions~\citep{xu2026niching}. Beyond inventory
control, \citet{maes2012formulas} search over short action-ranking
formulas using bandit-based evaluation allocation, while
\citet{hein2018interpretable} derive policy equations through
genetic programming and model-based batch reinforcement learning.
Neural approaches guide program search using a trained policy
as an oracle~\citep{verma2018programmatically}, or train an
expression generator to maximize the returns of symbolic
policies~\citep{landajuela2021symbolic}.

Broader research on LLM-based program and heuristic discovery
combines generation with numerical evaluation
\citep{romeraparedes2024funsearch,liu2024eoh}, reflective
feedback~\citep{ye2024reevo}, and multi-objective
selection~\citep{yao2025meoh}. \citet{chen2026teacher} additionally
use learned policies' action preferences as local feedback on
states visited by candidate programs. Related symbolic regression
methods search expression structures and numerical parameters
to fit observed relationships~\citep{cranmer2023pysr,shojaee2025llmsr};
LaSR further maintains a library of learned textual concepts
to guide search~\citep{grayeli2024concept}. These methods share
expression-search mechanisms with symbolic policy discovery,
although direct policy search evaluates expressions through
their sequential decision consequences.

Recent preprints explore LLM-based inventory policy discovery~\citep{yang2026aips,baek2026llms,huang2026invevolve}, as discussed in the introduction. In comparison, our work develops a more explicit connection between policy evaluation and the revision and reuse
of symbolic components.  State-specific guidance highlights costly ordering decisions, and contrasting replenishment responses to the same inventory, pipeline, and backlog conditions. The motif library preserves validated subexpressions from evaluated policies, allowing operational constructions such as shortfall-based adjustments, ordering limits, and state summaries to serve as reusable building blocks. The SOF representation gives these operations common structural units and supports consistent complexity assessment. Together, these designs guide structural revision toward specific replenishment responses while preserving useful operational constructions across iterations. We apply the system across a broader range of inventory settings (lead times, perishability, joint lead times and perishability, and lead times with limited customer patience) and
discover simple, operationally meaningful symbolic policies that achieve strong empirical performance.

\section{Problem Formulation and Preliminaries}\label{sec:problem-and-preliminaries}

We present our inventory model and problem formulation in
Section~\ref{sec:inventory-model-and-problem} and define the
\emph{simple operational formula} (SOF) policy class and SOF templates
in Section~\ref{sec:sec:SOF-policy-class}.

\subsection{Inventory model and problem formulation}\label{sec:inventory-model-and-problem}
We study a periodic-review, single-item inventory system with \emph{i.i.d.}~demand, deterministic lead times, perishability, and backlogging subject to a fixed patience horizon. Unmet demand remains eligible for service until its patience horizon expires, when any unserved quantity is permanently lost. The model includes immediate replenishment, nonperishable inventory, and immediate lost sales as special cases. We seek deterministic ordering policies that minimize finite-horizon or long-run average expected costs, accounting for inventory holding, carried backlog, lost demand, expiry, and ordering.

\medskip
\noindent\underline{\bf Model parameters.} The model is specified by
\begin{equation}
 \Theta=(\bm{\sigma},\mathcal{D},\bm{c}),
 \label{eq:model-parameters}
\end{equation}
where the configuration vector $\bm{\sigma}=(\sigma^{\rm lead},\sigma^{\rm life},\sigma^{\rm back})$ comprises the deterministic lead time $\sigma^{\rm lead}\in\mathbb{Z}_{\geq 0}$, the shelf life $\sigma^{\rm life}\in\mathbb{Z}_{+}\cup\{\infty\}$, and the patience horizon $\sigma^{\rm back}\in\mathbb{Z}_{\geq 0}$ for backlogged demand. The special cases $\sigma^{\rm lead}=0$, $\sigma^{\rm life}=\infty$, and $\sigma^{\rm back}=0$ correspond to immediate order arrival, nonperishable inventory, and lost sales, respectively. Demand is \emph{i.i.d.}~across periods with distribution $\mathcal{D}$ on the nonnegative integers. The cost vector $\bm{c}=(c^{\rm hold},c^{\rm backlog},c^{\rm lost},c^{\rm expiry},c^{\rm order})\ge\bm{0}$ assigns unit costs to holding, backlogged demand, lost demand, expiry, and ordering, respectively.

\medskip
\noindent\underline{\bf System state.}
Let $s_t$ denote the system state at the beginning of period $t$, belonging to the configuration-dependent state space $\mathcal S_{\bm\sigma}$:
\begin{equation}
 s_t=(\bm{x}_t,\bm{p}_t,\bm{b}_t)
 \in\mathcal S_{\bm\sigma}:=\mathbb Z_{\geq 0}^{\sigma^{\rm life,\natural}}\times\mathbb Z_{\geq 0}^{\sigma^{\rm lead}}\times\mathbb Z_{\geq 0}^{\sigma^{\rm back}}.
 \label{eq:system-state}
\end{equation}
We set $\sigma^{\rm life,\natural}=\sigma^{\rm life}$ for perishable inventory and $\sigma^{\rm life,\natural}=1$ for nonperishable inventory.
The vectors $\bm{x}_t=(x_{t,1},\dots,x_{t,\sigma^{\rm life,\natural}})$, $\bm{p}_t=(p_{t,1},\dots,p_{t,\sigma^{\rm lead}})$, and $\bm{b}_t=(b_{t,1},\dots,b_{t,\sigma^{\rm back}})$ record on-hand inventory, outstanding orders, and backlogged demand, respectively.
For perishable inventory, $x_{t,m}$ records the quantity with $m$ periods of remaining life; for nonperishable inventory, $x_{t}:= x_{t,1}$ records the total on-hand quantity.
The entry $p_{t,r}$ records supply arriving in $(r-1)$ periods (i.e., period $(t + r-1)$), and $b_{t,\kappa}$ records the quantity in backlog cohort $\kappa$ (i.e., the unmet demand originating in period $t-\kappa$).
A vector with dimension zero is empty.
When writing policy formulas, we use $s=(\bm x,\bm p,\bm b)\in\mathcal S_{\bm\sigma}$ as a generic state input, with the time index suppressed. 

\medskip
\noindent\underline{\bf Decision and transition.}
At the beginning of period $t$, the decision maker observes the system state $s_t$ and selects an order quantity $a_t=\pi(s_t, t)$ according to a deterministic inventory policy $\pi:\mathcal S_{\bm\sigma}\times\mathbb Z_{\geq0}\rightarrow\mathbb Z_{\geq0}$. The system then receives the oldest pipeline order or, when $\sigma^{\rm lead}=0$, the newly placed order. After demand $D_t\sim \mathcal{D}$ is realized, the system serves backlogged demand from the oldest cohort to the youngest, followed by current demand, using inventory in increasing order of remaining shelf life. The resulting updates to on-hand inventory, backlog, and the order pipeline follow Eqs.~(\ref{eq:expired-item-updated-inventory-perishable},\ref{eq:expired-item-updated-inventory-nonperishable},\ref{eq:backlogged-demand-update},\ref{eq:order-pipeline-update}), as detailed below.

\paragraph{Update of on-hand inventory.} We write $[x]^+=\max\{x,0\}$. The order arriving in period $t$ is 
\begin{align}
R_t = \begin{cases} 
a_t, & \sigma^{\rm lead}=0\\
p_{t,1} & \sigma^{\rm lead}> 0
\end{cases} .
\label{eq:period-receipt}
\end{align}
Arriving units enter the full-shelf-life coordinate, while other $x_{t,m}$ remains unchanged. The inventory available in period $t$ is therefore
\begin{equation}
\widetilde x_{t,m}
=x_{t,m}+R_t \cdot \mathbf 1\{m=\sigma^{\rm life,\natural}\},
\qquad m=1,\ldots,\sigma^{\rm life,\natural}.
\label{eq:inventory-after-receipt}
\end{equation}
The quantity served is the smaller of available inventory and total demand, including current demand and outstanding backlog:
\begin{equation}
V_t=\min\left\{
\sum_{m=1}^{\sigma^{\rm life,\natural}}\widetilde x_{t,m},
D_t+\sum_{\kappa=1}^{\sigma^{\rm back}}b_{t,\kappa}
\right\}. \label{eq:total-quantity-served}
\end{equation}
Inventory is consumed in increasing order of remaining shelf life. After accounting for shorter-lived inventory, the residual quantity in cohort $m$ is
\begin{align}
&\overline x_{t,m}
=\left[
\widetilde x_{t,m}
-\left[V_t-\sum_{j=1}^{m-1}\widetilde x_{t,j}\right]^+
\right]^+,
\qquad m=1,\ldots,\sigma^{\rm life,\natural}. \label{eq:inventory-after-service}
\end{align}
For perishable items, the expired items and updated inventory become
\begin{align}
W_t = \overline{x}_{t,1}, \qquad x_{t+1, \sigma^{\rm life}} = 0, \qquad x_{t+1, m} = \overline{x}_{t, m+1}, \qquad m=1,\ldots,\sigma^{\rm life} - 1. \label{eq:expired-item-updated-inventory-perishable}
\end{align}
For nonperishable items, we directly have
\begin{align}
W_t = 0, \qquad x_{t+1,1}= \overline{x}_{t,1}. \label{eq:expired-item-updated-inventory-nonperishable}
\end{align}

\paragraph{Update of backlogged demand.} Let $b_{t,0}=D_t$ for notational convenience.
Under the stated service priority, demand is served from the oldest cohort to the youngest, yielding the unmet quantity in cohort $\kappa$:
\begin{equation}
\overline b_{t,\kappa}
=\left[
b_{t,\kappa}
-\left[V_t-\sum_{j=\kappa+1}^{\sigma^{\rm back}}b_{t,j}\right]^+
\right]^+,
\qquad \kappa=0,\ldots,\sigma^{\rm back}.
\label{eq:unserved-demand-cohort}
\end{equation}
The permanently lost demand and the updated backlogged demand are:
\begin{equation}
A_t=\overline b_{t,\sigma^{\rm back}},
\qquad
b_{t+1,\kappa}=\overline b_{t,\kappa-1},
\quad \kappa=1,\ldots,\sigma^{\rm back}. \label{eq:backlogged-demand-update}
\end{equation}

\paragraph{Update of order pipeline.} For $\sigma^{\rm lead}\geq1$, the new order enters the pipeline and the updated pipeline is:
\begin{align}
p_{t+1,r}=\begin{cases}
    p_{t,r+1}, & r=1,\ldots,\sigma^{\rm lead}-1\\
    a_t, & r = \sigma^{\rm lead}
\end{cases} . \label{eq:order-pipeline-update}
\end{align}
When $\sigma^{\rm lead}=0$, the pipeline is empty and the current order directly enters the inventory.

\medskip
\noindent\underline{\bf Cost and evaluation objective.} Recall that $W_t$ is expired inventory (Eqs.~(\ref{eq:expired-item-updated-inventory-perishable},\ref{eq:expired-item-updated-inventory-nonperishable})), and $A_t$ is permanently lost demand (Eq.~\eqref{eq:backlogged-demand-update}). Also let $X_{t+1}= \sum_{m=1}^{\sigma^{\rm life, \natural}}x_{t+1,m}$ be the inventory carried to the next period, let $B_{t+1} = \sum_{\kappa=1}^{\sigma^{\rm back}} b_{t+1,\kappa}$ be the backlogged demand retained for the next period. Using the cost vector $\bm{c}$ defined in model parameters, the realized one-period cost is
\begin{equation}
 G_{t}=c^{\rm hold} \cdot X_{t+1}+ c^{\rm backlog} \cdot B_{t+1}+ c^{\rm lost} \cdot A_t + c^{\rm expiry} \cdot W_t + c^{\rm order} \cdot a_t.
\label{eq:per-period-cost}
\end{equation}

For a fixed initial state $s_0=\bm{0}$, the finite-horizon (with horizon $H$) and average-cost objectives are
\begin{equation}
J_H(\pi;\Theta):=\frac1H\,\mathbb E_{\pi,\Theta} \left[\left.\sum_{t=0}^{H-1}G_t\,\right|\,s_0=\bm{0}\right],\qquad J(\pi;\Theta)
:=\limsup_{H\to\infty}J_H(\pi;\Theta).
\label{eq:cost-objectives}
\end{equation}
The decision maker seeks to minimize the finite-horizon or average-cost objective over deterministic policies.

\subsection{The simple operational formula (SOF) policies and policy templates}
\label{sec:sec:SOF-policy-class}

We focus on simple symbolic policies for their computational efficiency, ease of implementation, and transparent decision logic. We first define \emph{simple operational formula} (SOF) policies for a fixed system configuration, then introduce \emph{SOF templates}, which use a shared symbolic structure to generate SOF policies uniformly across system
configurations.

\medskip
\noindent\underline{\bf Definition of SOF policies.}
For a fixed configuration vector $\bm{\sigma}$, a \emph{simple operational
formula} (SOF) is a pair $\mathfrak{p}=(g,\bm{\alpha})$, where $g$ is a
finite scalar-valued symbolic expression and $\bm{\alpha}$ is a vector
of parameter symbols. The set of all corresponding SOFs are $\mathcal H_{\bm\sigma}^{\rm SOF}$ with $\mathfrak{p}\in\mathcal H_{\bm\sigma}^{\rm SOF}$.  The numerical values of parameters in $\bm{\alpha}$
are determined from the demand distribution $\mathcal{D}$ and cost vector $\bm{c}$
before deployment and remain fixed during deployment. The expression
is represented by a typed abstract syntax tree (AST) constructed from
the following terminals:
\begin{itemize}
\item \emph{State coordinates and time index:}
$x_m$ for $m=1,\ldots,\sigma^{\rm life,\natural}$,
$p_r$ for $r=1,\ldots,\sigma^{\rm lead}$, and
$b_\kappa$ for $\kappa=1,\ldots,\sigma^{\rm back}$.
These coordinates retain the meanings and orderings defined in the system state, with the time index suppressed. The time index $t$ itself can also be used as a terminal in the AST.
\item \emph{Constants:}
entries of the fitted parameter vector $\boldsymbol\alpha$, finite numerical constants, and finite instance statistics computed from $\mathcal{D}$ and $\bm c$ before deployment.
These statistics may include moments or quantiles, and all constant terminals are independent of the current state.
\end{itemize}

The admissible operators are addition, subtraction (negation), multiplication, the guarded ratio ($u/(1+|v|)$), absolute value, minimum, maximum, clipping, rounding, and the indicator $\mathbf1\{u>v\}$. Here $u$, $v$, and $w$ denote admissible scalar subexpressions, and the indicator takes values in $\{0,1\}$.
Clipping is defined by
$\operatorname{clip}(u;v,w)=\min\{\max\{u,v\},w\}$ for bounds satisfying $v\leq w$.

Given system state $s = (\bm{x}, \bm{p}, \bm{b})$, the expression $g$ and parameter vector $\bm{\alpha} = \bm{\alpha}(\mathcal{D}, \bm{c})$ determine the order quantity
\begin{equation}
\pi_{g}(s, t;\bm{\alpha})
=
\max\left\{
0,\left\lfloor g(s, t;\bm{\alpha})\right\rceil
\right\},
\label{eq:definition-sof}
\end{equation}
where $\lfloor\cdot\rceil$ denotes nearest-integer rounding with ties to the even integer, ensuring a nonnegative integer order quantity. We refer to the ordering rule represented by $g$ as an SOF policy, with its order quantity obtained through Eq.~\eqref{eq:definition-sof}.

\medskip
\noindent\underline{\textbf{Definition of SOF templates.}}
Let $\mathfrak{S}$ be a set of system configurations. A
\emph{simple operational formula (SOF) template} is a triple
$\mathfrak{t}=(g,\Phi,\bm{\alpha})$. The set of all corresponding SOF templates are $\mathcal H_{\mathfrak S}^{\rm template}$ with $\mathfrak{t}\in\mathcal H_{\mathfrak S}^{\rm template}$. The expression $g$ is represented by an AST using the operators permitted for SOF policies and the following terminals:
\begin{itemize}
\item \emph{State summaries and time index:}
The template uses a finite collection of state summaries. Each summary
$\varphi \in \Phi$ is a linear combination of the coordinates of
$s=(\bm{x},\bm{p},\bm{b})$:
\begin{align}
\varphi(s;\bm{\sigma}, \bm{\alpha}) = \sum_{m=1}^{\sigma^{\rm life, \natural}} w_\varphi^{\rm life}(m; \bm{\sigma}, \bm{\alpha}) x_m + \sum_{r=1}^{\sigma^{\rm lead}} w_\varphi^{\rm lead}(r; \bm{\sigma}, \bm{\alpha}) p_r + \sum_{\kappa=1}^{\sigma^{\rm back}} w_\varphi^{\rm back}(\kappa; \bm{\sigma}, \bm{\alpha}) b_\kappa,
\label{eq:state-summary}
\end{align}
Weight functions $w_\varphi^{\rm life}(\cdot)$, $w_\varphi^{\rm lead}(\cdot)$, and $w_\varphi^{\rm back}(\cdot)$ are finite symbolic expressions in their coordinate
index, $\bm{\sigma}$, and $\bm{\alpha}$, constructed using numerical
constants, arithmetic operations, floor and ceiling operators, and
indicator functions. The time index $t$ itself can also be used as a terminal in the AST.
\item \emph{Constants:}
Constant terminals are the same as those permitted in the
definition of SOF policies.
\end{itemize}

An SOF template $\mathfrak{t}=(g,\Phi,\bm{\alpha})$ defines a shared symbolic structure across configurations. For each
$\bm{\sigma}\in\mathfrak{S}$, expanding its state summaries according to \eqref{eq:state-summary} yields an expression $g_{\bm{\sigma}}$ in the
state coordinates of that configuration. The resulting SOF
$(g_{\bm{\sigma}},\bm{\alpha})$ is the configuration-specific
instantiation of the template $\mathfrak{t}$.

\medskip
\noindent\underline{\textbf{Examples and relation to existing policies.}}
The following SOF templates represent integer-valued base-stock (BS),
capped base-stock (CBS), and constant-order (CO) policies and illustrate
how state summaries, parameters, and admissible operators define an
ordering rule. For a state $s=(\bm{x},\bm{p},\bm{b})$, define the net
inventory position through the state summary
\begin{equation}
\varphi_{\rm IP}(s;\bm{\sigma})
=
\sum_{m=1}^{\sigma^{\rm life,\natural}}x_m
+
\sum_{r=1}^{\sigma^{\rm lead}}p_r
-
\sum_{\kappa=1}^{\sigma^{\rm back}}b_\kappa.
\label{eq:sof-inventory-position}
\end{equation}
For nonnegative integer parameter $S$ and nonnegative real parameter $C$, define
the symbolic expressions
\begin{align}
g_{\rm BS}=S-\varphi_{\rm IP},\qquad
g_{\rm CBS}=\min\{\xi(C,t),[S-\varphi_{\rm IP}]^+\},\qquad
g_{\rm CO}=\xi(C,t).
\label{eq:sof-bs-cbs-co}
\end{align}
Here, $S$ denotes the base-stock level, while $C$ denotes the order cap
for CBS and the constant order quantity for CO. The transformation $\xi(C, t) :=  \lfloor C(t+1)\rceil-\lfloor Ct\rceil$ accommodates fractional values of $C$ through integer quantities
in $\{\lfloor C\rfloor,\lceil C\rceil\}$ whose time average
converges to $C$.\footnote{For any fractional order rate \(C\), with \(C<\mathbb E[D]\) in the nonperishable case, the long-run average cost of ordering exactly \(C\) units per period can be approximated arbitrarily closely by the integer-valued ordering sequence \(\xi(C',t)\), where \(C'\) can be chosen arbitrarily close to \(C\). Thus, the family \(\{\xi(C',\cdot):C'\ge0\}\) provides integer-valued approximations to the long-run performance of fractional constant-order policies.} The corresponding SOF
templates are
$
\mathfrak{t}_{\rm BS}
=(g_{\rm BS},\{\varphi_{\rm IP}\},(S))$,
$\mathfrak{t}_{\rm CBS}
=(g_{\rm CBS},\{\varphi_{\rm IP}\},(C,S))$, and
$\mathfrak{t}_{\rm CO}
=(g_{\rm CO},\emptyset,(C))$.

To illustrate how an SOF template yields an SOF under a specific system configuration, consider $\bm\sigma=(2,\infty,0)$, corresponding to nonperishable inventory with a two-period lead time and lost sales. Under this configuration, the state summary expands to $\varphi_{\rm IP}(s;\bm\sigma)=x_1+p_1+p_2$, and the CBS template yields the configuration-specific SOF expression
$g_{{\rm CBS},\bm\sigma}(s, t;C,S)=\min\{\xi(C,t),[S-x_1-p_1-p_2]^+\}$ over the parameter vector $(C, S)$.

Other policies typically require additional state-dependent calculations.
Projected-inventory-level (PIL) policies
\citep{vanjaarsveld2024projected} determine orders using conditional
expectations of future inventory, while fixed non-stockout-probability
(FP3) policies \citep{dekok2024fixed} use projected service probabilities.
Myopic and Myopic-2 policies \citep{morton1971near, zipkin2008old} determine orders by minimizing expected costs over their respective look-ahead horizons.
These policies fall outside the SOF representation in their standard
computational form because conditional expectations, service-probability
evaluation and inversion, and expected-cost minimization are not
admissible primitives. These calculations may require distributional
recursions or numerical procedures. In contrast, compact SOFs use few elementary
operations, can be implemented directly from their expressions, and
allow managers to inspect and validate the ordering logic.

\section{Agentic Inventory-Policy Discovery System}\label{sec:methods}
We now present our agentic evolutionary algorithm for discovering compact SOF policies and templates with low inventory costs across a prescribed set of instances. An LLM agent proposes structural edits within the admissible symbolic classes to generate offspring, while a numerical engine fits each candidate's parameters and evaluates its induced policy separately for each instance. Three mechanisms organize the structural search: multicriteria parent selection balances cost and simplicity, preserves diverse ordering behaviors, and favors parents that generate useful offspring; state-specific operational diagnostics guide edits toward costly or contrasting replenishment responses; and a reflection agent maintains a library of validated symbolic subexpressions, preserving operational constructions for reuse and recombination across iterations. To manage the cost of searching across many instances, we embed this evolutionary search in a three-stage framework that progressively broadens evaluation from a single instance to multiple demand distributions and cost vectors at a fixed configuration, and then to multiple configurations. Each stage runs the same evolutionary procedure, with the later stages initialized from promising structures found in the preceding stage. In the final stage, the algorithm jointly revises ordering expressions and state summaries to discover shared templates that accommodate varying state dimensions.

In the rest of this section, Section~\ref{sec:methodoverview} first defines the search objects and their evaluation; Section~\ref{sec:methodclosedloop} then explains how the three mechanisms work together in the core evolutionary LLM search; Section~\ref{sec:methodstages} further describes the outer three-stage progression, and Section~\ref{sec:methodframework} finally integrates the inner search and outer progression into the complete algorithm. Implementation details appear in Section~\ref{app:discoveryimplementation}.

\subsection{Search objects and evaluation}
\label{sec:methodoverview}
A search run is specified by four inputs: the admissible candidate class $\mathcal H$, an evaluation set $\mathcal T$, a scoring rule $f$, and an initial population $\mathcal G^{\rm pop}_0\subseteq\mathcal H$. Its objective is to identify candidates whose induced ordering rules achieve low costs across the evaluation instances while keeping their symbolic structures simple. The finite, nonempty initial population $\mathcal G^{\rm pop}_0$ contains admissible candidates and may be inherited from an earlier stage (see Section~\ref{sec:methodstages} for details).

\medskip
\noindent\underline{\textbf{Candidate structures and evaluation instances.}}
Candidates are complete symbolic objects, which could be SOFs $\mathfrak p=(g,\bm\alpha)$ at a fixed configuration or SOF templates $\mathfrak t=(g,\Phi,\bm\alpha)$ across configurations. We use $\mathfrak c\in\mathcal H$ for a generic candidate of either type, reserving $g$ for its ordering expression and $\bm\alpha$ for its parameter vector.
For a fixed configuration $\bm\sigma$, the class $\mathcal H = \mathcal H_{\bm\sigma}^{\rm SOF}$ contains SOF pairs $\mathfrak p=(g,\bm\alpha)$ using that configuration's state coordinates. Across a configuration set $\mathfrak S$, the class $\mathcal H = \mathcal H_{\mathfrak S}^{\rm template}$ contains SOF templates $\mathfrak t=(g,\Phi,\bm\alpha)$. Instantiating a template at $\bm\sigma$ expands its summaries to produce the SOF pair $\mathfrak p_{\bm\sigma}=(g_{\bm\sigma},\bm\alpha)$.

The finite, nonempty evaluation set is $\mathcal T=\{\Theta_1,\ldots,\Theta_m\}$, where $\Theta_j=(\bm\sigma_j,\mathcal D_j,\bm c_j)$ is the $j$-th evaluation instance, $j=1,\ldots,m$. In the SOF search over $\mathcal H_{\bm\sigma}^{\rm SOF}$, $\bm\sigma_j=\bm\sigma$ for every instance, while demand distributions $\mathcal D_j$ and cost vectors $\bm c_j$ can differ. In the SOF-template search over $ \mathcal H_{\mathfrak S}^{\rm template}$, the instances can span configurations in $\mathfrak S$. A template candidate $(g,\Phi,\bm\alpha)$ is instantiated at $\bm\sigma_j\in \mathfrak S$ to obtain the SOF pair $(g_{\bm\sigma_j},\bm\alpha)$, whereas an SOF candidate $(g,\bm\alpha)$ is used directly, with $g_{\bm\sigma_j}=g_{\bm\sigma}=g$ for a unified evaluation notation. 

\medskip
\phantomsection
\noindent\underline{\textbf{From a structure to its search score.}}
\label{sec:methodfitting}
To evaluate a candidate $\mathfrak c$, the numerical engine fits the parameters $\bm\alpha$ of its SOF pair $(g_{\bm\sigma_j},\bm\alpha)$ separately for each $\Theta_j\in\mathcal T$, seeking to minimize the finite-horizon average cost $J_H$ in Eq.~\eqref{eq:cost-objectives}. We write $J^{\rm eval}$ for its numerical evaluation: it equals $J_H$ under exact state propagation and is a Monte Carlo estimate when the reachable state space is too large. The fitted vector $\widehat{\bm\alpha}(\mathfrak c,\Theta_j)$ contains the lowest-cost numerical values found within the allocated computational budget. Together with $g_{\bm\sigma_j}$, it determines the policy $\pi_{g_{\bm\sigma_j}}(\,\cdot\,;\widehat{\bm\alpha}(\mathfrak c,\Theta_j))$ through Eq.~\eqref{eq:definition-sof}.

We compare each candidate with a prespecified set of baseline policies (such as Base Stock and Capped Base Stock), each fitted separately on the same instance and evaluated under the same objective. Let $J^{\rm best}(\Theta_j)$ denote the lowest evaluated cost among these baselines on instance $\Theta_j$. Assuming $J^{\rm best}(\Theta_j)>0$, the fitted relative cost gap is
\begin{equation}
\Delta(\mathfrak c;\Theta_j)
=\frac{J^{\rm eval}\bigl(\pi_{g_{\bm\sigma_j}}(\,\cdot\,;\widehat{\bm\alpha}(\mathfrak c,\Theta_j));\Theta_j\bigr)-J^{\rm best}(\Theta_j)}{J^{\rm best}(\Theta_j)}.
\label{eq:methodgap}
\end{equation}

A low average gap can conceal poor performance on individual instances. We therefore aggregate fitted gaps using a common mean--CVaR measure \citep{rockafellar2000optimization, bayraksan2015datadriven}. For a nonempty vector $\mathbf x=(x_1,\ldots,x_n)$, define the aggregation function as 
\begin{equation}
\rho_{\alpha,\lambda}(\mathbf x)
=(1-\lambda)\operatorname{Mean}(\mathbf x)
+\lambda\operatorname{CVaR}_{\alpha}(\mathbf x),
\qquad 0\le\alpha<1,\quad 0\le\lambda\le1,
\label{eq:methodmeancvar}
\end{equation}
where
$
\operatorname{CVaR}_{\alpha}(\mathbf x)
=\min_{u\in\mathbb R}
\left\{u+\frac{1}{n(1-\alpha)}
\sum_{i=1}^n[x_i-u]^+\right\}
$
measures the average loss in the worst $1-\alpha$ fraction of probability mass, and $\lambda$ controls its weight relative to the mean. 

For any evaluation set $\mathcal T$, define its configuration set and corresponding instance subsets by
\begin{equation}
\mathfrak S(\mathcal T)=\{\bm\sigma:\exists\,\mathcal D,\bm c\text{ such that }(\bm\sigma,\mathcal D,\bm c)\in\mathcal T\},
\qquad
\mathcal T_{\bm\sigma}=\{\Theta\in\mathcal T:\Theta=(\bm\sigma,\mathcal D,\bm c)\},
\quad \bm\sigma\in\mathfrak S(\mathcal T).
\label{eq:methodscorepartition}
\end{equation}
The scoring function $f(\mathfrak c;\mathcal T)$ takes the complete candidate and the evaluation set as its arguments. It aggregates fitted gaps within each configuration, then aggregates the resulting configuration scores:
\begin{equation}
R_{\bm\sigma}(\mathfrak c;\mathcal T)
=\rho_{\alpha,\lambda}\!\left(
\bigl(\Delta(\mathfrak c;\Theta)\bigr)_{\Theta\in\mathcal T_{\bm\sigma}}\right),
\qquad
f(\mathfrak c;\mathcal T)
=\rho_{\alpha,\lambda}\!\left(
\bigl(R_{\bm\sigma}(\mathfrak c;\mathcal T)\bigr)_{\bm\sigma\in\mathfrak S(\mathcal T)}\right).
\label{eq:methodsearchscore}
\end{equation}
We use $(\alpha,\lambda)=(0.80,0.40)$ at both aggregation levels across all experimental settings, during discovery and benchmark evaluation. For a fixed search set $\mathcal T$, we abbreviate $f(\mathfrak c;\mathcal T)$ as $f(\mathfrak c)$. 

\medskip
\noindent\underline{\textbf{Trade-offs between performance and simplicity.}}
$f(\mathfrak c)$ measures cost performance across the evaluation instances, with lower scores indicating better performance. We also need to account for the simplicity of its symbolic structure when selecting candidates. 
We assess simplicity in two ways: $\zeta_{\rm par}(\mathfrak c)$ counts independently fitted parameter symbols, and $\zeta_{\rm cmp}(\mathfrak c)$ measures structural complexity after reusable subexpressions have been expanded (Section~\ref{app:discoverycomplexity}). For $k\in\{\mathrm{par},\mathrm{cmp}\}$, write $\mathfrak c'\prec_k\mathfrak c$ if $\mathfrak c'$ has no higher search score and no greater $\zeta_k$, with a strict improvement in at least one. {\color{black}Let $\mathcal G_r^{\rm pop}\subseteq\mathcal H$ denote the population of candidates evaluated by the start of iteration $r$.} The Pareto frontier is
\begin{equation}
 \mathcal F^k({\color{black}\mathcal G_r^{\rm pop}})=\{\mathfrak c\in{\color{black}\mathcal G_r^{\rm pop}}:\nexists \mathfrak c'\in{\color{black}\mathcal G_r^{\rm pop}}\text{ with }\mathfrak c'\prec_k \mathfrak c\},\qquad k\in\{\mathrm{par},\mathrm{cmp}\}.
 \label{eq:methodbilevel}
\end{equation}
It gives empirical trade-offs among evaluated structures.

\subsection{Feedback-guided LLM search with reusable subexpressions}
\label{sec:methodclosedloop}
\begin{figure}[!htbp]
\centering
\includegraphics[width=\textwidth]{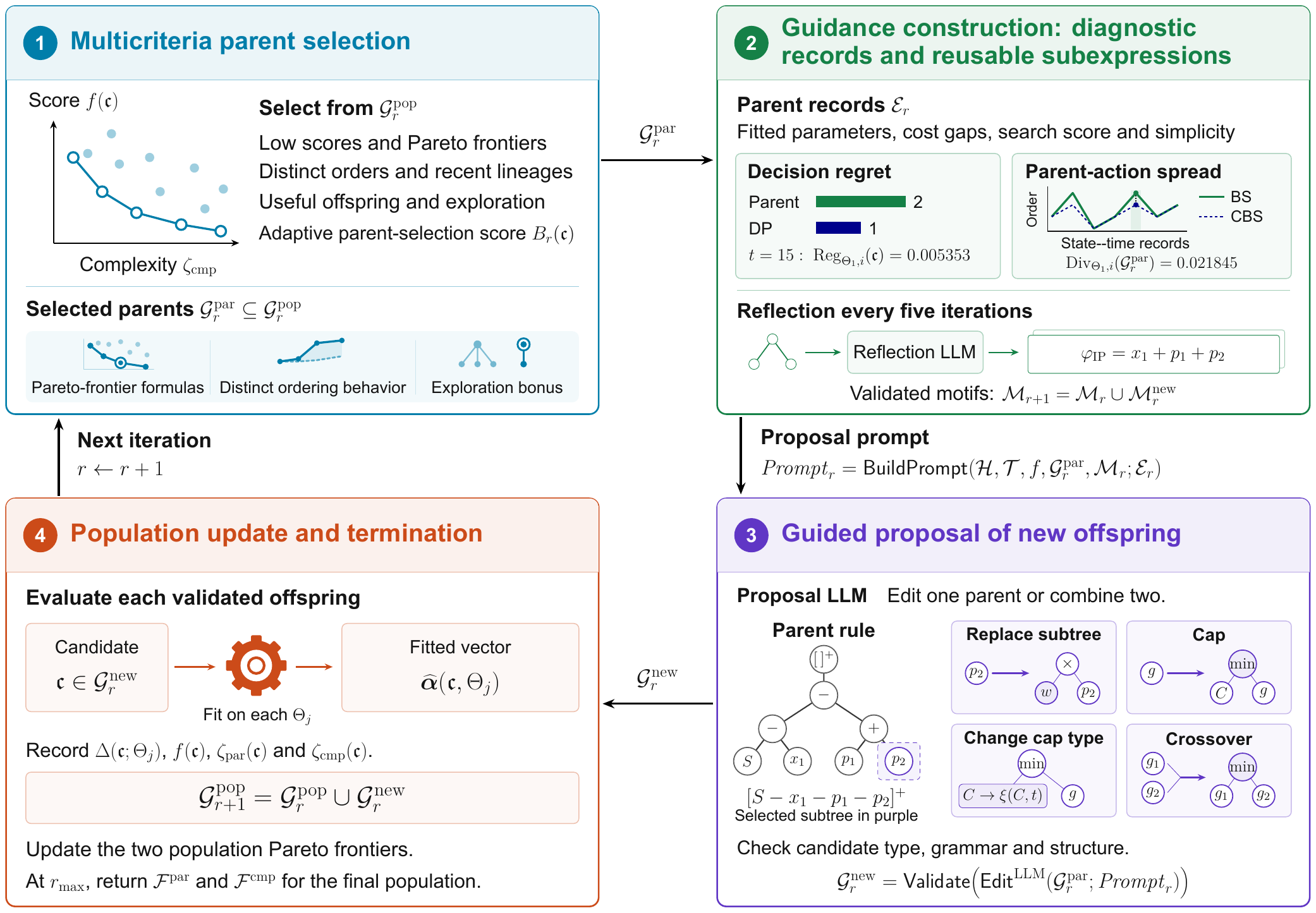}
\caption{Feedback-guided LLM search.}
\label{fig:discoveryloop}
\end{figure}

Our feedback-guided LLM search runs in iterations. For fixed search inputs $\mathcal H$, $\mathcal T$, $f$, and $\mathcal G_0^{\rm pop}$, $r$ indexes the inner iterations of structural search. The search begins by fitting the parameters of the initial candidates in $\mathcal G_0^{\rm pop}$ and evaluating them on $\mathcal T$. Each iteration $r$ comprises four steps, as illustrated in Figure~\ref{fig:discoveryloop}: (1) select a parent set $\mathcal G_r^{\rm par}\subseteq\mathcal G_r^{\rm pop}$ using the multicriteria procedure; (2) construct diagnostic records and extract reusable symbolic subexpressions to guide structural edits; (3) use this guidance to generate offspring from the selected parents; and (4) fit and evaluate the offspring and add them to the population. The search terminates when it reaches the predefined iteration limit. In the following paragraphs, we explain each step in more detail.

\medskip
\noindent\underline{\textbf{Step 1: multicriteria parent selection.}} The selection mechanisms in AIPS and InvEvolve primarily emphasize each policy's evaluated performance~\citep{yang2026aips,huang2026invevolve}. In contrast, we construct $\mathcal G_r^{\rm par}$ through a multicriteria procedure that balances performance and simplicity,
retains parents with distinct ordering behavior, and uses offspring quality to identify productive structures for further search. Joint consideration of these parents allows the agent to compare contrasting replenishment responses and recombine symbolic components across formulas, broadening the structural alternatives available for subsequent proposals. Specifically, we reserve places for formulas with the best current cost scores and for Pareto-frontier formulas that balance inventory cost against policy complexity. We select formulas with distinct ordering behavior to explore alternative replenishment responses to the same inventory conditions, and formulas from recent lineages to give newly emerging structures opportunities for further refinement. For the remaining places, we sample candidates $\mathfrak c\in\mathcal G_r^{\rm pop}$ using an adaptive parent-selection score $B_r(\mathfrak c)$. This score rewards a parent when its offspring improves on the best score among its contributing parents or joins a Pareto frontier of the candidates evaluated so far, favoring parents that generate useful structural improvements. An exploration bonus gives candidates with few evaluated offspring an opportunity to demonstrate their value as parents. Section~\ref{app:discoverysearchcontrol} provides the parent-selection details.

\medskip
\phantomsection
\noindent\underline{\textbf{Step 2: guidance construction: diagnostic records and reusable subexpressions.}}
We construct \emph{state-specific operational diagnostics} that highlight costly ordering decisions and contrasting parent actions under the same inventory conditions, giving the LLM more targeted guidance for structural revision than aggregate performance scores alone. Specifically, for each instance $\Theta_j\in\mathcal T$, we compare the order prescribed by each parent $\mathfrak c\in\mathcal G_r^{\rm par}$ with an optimal order, obtained through dynamic programming (DP). We make this comparison at selected periods $t_i$ and states $s_i$ reachable under the DP-optimal policy. The weight $\omega_i$ is proportional to the probability of reaching $s_i$ at $t_i$ under that policy. The state-specific \emph{decision regret} $\operatorname{Reg}_{\Theta_j,i}(\mathfrak c)$ is $\omega_i$ times the additional expected cost of placing the parent's order at $(t_i,s_i)$ and acting optimally thereafter. For each parent, we record the periods and states with the largest regrets, together with the corresponding orders and cost differences. These records highlight costly ordering decisions and identify opportunities for structural revision.  We also record the states and periods where the selected parents disagree most. We rank them by state-specific \emph{parent-action spread} $\operatorname{Div}_{\Theta_j,i}(\mathcal G_r^{\rm par})$, the difference between the largest and smallest orders the parents prescribe at $(t_i,s_i)$, weighted by $\omega_i$. These records help the LLM identify and edit the parts of the parent formulas that produce different orders under the same inventory conditions. We denote the evaluation and diagnostic records selected for iteration $r$ by $\mathcal E_r$ and use them to guide the proposal of new offspring. The evaluation records include the fitted parameters, cost gaps, search scores, parameter counts, and structural complexity measures obtained in Step~4. Section~\ref{app:discoverydiagnostics} provides the diagnostic
definitions and selection details.

\begin{figure}[!htbp]
\centering
\includegraphics[width=\textwidth]{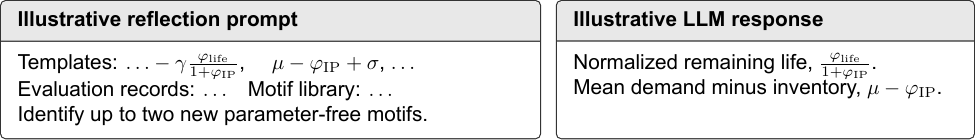}
\caption{Extracting reusable motifs from perishable-inventory templates.}
\label{fig:reflectionprompt}
\end{figure}

Another way to guide the LLM search is to build a library of useful subexpressions (called motifs) from formulas evaluated in earlier iterations. A reflection agent identifies these motifs for reuse in later proposals. Every five iterations, it reviews {\color{black}the parents $\mathcal G_r^{\rm par}$ selected for that iteration}. Its prompt contains their expressions, their evaluation records ${\color{black}\mathcal E_r}$, the current motif library, and the grammar and structural constraints on motifs. It asks for at most two new motifs, each with a name and a parameter-free AST. At a reflection iteration, the accepted motifs are generated by
\begin{equation}
\mathcal M_r^{\rm new}
=\mathsf{Validate}\!\left(\mathsf{Reflect}^{\rm LLM}
  ({\color{black}\mathcal G_r^{\rm par}},{\color{black}\mathcal E_r},\mathcal M_r;\mathcal H)\right).
\label{eq:methodmotifupdate}
\end{equation}
Here, $\mathsf{Validate}$ retains only motifs passing the candidate type, grammar, and structural checks. The motif library is updated by $\mathcal M_{r+1}=\mathcal M_r\cup\mathcal M_r^{\rm new}$ and remains unchanged between reflection iterations. Providing accepted motifs in subsequent proposal prompts lets the LLM reuse a validated subexpression by name in a new candidate, saving it the work of reconstructing the same expression from primitive operators. Section~\ref{app:discoveryprompts} describes the reflection prompt inputs and response formats. Figure~\ref{fig:reflectionprompt} gives an example.

\begin{figure}[!htbp]
\centering
\includegraphics[width=\textwidth]{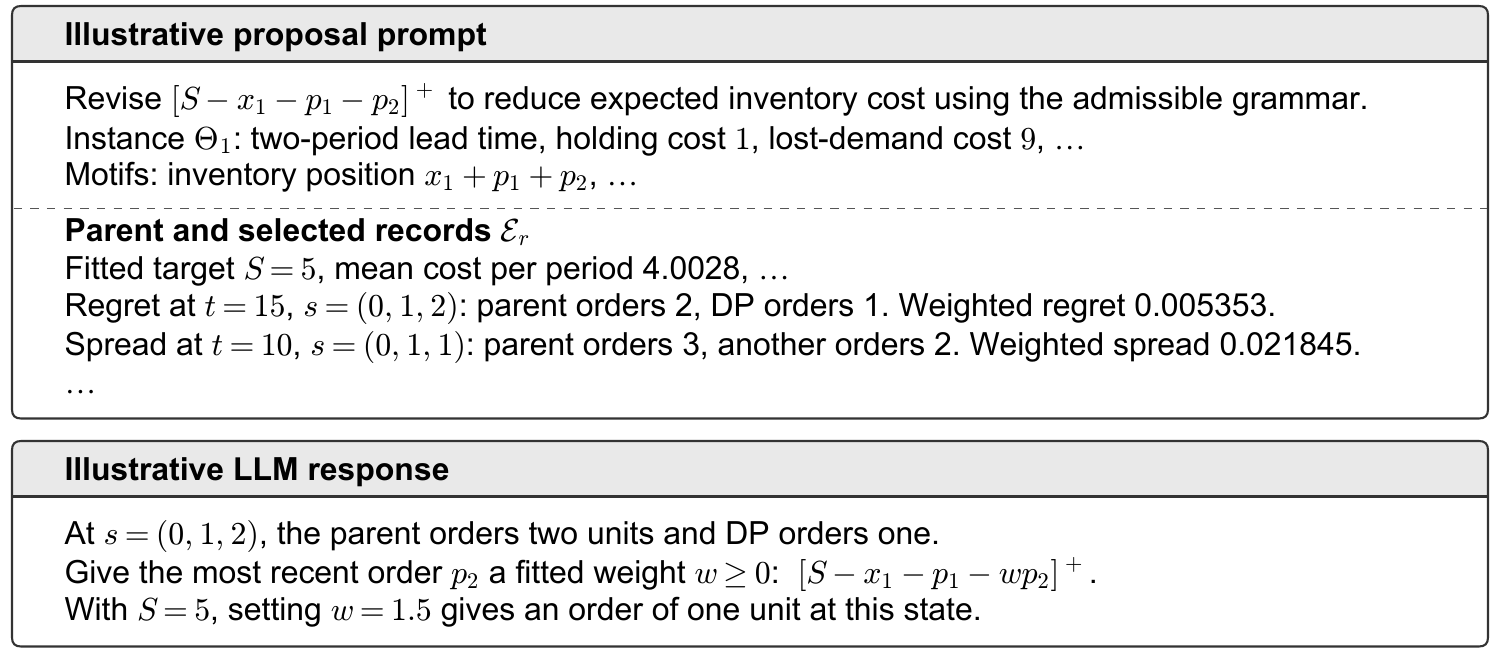}
\caption{Illustrative proposal prompt with computed records and a possible LLM response.}
\label{fig:proposalprompt}
\end{figure}

\medskip
\noindent\underline{\textbf{Step 3: guided proposal of new offspring.}}
A proposal LLM agent generates offspring by editing the selected parent formulas. The agent receives the selected parents $\mathcal G_r^{\rm par}$, their evaluation and diagnostic records $\mathcal E_r$, the grammar defining $\mathcal H$, and the motif library $\mathcal M_r$. For SOF candidates, edits modify the ordering expression; for SOF-template candidates, they can modify both the ordering expression and the state-summary definitions. Figure~\ref{fig:proposalprompt} illustrates a proposal prompt and a possible symbolic edit for a parent ordering rule.

With these inputs specified, we construct $\mathit{Prompt}_r=\mathsf{BuildPrompt}(\mathcal H,\mathcal T,f,\mathcal G_r^{\rm par},\mathcal M_r;\mathcal E_r)$. The proposal LLM agent uses this prompt to propose symbolic modifications to the selected parents. Applying these modifications and validating the resulting offspring candidates gives
\begin{equation}
\mathcal G_r^{\rm new}
=\mathsf{Validate}\!\left(
\mathsf{Edit}^{\rm LLM}\!\left(\mathcal G_r^{\rm par};
\mathit{Prompt}_r\right)\right)
\subseteq\mathcal H.
\label{eq:methodpromptgeneration}
\end{equation}
Here, $\mathsf{Edit}^{\rm LLM}$ obtains symbolic modifications from the LLM and constructs each offspring from one or two designated parents. The operator $\mathsf{Validate}$ retains the resulting candidates passing the type, grammar, and structural checks. The permitted edits comprise subtree
replacement, capping, crossover, changes to the cap type, and, for templates, summary replacement. detailed in Section~\ref{app:discoveryedits} details these actions, while Section~\ref{app:discoveryprompts} specifies the prompt inputs and response formats. As an illustration of the LLM edits, in Figure~\ref{fig:proposalprompt}, the parent rule is $[S-x_1-p_1-p_2]^+$. Introducing a fitted weight for $p_2$ yields $[S-x_1-p_1-wp_2]^+$, where $w\ge0$ is a fitted weight for the most recently placed pipeline order. For a capped rule, the LLM can also replace the original shortfall with the smaller of the inventory-position and far-pipeline shortfalls, yielding $\min\{\xi(C,t),[S-x_1-p_1-p_2]^+,[F-p_2]^+\}$, where $F\ge0$ is a fitted far-pipeline target. 

\medskip
\noindent\underline{\textbf{Step 4: {\color{black}population update} and termination.}}
For each candidate $\mathfrak c\in\mathcal G_r^{\rm new}$, we fit and evaluate them using the process in Section~\ref{sec:methodoverview}. This generates the evaluation record including fitted parameters $\widehat{\bm\alpha}(\mathfrak c,\Theta_j)$ and cost gap $\Delta(\mathfrak c;\Theta_j)$ for each $\Theta_j\in\mathcal T$, along with the {\color{black}candidate's} search score $f(\mathfrak c)$, parameter count $\zeta_{\rm par}(\mathfrak c)$, and structural complexity $\zeta_{\rm cmp}(\mathfrak c)$. 
{\color{black}After evaluation, we add the offspring to the population, giving $\mathcal G_{r+1}^{\rm pop}=\mathcal G_r^{\rm pop}\cup\mathcal G_r^{\rm new}$.} Their scores and complexity measures update the frontiers in Eq.~\eqref{eq:methodbilevel}. The run stops after a prespecified number of iterations $r_{\max}$ (100 by default) and returns the two frontiers of the final {\color{black}population}. A final formula can be chosen according to the desired trade-off between cost and simplicity.

\subsection{Progressive expansion across instances and configurations}
\label{sec:methodstages}
Evaluating every candidate across many configurations is expensive because each instance requires its own parameter fit. We therefore begin with a single instance and progressively expand the evaluation set $\mathcal T$, using the structures selected in one stage to initialize the next. The three stages form an outer progression: each executes the same inner search iterations in Section~\ref{sec:methodclosedloop}, while changing the candidate class, evaluation set, and initialization. All stages use the same scoring function $f$ in Eq.~\eqref{eq:methodsearchscore}. Let $\mathsf{Search}$ denote a complete run of the inner search in Section~\ref{sec:methodclosedloop}. For stage $\ell\in\{\mathrm I,\mathrm{II},\mathrm{III}\}$, we write
\begin{equation}
(\mathcal F_{\ell}^{\rm par},\mathcal F_{\ell}^{\rm cmp})
=\mathsf{Search}(\mathcal H_{\ell},\mathcal T_{\ell},f,\mathcal G_{\ell,0}^{\rm pop}).
\label{eq:methodstagesearch}
\end{equation}
Here, the candidate class $\mathcal H_{\ell}$, evaluation set $\mathcal T_{\ell}$, and initial population $\mathcal G_{\ell,0}^{\rm pop}$ are specified for stage $\ell$, while the scoring function $f$ is shared across stages and applied to the stage-specific arguments $(\mathfrak c,\mathcal T_\ell)$. These are the four inputs defined in Section~\ref{sec:methodoverview}; the outputs are the stage's two Pareto frontiers.

\begin{figure}[!htbp]
\centering
\includegraphics[width=\textwidth]{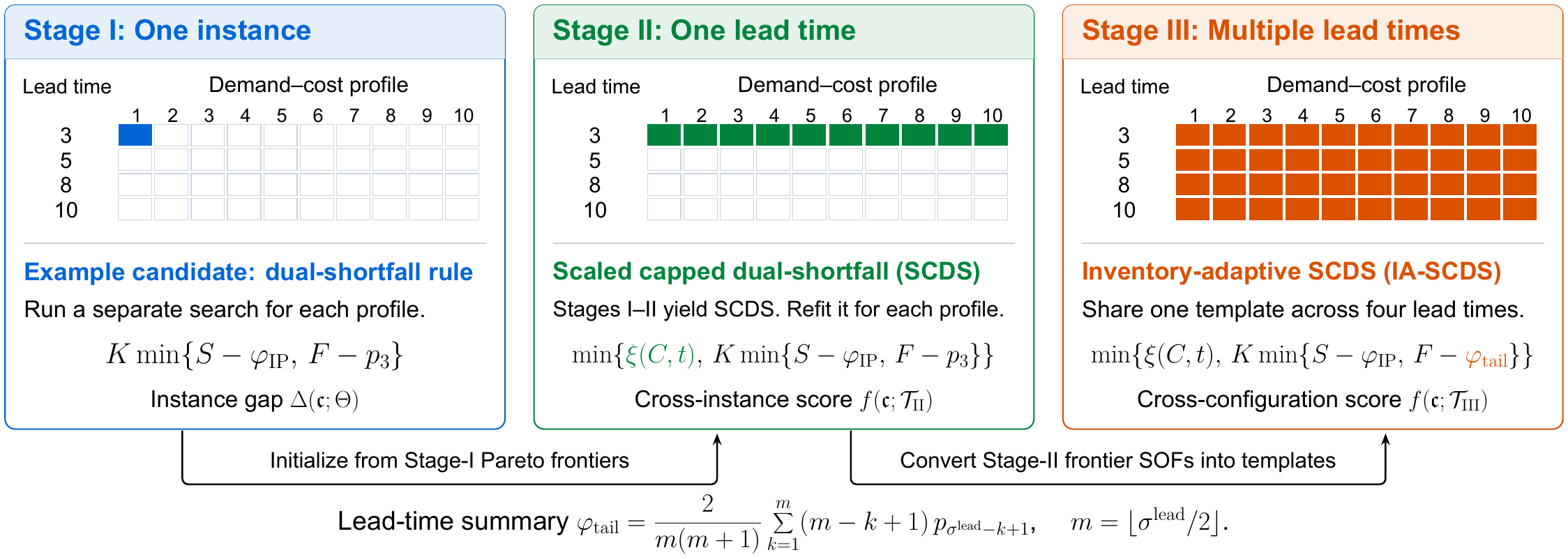}
\caption{Progressive expansion across the three search stages.}
\label{fig:threestages}
\end{figure}

Stage I searches for an SOF on a single instance at a fixed configuration, without imposing cross-instance sharing. Stage II shares an SOF across demand distributions and cost vectors at a fixed configuration. Stage III shares an SOF template that adapts its state summaries to different configurations. Figure~\ref{fig:threestages} illustrates this progression for the lead-time experiment.

\medskip
\noindent\underline{\textbf{Stage I: search for SOFs on one instance.}}
For a selected instance $\Theta=(\bm\sigma,\mathcal D,\bm c)$, set $\mathcal H_{\rm I}=\mathcal H_{\bm\sigma}^{\rm SOF}$ and $\mathcal T_{\rm I}=\{\Theta\}$. The initial population contains admissible SOFs, and the common score reduces to $f(\mathfrak c;\{\Theta\})=\Delta(\mathfrak c;\Theta)$ because both aggregation levels have a single entry. With only one instance to fit per candidate, this stage supports rapid structural exploration under the same search budget used in the later stages. It returns the frontiers $\mathcal F_{\rm I}^{\rm par}$ and $\mathcal F_{\rm I}^{\rm cmp}$ for that instance.

\medskip
\noindent\underline{\textbf{Stage II: search for SOFs shared across instances at one configuration.}}
Keep $\mathcal H_{\rm II}=\mathcal H_{\bm\sigma}^{\rm SOF}$ and expand $\mathcal T_{\rm II}$ to $m\ge2$ instances with different demand distributions or cost vectors at the same configuration $\bm\sigma$. The initial population is $\mathcal G_{{\rm II},0}^{\rm pop} =\mathcal F_{\rm I}^{\rm par}\cup\mathcal F_{\rm I}^{\rm cmp}$. Every inherited candidate has its parameters refitted on every instance before it receives a Stage-II score. The search then modifies these structures to find formulas that perform well across the expanded instance set.

Because every instance in Stage II has configuration $\bm\sigma$, the outer aggregation in Eq.~\eqref{eq:methodsearchscore} has a single entry. The score reduces to
\begin{equation}
f(\mathfrak c;\mathcal T_{\rm II})
=R_{\bm\sigma}(\mathfrak c;\mathcal T_{\rm II})
=\rho_{\alpha,\lambda}\!\left(
\bigl(\Delta(\mathfrak c;\Theta)\bigr)_{\Theta\in\mathcal T_{\rm II}}\right).
\label{eq:methodstage2score}
\end{equation}

{\color{black}Let $\mathcal G_{{\rm II},r_{\max}}^{\rm pop}$ denote the final Stage-II population.} The stage returns the two Pareto frontiers $\mathcal F_{\rm II}^{k}=\mathcal F^k({\color{black}\mathcal G_{{\rm II},r_{\max}}^{\rm pop}})$, $k\in\{\mathrm{par},\mathrm{cmp}\}$, using $\mathcal F^k$ in Eq.~\eqref{eq:methodbilevel}.

\medskip
\noindent\underline{\textbf{Stage III: search for SOF templates across configurations.}}
Stage III expands evaluation across a finite configuration set $\mathfrak S$, using $\mathcal T_{\rm III}=\bigcup_{\bm\sigma\in\mathfrak S}\mathcal T_{\bm\sigma}$ with a nonempty instance set at each configuration. In the lead-time setting, Stage II varies demand and costs at a fixed lead time. Stage III also varies the lead time and searches for a single ordering template that applies across all these instances, so $\mathcal H_{\rm III}=\mathcal H_{\mathfrak S}^{\rm template}$. 

To initialize this stage, the LLM converts SOFs in $\mathcal F_{\rm II}^{\rm par}\cup\mathcal F_{\rm II}^{\rm cmp}$ into templates by replacing fixed-coordinate state terms with configuration-dependent summaries. The initialization prompt lists the Stage-II frontier formulas, explains their state coordinates, and specifies the target configurations and template grammar. It asks the LLM to return each proposed template as a triple $(g,\Phi,\bm\alpha)$. Writing this prompt as $\mathit{Prompt}_{\rm III}^{\rm init}$, the initial population is
\begin{equation}
\mathcal G_{{\rm III},0}^{\rm pop}
=\mathsf{Validate}\!\left(
\mathsf{Generalize}^{\rm LLM}\!\left(
\mathcal F_{\rm II}^{\rm par}\cup\mathcal F_{\rm II}^{\rm cmp};
\mathit{Prompt}_{\rm III}^{\rm init}\right)\right)
\subseteq\mathcal H_{\rm III}.
\label{eq:methodstage3initialization}
\end{equation}
Here, $\mathsf{Generalize}^{\rm LLM}$ denotes the LLM call that constructs SOF templates from the retained SOFs, and $\mathsf{Validate}$ retains admissible templates using the checks described in Section~\ref{app:discoveryedits}. The initialization prompt and its required output are described in Section~\ref{app:discoveryprompts}.

For example, suppose a Stage-II frontier contains an SOF pair $(g_0,(C,S,w))$ at $\bm\sigma=(2,\infty,0)$, where $g_0(s;C,S,w)=\min\{\xi(C,t),[S-x_1-p_1-wp_2]^+\}$ and $w\ge0$. For target configurations $\bm\sigma=(\sigma^{\rm lead},\infty,0)$ with $\sigma^{\rm lead}\ge1$, applying $\mathsf{Generalize}^{\rm LLM}$ to this pair may produce a template $(g,\{\varphi_w\},(C,S,w))$ with ordering expression $g=\min\{\xi(C,t),[S-\varphi_w]^+\}$ and summary definition $\varphi_w(s;\bm\sigma,w)=x_1+\sum_{k=1}^{\sigma^{\rm lead}-1}p_k+w p_{\sigma^{\rm lead}}$. Instantiating this template at $\sigma^{\rm lead}=2$ recovers $g_0$; at $\sigma^{\rm lead}=3$, it yields the SOF pair $(g_{(3,\infty,0)},(C,S,w))$ with $g_{(3,\infty,0)}(s;C,S,w)=\min\{\xi(C,t),[S-x_1-p_1-p_2-wp_3]^+\}$. The search can modify both the ordering expression and the summary definitions, while each template must yield a valid SOF throughout $\mathfrak S$.

Stage III uses the full two-level score in Eq.~\eqref{eq:methodsearchscore}. Writing $\mathcal T_{\rm III}$ for its evaluation set, we obtain
\begin{equation}
f(\mathfrak c;\mathcal T_{\rm III})
=\rho_{\alpha,\lambda}\!\left(
\bigl(R_{\bm\sigma}(\mathfrak c;\mathcal T_{\rm III})\bigr)_{\bm\sigma\in\mathfrak S(\mathcal T_{\rm III})}\right),
\label{eq:methodstage3cross}
\end{equation}
Each configuration score is computed from its fitted instance gaps using the same risk measure. Stage III returns the parameter-count and structural-complexity Pareto frontiers using this score.

\subsection{The complete algorithm}
\label{sec:methodframework}
We embed the LLM-guided search in the three-stage expansion framework to obtain the complete symbolic policy discovery algorithm, presented in Algorithm~\ref{alg:discoveryframework}. The outer loop, indexed by $\ell$, progressively expands the search from a single instance to multiple instances at a fixed configuration and then to multiple configurations. Within each stage, $\mathcal H$, $\mathcal T$, and $f$ remain fixed while the inner loop, indexed by $r$, updates the candidate population. Parent records and motifs guide $\mathsf{Edit}^{\rm LLM}$, while numerical evaluation provides the scores and diagnostics for subsequent iterations. In the algorithm, Stage-I frontiers initialize Stage II directly; Stage-II frontiers are converted to SOF templates by $\mathsf{Generalize}^{\rm LLM}$ before Stage III. Transferred structures are re-evaluated under the new stage's inputs. Each inner search stops after $r_{\max}$ iterations, and the complete procedure returns $(\mathcal F_{\rm III}^{\rm par},\mathcal F_{\rm III}^{\rm cmp})$, preserving performance--simplicity trade-offs for final selection.

\begin{algorithm}[htb]
\caption{The complete agentic inventory-policy discovery algorithm}
\label{alg:discoveryframework}
\small
\renewcommand{\baselinestretch}{1.02}\selectfont
\begin{algorithmic}[1]
\Require Common scoring function $f$ (Eq.~\eqref{eq:methodsearchscore}); $(\mathcal H_\ell,\mathcal T_\ell)$ for $\ell\in\{\mathrm I,\mathrm{II},\mathrm{III}\}$ (Section~\ref{sec:methodstages})
\Require $\mathcal G_{{\rm I},0}^{\rm pop}$; initial motif libraries $\mathcal M_{\ell,0}$; $r_{\max}$; $\mathit{Prompt}_{\rm III}^{\rm init}$
\Ensure $(\mathcal F_{\rm III}^{\rm par},\mathcal F_{\rm III}^{\rm cmp})$
\For{$\ell\in\{\mathrm I,\mathrm{II},\mathrm{III}\}$} \Comment{outer loop}
  \If{$\ell=\mathrm{II}$}
    \State $\mathcal G_{{\rm II},0}^{\rm pop}\gets\mathcal F_{\rm I}^{\rm par}\cup\mathcal F_{\rm I}^{\rm cmp}$
  \ElsIf{$\ell=\mathrm{III}$}
    \State $\mathcal G_{{\rm III},0}^{\rm pop}\gets\mathsf{Validate}\!\left(\mathsf{Generalize}^{\rm LLM}\!\left(\mathcal F_{\rm II}^{\rm par}\cup\mathcal F_{\rm II}^{\rm cmp};\mathit{Prompt}_{\rm III}^{\rm init}\right)\right)$
  \EndIf
  \State $\mathcal H\gets\mathcal H_\ell$, $\mathcal T\gets\mathcal T_\ell$, $\mathcal G_0^{\rm pop}\gets\mathcal G_{\ell,0}^{\rm pop}$, $\mathcal M_0\gets\mathcal M_{\ell,0}$
  \State Evaluate $\mathcal G_0^{\rm pop}$ using Eqs.~\eqref{eq:methodgap}--\eqref{eq:methodsearchscore}
  \For{$r=0,\ldots,r_{\max}-1$} \Comment{inner loop}
    \State Select $\mathcal G_r^{\rm par}\subseteq\mathcal G_r^{\rm pop}$ (Section~\ref{app:discoverysearchcontrol})
    \State Collect $\mathcal E_r$ for $\mathcal G_r^{\rm par}$ (Section~\ref{app:discoverydiagnostics})
    \State $\mathit{Prompt}_r\gets\mathsf{BuildPrompt}(\mathcal H,\mathcal T,f,\mathcal G_r^{\rm par},\mathcal M_r;\mathcal E_r)$
    \State $\mathcal G_r^{\rm new}\gets\mathsf{Validate}\!\left(\mathsf{Edit}^{\rm LLM}(\mathcal G_r^{\rm par};\mathit{Prompt}_r)\right)$
    \For{$\mathfrak c\in\mathcal G_r^{\rm new}$}
      \State $\forall\Theta_j\in\mathcal T:$ Fit $\widehat{\bm\alpha}(\mathfrak c,\Theta_j)$ and compute $\Delta(\mathfrak c;\Theta_j)$ \Comment{Eq.~\eqref{eq:methodgap}}
      \State $f(\mathfrak c)\gets f(\mathfrak c;\mathcal T)$ using Eq.~\eqref{eq:methodsearchscore}
      \State Compute $(\zeta_{\rm par}(\mathfrak c),\zeta_{\rm cmp}(\mathfrak c))$ (Section~\ref{sec:methodoverview})
    \EndFor
    \State ${\color{black}\mathcal G_{r+1}^{\rm pop}\gets\mathcal G_r^{\rm pop}\cup\mathcal G_r^{\rm new}}$
    \State {\color{black}Update the frontiers $\mathcal F^k(\mathcal G_{r+1}^{\rm pop})$}, $k\in\{\mathrm{par},\mathrm{cmp}\}$ (Eq.~\eqref{eq:methodbilevel})
    \State $\mathcal M_{r+1}\gets\mathcal M_r$
    \If{reflection is scheduled} \Comment{Extract reusable motifs}
      \State $\mathcal M_r^{\rm new}\gets\mathsf{Validate}\!\left(\mathsf{Reflect}^{\rm LLM}({\color{black}\mathcal G_r^{\rm par}},{\color{black}\mathcal E_r},\mathcal M_r;\mathcal H)\right)$
      \State $\mathcal M_{r+1}\gets\mathcal M_r\cup\mathcal M_r^{\rm new}$
    \EndIf
  \EndFor
  \State $\mathcal F_\ell^k\gets\mathcal F^k({\color{black}\mathcal G_{r_{\max}}^{\rm pop}})$, $k\in\{\mathrm{par},\mathrm{cmp}\}$
\EndFor
\State Refit selected $\mathfrak c\in\mathcal F_{\rm III}^{\rm par}\cup\mathcal F_{\rm III}^{\rm cmp}$ with a larger budget
\State \Return $(\mathcal F_{\rm III}^{\rm par},\mathcal F_{\rm III}^{\rm cmp})$
\end{algorithmic}
\end{algorithm}

\section{Numerical Results}\label{sec:results}
\label{sec:methodbenchmark}
\medskip

\noindent\underline{\textbf{Experiment design.}}
We run the discovery framework separately for four inventory settings. Each run follows the three-stage procedure in Section~\ref{sec:methodstages}. The training configurations and demand--cost profiles are specified in Section~\ref{app:searchinstances}. All four settings use a common catalog $\mathfrak D^{\rm bench}$ of 50 demand distributions, combined with setting-specific configurations and cost vectors; these benchmark sets are specified together in Section~\ref{app:demandcatalog}. The compared literature baselines are described within each setting. Each setting also includes an independently generated SOF baseline from \emph{Codex (gpt-6-astra)}, with all formulas given in Section~\ref{app:codexbaseline}.  For each setting, let $\mathfrak S^{\rm bench}$ denote its set of benchmark configurations and $\mathcal C_{\bm\sigma}^{\rm bench}$ its set of cost vectors at configuration $\bm\sigma$. Using the instance definition $\Theta=(\bm\sigma,\mathcal D,\bm c)$, the complete benchmark instance set and its configuration-specific subsets are
\begin{equation}
\mathcal I=\{(\bm\sigma,\mathcal D,\bm c):\bm\sigma\in\mathfrak S^{\rm bench},\ \mathcal D\in\mathfrak D^{\rm bench},\ \bm c\in\mathcal C_{\bm\sigma}^{\rm bench}\},
\qquad
\mathcal I_{\bm\sigma}=\{\bm\sigma\}\times\mathfrak D^{\rm bench}\times\mathcal C_{\bm\sigma}^{\rm bench}.
\label{eq:benchmarksets}
\end{equation}

\medskip
\noindent\underline{\textbf{Reporting metrics.}}
For benchmark reporting, the reference cost on each instance is the lowest measured cost among all evaluated policies. Benchmark gaps use Eq.~\eqref{eq:methodgap} with this empirical reference in place of the prespecified baseline reference used during discovery. For a selected candidate $\mathfrak c$, configuration-level fitness $R_{\bm\sigma}^{\rm bench}$ and aggregate fitness $f^{\rm bench}$ apply the common aggregations in Eq.~\eqref{eq:methodsearchscore} to these benchmark gaps on $\mathcal I$:
\begin{equation}
R_{\bm\sigma}^{\rm bench}(\mathfrak c)=R_{\bm\sigma}(\mathfrak c;\mathcal I),
\qquad
f^{\rm bench}(\mathfrak c)=f(\mathfrak c;\mathcal I).
\label{eq:benchmarkfitness}
\end{equation}
Configuration-level fitness summarizes mean and upper-tail relative gaps within a configuration, while aggregate fitness applies the same risk measure across configurations. All reported fitness values use a percentage scale; lower values indicate better performance.

The main text focuses on two inventory settings. Section~\ref{sec:resultsleadtime} presents detailed results for the lead-time setting, and Section~\ref{sec:resultsperishable} reports results for the perishable setting. Additional results for the joint lead-time--perishability and lead-time--backlog settings are presented in Section~\ref{app:additionalsettings}.

\subsection{Lead-time setting}
\label{sec:resultsleadtime}

We first consider nonperishable inventory with lost sales, so that
$\bm\sigma=(\sigma^{\rm lead},\infty,0)$.
We compare the discovered rules with ten baselines: seven SOF policies
and three non-SOF policies from the literature.
The SOF baselines comprise Constant Order \citep{goldberg2016constant},
Base Stock \citep{zipkin2008old}, Capped Base Stock
(CBS) \citep{xin2021capped}, Tilted PIC and Weighted Tilted PIC
\citep{huang2026invevolve}, AIPS Hybrid \citep{yang2026aips}, and the one generated by Codex.
For reference, we also include three non-SOF baselines:
Myopic-2 \citep{morton1971near,zipkin2008old},
Projected Inventory Level (PIL) \citep{vanjaarsveld2024projected},
and LLM-OR Projected MV \citep{baek2026llms}.

\begin{figure}[!htbp]
\centering
\includegraphics[width=\textwidth]{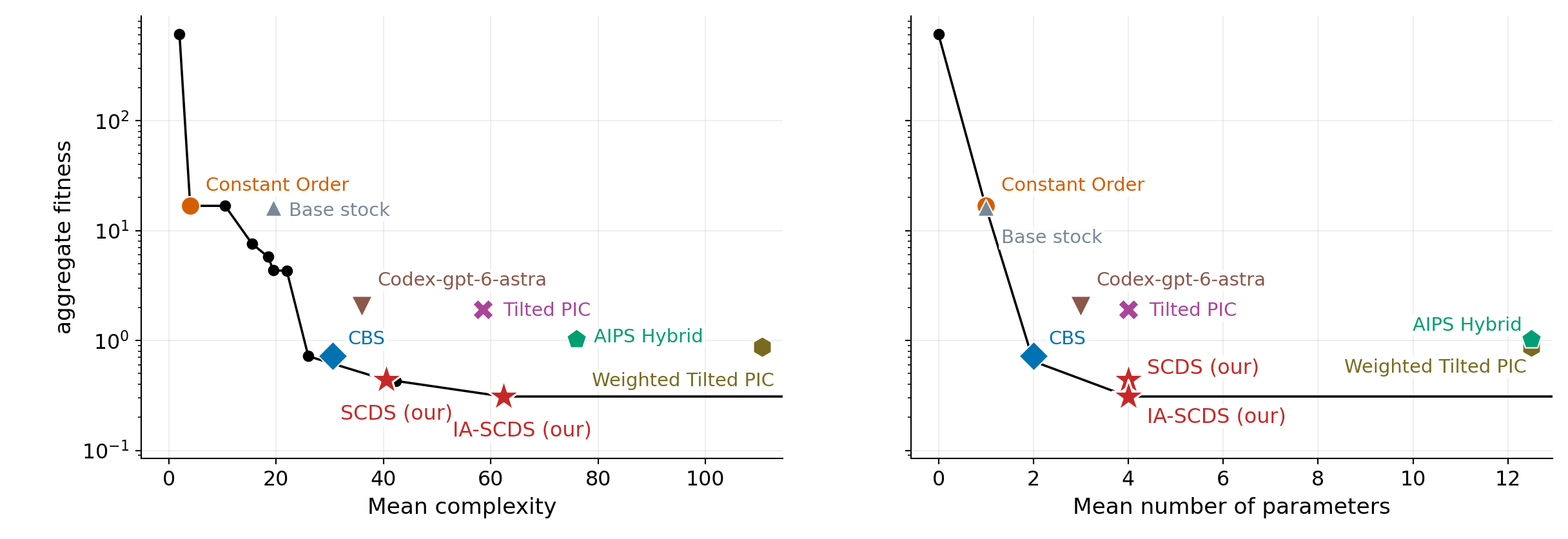}
\caption{Evaluation Pareto frontiers of the framework's final output candidates and comparison policies in the lead-time setting, with SCDS, IA-SCDS, and SOF baselines highlighted. Both panels report aggregate fitness $f^{\rm bench}$ and the horizontal axes report mean expanded AST complexity (left) and mean parameter count (right).}
\label{fig:resultspareto}
\label{fig:resultsparetoparam}
\label{fig:resultsparetocomplexity}
\end{figure}

\medskip
\noindent\underline{\textbf{Pareto frontier.}} Figure~\ref{fig:resultspareto} shows the evaluation Pareto frontiers alongside the SOF baselines, using aggregate fitness $f^{\rm bench}$ from Eq.~\eqref{eq:benchmarkfitness} over $\sigma^{\rm lead}=2,\ldots,15$. The horizontal axes measure expanded AST complexity and parameter count, averaged over the benchmark configurations. Constant Order, Base Stock, and CBS essentially coincide with frontier rules of the same form, whereas Tilted PIC, Weighted Tilted PIC, and AIPS Hybrid are dominated by discovered candidates with both lower fitness and smaller expressions. The complete frontier formulas are listed in Section~\ref{app:leadtimefrontiers}.

\medskip
\noindent\underline{\textbf{Representative policy templates.}} From this frontier, we select two representative templates:
\emph{scaled capped dual-shortfall} (SCDS) and
\emph{inventory-adaptive SCDS} (IA-SCDS).
SCDS is discovered in Stages~I--II with $\sigma^{\rm lead}=3$,
while IA-SCDS emerges in Stage~III.
Both formulas have four free parameters, $S$, $F$, $K$, and $C$:
\begin{align}
q_{\rm SCDS}
&=\min\!\left\{\xi(C,t),
K\min\{S-\varphi_{\rm IP},F-p_{\sigma^{\rm lead}}\}\right\},
\label{eq:resultscappeddual}\\
q_{\text{IA-SCDS}}
&=\min\!\left\{\xi(C,t),
K\min\{S-\varphi_{\rm IP},F-\varphi_{\rm tail}\}\right\}.
\label{eq:resultsimbalance}
\end{align}
Here, $\varphi_{\rm IP}=x_1+\sum_{r=1}^{\sigma^{\rm lead}}p_r$
is the inventory position, $p_{\sigma^{\rm lead}}$ is the
quantity scheduled to arrive latest in the pipeline, and
$
\varphi_{\rm tail}
=\frac{2}{m(m+1)}
\sum_{k=1}^{m}(m-k+1)p_{\sigma^{\rm lead}-k+1}$ (where $
m=\max\{1,\lfloor\sigma^{\rm lead}/2\rfloor\}$).
Both policies supplement the inventory-position shortfall with
an additional shortfall based on pipeline supply.
SCDS uses only the farthest pipeline coordinate, while IA-SCDS
uses the weighted summary $\varphi_{\rm tail}$, which adapts
to the lead time and assigns greater weight to later arrivals.
Each formula scales the smaller of the two shortfalls by $K$
and caps the result at $\xi(C,t)$.

\begin{figure}[!htbp]
\centering
\includegraphics[width=\textwidth]{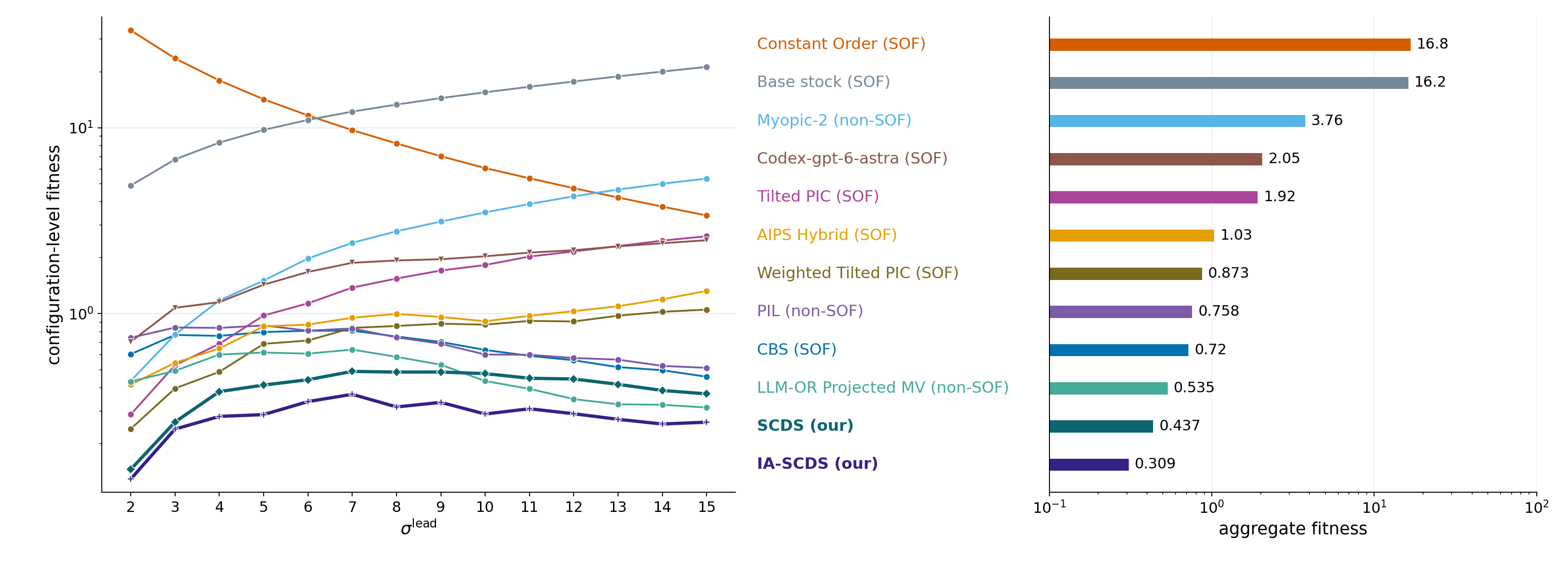}
\caption{Main benchmark comparison of SCDS, IA-SCDS, and the displayed baselines. The left panel reports configuration-level fitness for $\sigma^{\rm lead}=2,\ldots,15$; the right panel reports aggregate fitness across these lead times.}
\label{fig:resultsmcfitness}
\end{figure}

\medskip
\noindent\underline{\textbf{Evaluation of the representative templates.}} Figure~\ref{fig:resultsmcfitness} compares SCDS and IA-SCDS against the displayed baselines. As shown in the left panel, both policies achieve substantial improvements at short lead times, reducing configuration-level fitness by approximately $34\%$--$46\%$ relative to the best displayed baseline, Weighted Tilted PIC, for $\sigma^{\rm lead}\in\{2,3\}$. Their nearly identical performance at these lead times is consistent with the fact that $\varphi_{\rm tail}=p_{\sigma^{\rm lead}}$, making the two policies equivalent in functional form up to reparameterization. As lead time increases, SCDS continues to outperform all displayed baselines through $\sigma^{\rm lead}=9$, but is eventually surpassed by LLM-OR Projected MV. In contrast, IA-SCDS exhibits a growing advantage over SCDS at longer lead times, where its weighted summary of outstanding orders better captures the distribution of supply across the pipeline. It consistently achieves the lowest configuration-level fitness, improving upon the best displayed baseline by at least $16.5\%$ at every lead time. The right panel confirms these gains in aggregate: SCDS and IA-SCDS reduce aggregate fitness by $18.3\%$ and $42.3\%$, respectively, relative to the strongest displayed baseline overall.

\medskip
\noindent\underline{\textbf{Operational insights behind the representative templates.}}
Both SCDS and IA-SCDS contain CBS as a special case and thus retain
its theoretical guarantees under the corresponding parameter
choices~\citep{xin2021capped}.
While CBS responds to the aggregate inventory-position shortfall,
SCDS introduces a complementary constraint based on the preceding
order, limiting large orders in consecutive periods.
This constraint can reduce holding costs while increasing the
risk of lost sales.
Section~\ref{sec:theoretical-analysis} theoretically examines this trade-off in a model with a three-period lead time and establishes demand
and cost conditions under which the additional constraint improves
performance over CBS.
IA-SCDS further extends SCDS by replacing the single preceding
order $p_{\sigma^{\rm lead}}$ with a weighted summary
$\varphi_{\rm tail}$ of recent outstanding orders.
This summary reflects the timing of pipeline supply by assigning
greater weight to more recently placed orders.
IA-SCDS therefore limits new orders in response to a weighted
aggregate of recent orders, with both the aggregation window
and the weights adapting to the lead time.
This adaptation can improve performance in settings with longer
lead times.

\subsection{Perishable-inventory setting}
\label{sec:resultsperishable}
We next consider a pure-perishable system with zero lead time and lost sales, so that $\bm\sigma=(0,\sigma^{\rm life},0)$. The main comparison includes eight baselines: five SOF policies and three non-SOF policies from the literature. The SOF baselines are Base Stock, Capped Base Stock, BSP-low \citep{haijema2016stock}, Weighted Base Stock \citep{haijema2019improved}, and the one generated by Codex. For reference, we also include non-SOF baselines: Truncated Balancing (TB-ZAW) \citep{zhang2023truncated}, Proportional Balancing (PB-Chao), and Dual Balancing (DB-Chao) \citep{chao2015approximation}.

\medskip
\noindent\underline{\textbf{Representative policy templates.}} We select two representative templates from the framework's final output Pareto frontier: \emph{expiry-adjusted base stock} (EABS), with three parameters, and \emph{remaining-life capped shortfall} (RLCS), which achieves lower aggregate fitness with four parameters. The evaluation Pareto frontiers and complete frontier formulas are reported in Section~\ref{app:perishablefrontiers}. Write $\varphi_{\rm IP}=\sum_{m=1}^{\sigma^{\rm life}}x_m$ for carried inventory and $\varphi_{\rm life}=\sum_{m=1}^{\sigma^{\rm life}}m x_m$ for its remaining-life-weighted sum. EABS adjusts the base-stock shortfall using total inventory and the cohort expiring at the end of the current period, while RLCS uses a normalized remaining-life summary and caps the resulting order:
\begin{equation}
\begin{aligned}
q_{\rm EABS} = S-\omega_1\varphi_{\rm IP}+\omega_2 x_1,
\qquad
q_{\rm RLCS} = \min\!\left\{\xi(C,t),S-\omega_1\varphi_{\rm IP}-\gamma\frac{\varphi_{\rm life}}{1+\varphi_{\rm IP}}\right\},
\end{aligned}
\label{eq:resultsperishable}
\end{equation}
where $S, C, \omega_1,\omega_2, \gamma\ge0$ are fitted parameters. The normalized summary $\varphi_{\rm life}/(1+\varphi_{\rm IP})$ uses a one-unit regularization to remain well-defined at zero inventory.

\begin{figure}[!htbp]
\centering
\includegraphics[width=\textwidth]{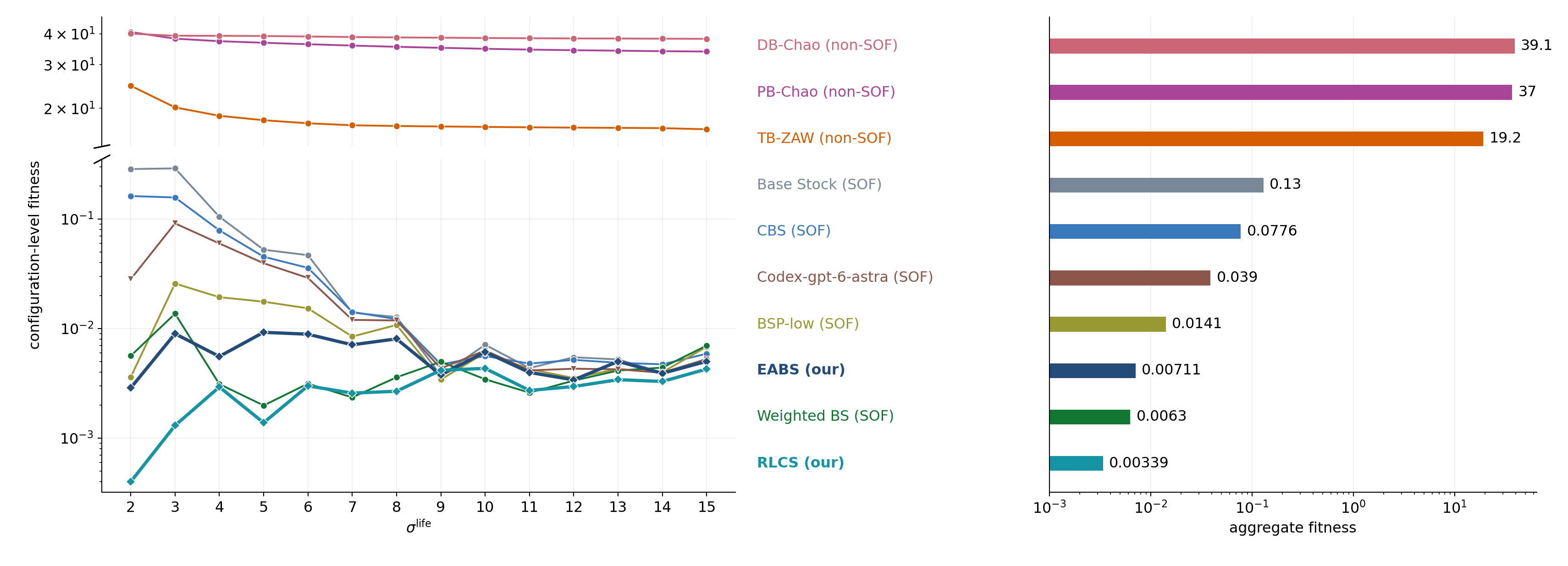}
\caption{Main benchmark comparison of EABS, RLCS, and the displayed baselines. The left panel reports configuration-level fitness for $\sigma^{\rm life}=2,\ldots,15$; the right panel reports aggregate fitness across these shelf lives.}
\label{fig:resultsperishablemc}
\end{figure}

\medskip
\noindent\underline{\textbf{Evaluation of the representative templates.}} Figure~\ref{fig:resultsperishablemc} compares the two discovered policy templates with other baselines over $\sigma^{\rm life}=2,\ldots,15$. The advantages of RLCS are most pronounced at short shelf lives. As shelf life increases, the performance gaps among the strongest SOF policies generally narrow, consistent with the diminishing impact of expiration and the optimality of base-stock control in the nonperishable, zero-lead-time setting~\citep{snyder2019fundamentals}. In the right panel, RLCS reduces aggregate fitness by $46.2\%$ relative to Weighted Base Stock, the strongest baseline. EABS has $12.9\%$ higher aggregate fitness than this baseline but uses a simpler three-parameter expression.

\medskip
\noindent\underline{\textbf{Operational insights behind the representative templates.}}
EABS assigns less weight to inventory approaching expiration
when determining replenishment needs.
Although expiring stock can satisfy current demand, any unused
quantity cannot protect against future shortages.
EABS captures this distinction by assigning the expiring cohort
an effective weight of $\omega_1-\omega_2$, compared with
$\omega_1$ for longer-lived stock.
This mechanism follows the age-dependent weighting principle
underlying Weighted Base Stock \citep{haijema2019improved},
but uses a particularly simple distinction between expiring
and longer-lived inventory.
Our experiments show that EABS achieves performance comparable
to Weighted Base Stock.
RLCS follows the same principle but incorporates the remaining
life of every inventory cohort.
For $\varphi_{\rm IP}>0$, its remaining-life adjustment can be
written as
$
-\frac{\gamma\varphi_{\rm life}}{1+\varphi_{\rm IP}}
=
-\gamma\frac{\varphi_{\rm IP}}{1+\varphi_{\rm IP}}\,\bar m$,
where $
\bar m=\frac{\varphi_{\rm life}}{\varphi_{\rm IP}}$ is the average remaining life of carried stock
and $\gamma>0$ is a fitted coefficient controlling the
replenishment reduction.
At a fixed inventory quantity, fresher stock has a larger
$\bar m$ and produces a greater deduction from the replenishment
shortfall, reflecting its ability to satisfy demand over a longer
horizon. The factor $\varphi_{\rm IP}/(1+\varphi_{\rm IP})$ also introduces
a nonlinear dependence on total inventory.
For a fixed average remaining life, the deduction increases
with inventory quantity at a diminishing rate, making
replenishment more conservative as more stock becomes available.

\section{Theoretical Analysis and Operational Insights}\label{sec:theoretical-analysis}
Our agentic policy discovery system has identified several compact
SOF policies and templates that achieve low costs in empirical
evaluation. In this section, we isolate the additional shortfall in
the \emph{scaled capped dual-shortfall} (SCDS, Eq.~\eqref{eq:resultscappeddual}) policy.
Compared with the capped base-stock policy
(CBS, Eq.~\eqref{eq:sof-bs-cbs-co}), this component adds a bound
$[F-p_{\sigma^{\rm lead}}]^+$ that depends on the most recently
placed order in the pipeline.
Our theoretical analysis examines how this bound adjusts
replenishment decisions and can reduce inventory costs. 

To obtain tractable theoretical comparisons, we specialize the model
to nonperishable inventory, lost sales, and a $3$-period lead time
(Section~\ref{sec:model-set-theoretical}).
We then develop an average-cost decomposition that links individual
order changes to their long-run cost effects
(Section~\ref{sec:sharedstatements}).
In Section~\ref{sec:cds}, we examine how the additional shortfall
delays replenishment when a preceding order can help meet demand.
This can improve the trade-off between lost-sales and holding costs
by limiting inventory accumulation from consecutive orders.
Although the complexity of the inventory dynamics limits our
analysis to specialized model and demand classes, the resulting
operational insights may also inform policy design in more
general inventory settings. We would like to mention that SCDS also introduces a scaling factor $K$.
A separate theoretical analysis indicates that this factor can improve the cost trade-off by adapting replenishment timing to current inventory conditions; we omit that analysis due to space constraints.

\subsection{Setting of the specialized model}
\label{sec:model-set-theoretical}

We specialize the model in Section~\ref{sec:inventory-model-and-problem} to
$\sigma^{\rm lead}=3$,
$\sigma^{\rm life}=\infty$,
$\sigma^{\rm back}=0$,
$c^{\rm order}=0$.
Thus, inventory is nonperishable, orders arrive $3$ periods after
placement, and unmet demand is immediately lost. Demand is
\emph{i.i.d.}~with distribution $\mathcal D$ on nonnegative integers.
Let $D\sim\mathcal D$ denote a generic demand, with finite mean
$\mu=\E[D]$. Throughout the theoretical analysis, we use the abbreviations
$h:=c^{\rm hold}>0$ and $b:=c^{\rm lost}>0$ for the unit holding and
lost-sales costs, respectively.

\medskip
\noindent\underline{\textbf{State and dynamics.}}
Writing $x_t=x_{t,1}$, we represent the state by $s_t=(x_t,p_{t,1},p_{t,2},p_{t,3})
\in\mathbb Z_{\geq0}^{4}$,
where $x_t$ is on-hand inventory before the current receipt and
$p_{t,r}$ is the order arriving in period $t+r-1$.
The decision maker observes $s_t$ and places an order
$a_t=\pi(s_t)$, which arrives in period $t+3$.
For a generic state $s=(x,p_1,p_2,p_3)$, the inventory position is
$\operatorname{IP}(s)=x+p_1+p_2+p_3$.
After receipt and demand realization, the system evolves according to
\begin{equation}
x_{t+1}=[x_t+p_{t,1}-D_t]^+,\qquad A_t=[D_t-x_t-p_{t,1}]^+,\qquad
s_{t+1}=(x_{t+1},p_{t,2},p_{t,3},a_t),
\label{eq:transition}
\end{equation}
where $A_t$ is lost demand.
The initial state is $s_0=\bm0$, unless otherwise specified for a
local comparison.

\medskip
\noindent\underline{\textbf{Costs.}}
For a fixed instance $\Theta$, let $\Prb_\pi$ and $\E_\pi$ denote
probability and expectation under policy $\pi$. In this theoretical
analysis, we use the lower average cost
\begin{equation}
 J^\pi:=\liminf_{H\to\infty}\frac1H
 \E_\pi\!\left[\sum_{t=0}^{H-1}(h x_{t+1}+bA_t)\right].
 \label{eq:theory-average-cost}
\end{equation}
When the expected average costs converge, this equals
$J(\pi;\Theta)$ in Eq.~\eqref{eq:cost-objectives}.
For stationary state-dependent policies with a finite reachable
state space, $\mu<\infty$ ensures that the expected average
inventory and lost demand converge, and we write
$(\overline X^\pi,\overline A^\pi)
=
\lim_{H\to\infty}\frac1H
\E_\pi\!\left[\sum_{t=0}^{H-1}(x_{t+1},A_t)\right]$.
In this case, the $\liminf$ in Eq.~\eqref{eq:theory-average-cost}
is an ordinary limit and $J^\pi=h\overline X^\pi+b\overline A^\pi$.

\medskip
\noindent\underline{\textbf{Policies and benchmark.}}
As our benchmark, we consider a generalized class of CBS policies with arbitrary time-dependent order caps:
\begin{equation}
 \pi_{\rm CBS}(s,t;\bm C,S)
 =\min\{C_t,[S-\operatorname{IP}(s)]^+\},\qquad\text{where }S\in{\mathbb N}_{\geq 0}, ~\bm C=(C_t)_{t\ge 0}, ~C_t\in{\mathbb N}_{\geq 0}.
 \label{eq:policy-cbs}
\end{equation}
In $\pi_{\rm CBS}(\cdot;\bm C,S)$, the dot denotes the state and
time arguments $(s,t)$. When $C_t\equiv C$, we write
$\pi_{\rm CBS}(\cdot;C,S)$ and omit the time argument, so the dot
denotes only the state $s$. Taking $C_t=\xi(C,t)$, with $\xi$ defined in
Eq.~\eqref{eq:sof-bs-cbs-co}, recovers the CBS policy introduced
there for every real $C\ge0$. We refer to the case with integer $C$ as the fixed-cap.

To examine the operational role of the additional shortfall, we
incorporate it into the fixed-cap CBS policy, yielding
the capped dual-shortfall (CDS) policy:
\begin{equation}
\pi_{\rm CDS}(s;C,S,F)
=
\left\lfloor
\min\{C,[S-\operatorname{IP}(s)]^+,[F-p_3]^+\}
\right\rfloor. \label{eq:policy-cds}
\end{equation}
Here $C,S,F$ are finite nonnegative parameters.
Our benchmark allows every
integer target and every exogenous adjacent-integer cap sequence (namely, the \emph{adjacent-integer CBS}):
\begin{equation}
 J^{\mathrm{AI-CBS},*}
 :=\inf_{S,m\in\mathbb Z_{\ge0}}\ \inf_{\bm C\in\{m,m+1\}^{\mathbb Z_{\ge0}}}
       J^{\pi_{\rm CBS}(\cdot;\bm C,S)}.
 \label{eq:cbsbenchmark}
\end{equation}
This class includes every sequence $C_t=\xi(C,t)$ from
Eq.~\eqref{eq:sof-bs-cbs-co} and every constant integer cap.

\medskip
\noindent\underline{\textbf{Demand tails and IFR.}}
Part of our analysis focuses on demand distributions with an
\emph{increasing failure rate} (IFR).
For the demand distribution $\mathcal D$, define
$\mathbb{S}_{\mathcal D}(j)
:=\Prb(D\geq j)$ for all $j\in\mathbb Z_{\geq0}$,
so $\mathbb{S}_{\mathcal D}(0)=1$.
The distribution $\mathcal D$ is IFR if
$\mathbb{S}_{\mathcal D}(j+1)/\mathbb{S}_{\mathcal D}(j)$
is nonincreasing in $j$ wherever $\mathbb{S}_{\mathcal D}(j)>0$.
The IFR assumption holds for several common discrete demand
distributions, including Poisson, binomial, geometric, and uniform
distributions on consecutive nonnegative integers.

\subsection{An average-cost decomposition lemma}
\label{sec:sharedstatements}

We introduce Lemma~\ref{lem:averaging}, which decomposes the average-cost
difference between a reference fixed-cap CBS policy $\varsigma$ and a candidate
stationary state-dependent policy $\pi$ using the effects of a single order deviation. Fix a state $s\in\mathbb Z_{\geq0}^4$ and an order quantity
$a\in\mathbb Z_{\geq0}$. Consider two systems starting from
$s_0^\varsigma=\widetilde s_0=s$ and facing the same
\emph{i.i.d.}~demand sequence with distribution $\mathcal D$.
The reference system follows $\varsigma$ throughout, while the
perturbed system orders $a$ in period zero and follows $\varsigma$
from its own state thereafter:
\[
a_t^\varsigma=\varsigma(s_t^\varsigma),\qquad
\widetilde a_t=
\begin{cases}
a,&t=0,\\
\varsigma(\widetilde s_t),&t\geq1.
\end{cases}
\]
Both systems evolve according to \eqref{eq:transition}.
Let $x_{t+1}^\varsigma$ and $\widetilde x_{t+1}$ denote their
respective ending inventories in period $t$, and let
$A_t^\varsigma$ and $\widetilde A_t$ denote their lost-sales quantities.
Whenever the sums below are integrable, define the expected cumulative
inventory reduction and additional lost sales caused by the perturbation as
\begin{equation}
\mathcal X^\varsigma(s,a)
=\E\!\left[\sum_{t=0}^{\infty}
(x_{t+1}^\varsigma-\widetilde x_{t+1})\right],\qquad
\mathcal A^\varsigma(s,a)
=\E\!\left[\sum_{t=0}^{\infty}
(\widetilde A_t-A_t^\varsigma)\right].
\label{eq:paired}
\end{equation}
We also define the long-run state-occupation measure under $\pi$,
starting from $s_0=\bm0$. For $s\in\mathbb Z_{\geq0}^4$ and
$E\subseteq\mathbb Z_{\geq0}^4$, let
\begin{equation}
d_\pi(s)
=\lim_{H\to\infty}\frac1H\sum_{t=0}^{H-1}\Prb_\pi(s_t=s),
\qquad
d_\pi(E)=\sum_{s\in E}d_\pi(s).
\label{eq:statelaws}
\end{equation}
These limits exist for stationary state-dependent policies
with a finite reachable state space, including the fixed-cap CBS
and CDS policies used in the local comparisons below.
We now state the main lemma of this subsection, whose proof is given
in Section~\ref{sec:averagingproof-new}.

\begin{lemma}[Average-cost decomposition]\label{lem:averaging}
Let $\varsigma$ be a fixed-cap CBS policy with nonnegative
integer parameters and $\pi$ a stationary state-dependent
candidate policy.
Suppose both policies preserve a common finite state set
containing $\bm0$, and the demand distribution satisfies $0<\Prb(D=0)<1$.
Then the long-run averages and the limits in Eq.~\eqref{eq:statelaws}
exist, and
\[
 \overline X^\varsigma-\overline X^\pi
   =\sum_{s:d_{\pi}(s)>0} d_\pi(s)\mathcal X^\varsigma(s,\pi(s)),\qquad
 \overline A^\pi-\overline A^\varsigma
   =\sum_{s:d_{\pi}(s)>0} d_\pi(s)\mathcal A^\varsigma(s,\pi(s)).
\]
Consequently,
\begin{equation}
 \frac{J^\varsigma-J^\pi}{h}
 =\sum_{s:d_{\pi}(s)>0} d_\pi(s)
       \left[\mathcal X^\varsigma(s,\pi(s))-\frac bh\mathcal A^\varsigma(s,\pi(s))\right].
 \label{eq:pd}
\end{equation}
\end{lemma}

\subsection{Theoretical interpretation of the additional shortfall}\label{sec:cds}
In this subsection, we compare the capped dual-shortfall policy
(CDS, Eq.~\eqref{eq:policy-cds}) with the fixed-cap CBS policy
(Eq.~{\eqref{eq:policy-cbs}} with $C_t=C_0$ for all $t$) to examine the role of the
additional shortfall $[F-p_3]^+$.
This shortfall limits the combined size of two consecutive
orders to $F$. A large preceding order can therefore constrain the current
order and delay part of the replenishment, preventing large
orders from accumulating in consecutive periods. As the analysis below shows, this mechanism can improve the
trade-off between holding and lost-sales costs by reducing
inventory carried while the preceding order helps meet demand
during the delay.

For each integer $C\geq1$, define the CBS benchmark and CDS policy
\begin{equation}
\varsigma_C=\pi_{\mathrm{CBS}}(\,\cdot\,;C,2C),
\qquad
\pi^{\rm dual}_C=\pi_{\mathrm{CDS}}(\,\cdot\,;C,2C,2C-1).
\label{eq:cdspolicies}
\end{equation}
Section~\ref{sec:cdsmain} compares these policies to examine the
holding savings and additional lost sales caused by delayed
replenishment.
Section~\ref{sec:cdsglobal-main} then shows that, under suitable
demand and cost conditions, $\pi^{\rm dual}_1$ outperforms every CBS
policy with nonnegative integer parameters.

\subsubsection{$\pi^{\rm dual}_C$ versus $\varsigma_C$:
delaying with a preceding order}
\label{sec:cdsmain}
We first compare $\pi^{\rm dual}_C$ with the CBS benchmark
$\varsigma_C$.
Among states reachable from an empty system, the two policies
differ only at $s_C^*=(0,0,0,C)$, where the preceding order of
$C$ units arrives in two periods.
Lemma~\ref{lem:averaging} reduces the comparison to paired
trajectories generated by a single order deviation at this state.
Taking the initial period as time $0$, the CDS and CBS
trajectories initially order $C-1$ and $C$ units, respectively,
and both subsequently follow $\varsigma_C$ from their own states.
Their orders in periods $0$ and $1$ are therefore $[C-1,1]$
and $[C,0]$, respectively.
Thus, CDS delays the arrival of one unit from period $3$
to period $4$.

Delaying the arrival creates a trade-off between holding-cost
savings and potential increases in lost sales.
If the delay causes no additional lost sale in period $3$,
CDS carries one unit less in that period.
The two systems then have identical orders, sales, and ending
inventories in period $4$, and their full states coincide at
the beginning of period $5$.
Thus, the delay saves exactly one unit-period of holding
in this case.
The preceding order helps meet demand before the delayed unit
arrives.
In period $3$, the CDS trajectory has $[C-D_2]^++C-1$ units
available, one fewer than the CBS trajectory.
An additional lost sale therefore occurs exactly when
$D_3\geq C$ and $D_2+D_3\geq2C$.
For independent demands $D_1,D_2\sim\mathcal D$, define
$
\psi_C
\coloneqq
\Prb(D_1\geq C,\ D_1+D_2\geq2C)$.
By the \emph{i.i.d.}~demand assumption, $\psi_C$ is the probability
of this additional lost sale.
The following lemma characterizes the resulting holding savings
and additional lost sales, including their effects in subsequent
periods; its proof appears in
Section~\ref{sec:cdsmechanism-new}.
\begin{lemma}
\label{lem:cdsmechanism}
Let $C\geq1$ be an integer, and let $\mathcal D$ be a demand
distribution with finite mean $\mu$, not necessarily IFR,
satisfying $q\coloneqq\mathbb S_{\mathcal D}(1)<1$.
For fixed $C$, uniformly as $\mu\downarrow0$, we have
\[
\overline X^{\varsigma_C}-\overline X^{\pi^{\rm dual}_C}
=\psi_C\bigl(1+O_C(\mu)\bigr),
\qquad
\overline A^{\pi^{\rm dual}_C}-\overline A^{\varsigma_C}
=\psi_C^2\bigl(1+O_C(\mu)\bigr).
\]
\end{lemma}
Consequently, when $\psi_C>0$ and $\mu$ is sufficiently small,
the CDS policy $\pi^{\rm dual}_C$ has a lower average cost than
the CBS benchmark $\varsigma_C$ whenever
$
h/b>
\psi_C\bigl(1+O_C(\mu)\bigr)$.

\subsubsection{Strict improvement over adjacent-integer CBS policies}
\label{sec:cdsglobal-main}
We now set $C=1$ and show that, for a class of IFR demand
distributions, the CDS policy $\pi^{\rm dual}_1$ achieves a strictly
lower average cost than every adjacent-integer CBS policy in Eq.~\eqref{eq:cbsbenchmark}.

\begin{theorem}
\label{thm:cds}
Let $\mathcal D$ be a nonnegative integer IFR demand distribution
with finite mean $\mu>0$.
Suppose
\begin{equation}
\theta\coloneqq\frac{b}{h}\mu^2\in(1/3,1/2),
\qquad
\mu\leq
\min\left\{
\frac1{100},
\frac{1-2\theta}{12},
\frac{3\theta-1}{48}
\right\}.
\label{eq:cdswindow}
\end{equation}
Then the policy
$\pi^{\rm dual}_1=\pi_{\mathrm{CDS}}(\,\cdot\,;1,2,1)$ satisfies
$
J^{\pi^{\rm dual}_1}
\leq
{J^{\mathrm{AI-CBS},*}}
-\frac{h}{4} \cdot (1-2\theta)\mu^2$.
\end{theorem}
The proof of Theorem~\ref{thm:cds} appears in
Section~\ref{sec:cdstheorem-new}.
It first uses cost bounds for CBS policies with
exogenous adjacent-integer caps to reduce the desired
inequality to a comparison with $\varsigma_1$.
A refined version of Lemma~\ref{lem:cdsmechanism} then
establishes the stated cost gap by quantifying the holding-cost
savings and additional lost-sales costs of delayed replenishment.
Although the theorem applies to a specific CDS policy, its
proof illustrates how the additional shortfall can improve
the trade-off between holding and lost-sales costs by delaying
replenishment after a sizable preceding order.
This mechanism may also reduce costs in more general
inventory settings.

\section{Conclusion}
\label{sec:conclusion}

We develop an agentic system that discovers compact SOF policies and templates through multicriteria selection, operational guidance, motif reuse, and progressive generalization. Experiments on 15,600 instances across four inventory settings show aggregate performance gains over SOF benchmarks and competitive performance against policies requiring more elaborate online computations.

For a lost-sales system with a three-period lead time, we prove that a capped dual-shortfall specialization strictly outperforms the considered CBS benchmark class under specified demand and cost conditions. The analysis explains how delaying replenishment after a preceding order can yield holding-cost savings that outweigh additional lost-sales costs. These findings make initial connection between our automated policy discovery with interpretable operational mechanisms and rigorous performance analysis.  A central and interesting direction for future research is to establish theoretical performance guarantees for more discovered policy families across the inventory settings studied.

\bibliographystyle{informs2014}
\bibliography{main}

\ECSwitch
\SetCompanionAnchors

\ECHead{Electronic Companion}
\section{{Cost decomposition by single-order changes}}\label{sec:cost-decomposition}
{
We first establish a finite-horizon comparison that allows time-dependent
orders. We then recover Lemma~\ref{lem:averaging} by stationary
occupation averaging and give a common method for computing the
single-order effects used in the cap-sequence comparisons.

\subsection{A finite-horizon cost decomposition}
\label{sec:finite-horizon-decomposition}
For the paired effects in Eq.~\eqref{eq:paired}, put
\[
 \Delta_B(s,a):=\frac bh\mathcal A^B(s,a)-\mathcal X^B(s,a).
\]

\begin{lemma}
\label{lem:finite-horizon-decomposition}
Suppose $\E D<\infty$ and $0<\Prb(D=0)<1$.
Let $B$ be a fixed-cap CBS with nonnegative integer parameters.
Let $\pi$ start empty, choose orders before observing current demand, and stay in
a given finite state set $\mathsf X$.
There is a constant $K<\infty$, uniform over $H$ and all such
nonanticipative policies $\pi$, such that
\begin{equation}
 \left|\frac1h\E_\pi\sum_{t<H}(h x_{t+1}+bA_t)
       -\frac{HJ^B}{h}
       -\sum_{t<H}\E_\pi\Delta_B(s_t,a_t)\right|\le K.
 \label{eq:finite-oneorder-bound}
\end{equation}
Consequently,
\[
 \inf_{H\ge1}\sum_{t<H}\E_\pi\Delta_B(s_t,a_t)>-\infty
 \quad\Longrightarrow\quad J^\pi\ge J^B.
\]
\end{lemma}

\noindent\emph{Proof.}
We first show that, after a single order change, the two paths
following $B$ meet with a geometrically bounded waiting time.
This makes the single-order effects well defined and bounds the
total error from truncating them at horizon $H$, uniformly over
the candidate policies in the lemma. We then replace the
candidate's orders by $B$'s orders one period at a time, obtaining
an exact finite-horizon telescoping identity. Combining this
identity with the truncation bound leaves only an $O(1)$ error,
which proves Eq.~\eqref{eq:finite-oneorder-bound}; dividing by $H$
gives the final comparison.

\noindent\underline{\textbf{Step 1: uniformly summable truncation errors.}}
For the paired paths from $(s,a)$ in Eq.~\eqref{eq:paired}, define
\[
 \begin{aligned}
 \mathcal X_n^B(s,a)
 =\E\sum_{j=0}^{n-1}(x_{j+1}^B-\widetilde x_{j+1}),~\mathcal A_n^B(s,a)=\E\sum_{j=0}^{n-1}(\widetilde A_j-A_j^B),~
 \Delta_{B,n}(s,a)=\frac bh\mathcal A_n^B(s,a)-\mathcal X_n^B(s,a).
 \end{aligned}
\]
Write $B(s)=\min\{C,[S-\operatorname{IP}(s)]^+\}$,
$S,C\in\mathbb Z_{\ge0}$.
Choose an integer $M\ge S$ such that
$\operatorname{IP}(s)+a\le M$ for every state-order pair $(s,a)$
used by $\pi$.
Subsequent reference paths remain bounded since
\[
 \operatorname{IP}(s_{t+1})\le
 \operatorname{IP}(s_t)+B(s_t)
 \le\max\{S,\operatorname{IP}(s_t)\}.
\]
Choose $d\ge1$ with $\Prb(D=d)>0$.
The demand block
\[
 (\underbrace{d,\ldots,d}_{M+4},
  \underbrace{0,\ldots,0}_{S+4})
 \quad\text{resets every reference path to}\quad
 \begin{cases}
 (S,0,0,0),&SC>0,\\
 \bm0,&SC=0.
 \end{cases}
\]
Thus, with $m=M+S+8$ and
$p=\Prb(D=d)^{M+4}\Prb(D=0)^{S+4}>0$, the paired coalescence time
$\tau=\inf\{t\ge1:s_t^B=\widetilde s_t\}$ satisfies
\[
 \Prb(\tau>1+km)\le(1-p)^k,\qquad
 \E\sum_{t\ge0}
 \bigl(|x_{t+1}^B-\widetilde x_{t+1}|+
       |\widetilde A_t-A_t^B|\bigr)\le2M\E\tau<\infty.
\]
Uniformly over possible $(s,a)$,
\[
 \begin{aligned}
 |\mathcal X_n^B-\mathcal X^B|,
 |\mathcal A_n^B-\mathcal A^B|\le M\E[(\tau-n)^+],\qquad
 |\Delta_{B,n}-\Delta_B|\le(1+b/h)M\E[(\tau-n)^+].
 \end{aligned}
\]
For $(F_n,F)=(\mathcal X_n^B,\mathcal X^B)$,
$(\mathcal A_n^B,\mathcal A^B)$, $(\Delta_{B,n},\Delta_B)$, set $\varepsilon_n:=(1+b/h)M\sup_{(s,a)}\E[(\tau-n)^+]$, and we have
\begin{align}
     \sup_{(s,a)}|F_n(s,a)-F(s,a)|&\le\varepsilon_n,
 \qquad \sum_{n\ge1}\varepsilon_n<\infty,
\end{align}
indicating
\begin{equation}
 \sup_\pi\left|
 \sum_{t=0}^{H-1}\E_\pi
 [F_{H-t}(s_t,a_t)-F(s_t,a_t)]\right|
 \le\sum_{n=1}^H\varepsilon_n
 \le\sum_{n\ge1}\varepsilon_n<\infty
 \qquad(H\ge1).
 \label{eq:oneorder-uniform-replacement}
\end{equation}

\medskip
\noindent\underline{\textbf{Step 2: telescope the single-order replacements.}}
Let $P_k$ be the policy that uses $\pi$ for $t<k$ and $B$ thereafter, for each $0\le k\le H$.
The pair $P_k,P_{k+1}$ has the same pre-demand history at $k$,
distributed as under $\pi$: the order changes from $B(s_k)$ to
$a_k$, and both policies follow $B$ thereafter.
Future demands are independent of this history.
Consequently,
\begin{align}
  \E_\pi\sum_{t<H}A_t-\E_B\sum_{t<H}A_t
 =\sum_{t<H}\E_\pi\mathcal A_{H-t}^B(s_t,a_t),\quad\E_B\sum_{t<H}x_{t+1}-\E_\pi\sum_{t<H}x_{t+1}
 =\sum_{t<H}\E_\pi\mathcal X_{H-t}^B(s_t,a_t).
 \label{eq:decomp-new-finite}
\end{align}

\medskip
\noindent\underline{\textbf{Step 3: remove truncation and initialization errors.}}
Apply Eq.~\eqref{eq:oneorder-uniform-replacement} to
$\mathcal X_{H-t}^B$, $\mathcal A_{H-t}^B$, and their cost
combination $\Delta_{B,H-t}$. Coupling its empty-start and stationary paths gives
\[
 \sup_{H\ge1}\left|
 \frac1h\E_B\sum_{t<H}(h x_{t+1}+bA_t)-\frac{HJ^B}{h}
 \right|\le K_B<\infty.
\]
Hence Eq.~\eqref{eq:finite-oneorder-bound} holds with
$K=K_B+\sum_{n\ge1}\varepsilon_n$.
Dividing by $H$ and taking the lower limit proves the implication.

The same telescoping and truncation bounds, applied separately to
ending inventory and lost sales, give
\begin{align}
 \E_B\sum_{t<H}x_{t+1}-\E_\pi\sum_{t<H}x_{t+1}
 &=\sum_{t<H}\E_\pi\mathcal X^B(s_t,a_t)+O(1),
 \label{eq:finite-holding-bound}\\
 \E_\pi\sum_{t<H}A_t-\E_B\sum_{t<H}A_t
 &=\sum_{t<H}\E_\pi\mathcal A^B(s_t,a_t)+O(1),
 \label{eq:finite-loss-bound}
\end{align}
uniformly in $H$ and $\pi$ on $\mathsf X$.
\hfill$\square$

\subsection{Proof of Lemma~\ref{lem:averaging}}
\label{sec:averagingproof-new}
We apply the component comparisons above with $B=\varsigma$.
Both stationary policies preserve the common finite state set
$\mathsf X$ in the lemma, so their state processes are finite
Markov chains. Their Ces\`aro state probabilities converge:
\[
 \frac1H\sum_{t<H}\Prb_\xi(s_t=s)\longrightarrow d_\xi(s),
 \qquad \xi\in\{\pi,\varsigma\},\quad s\in\mathsf X.
\]
With $M=\max_{s\in\mathsf X}\operatorname{IP}(s)$, the conditional expected ending inventory and lost sales are bounded by $M$ and $\E D$, respectively. State averaging therefore proves
existence and finiteness of $\overline X^\xi$ and
$\overline A^\xi$. Absolute integrability of the paired effects
follows from Step 1 of
Section~\ref{sec:finite-horizon-decomposition}.

Since $a_t=\pi(s_t)$, divide Eq.~\eqref{eq:finite-holding-bound}
by $H$ and group its right side by state:
\[
 \frac1H\sum_{t<H}\E_\pi\mathcal X^\varsigma(s_t,\pi(s_t))
 =\sum_{s\in\mathsf X}
   \left(\frac1H\sum_{t<H}\Prb_\pi(s_t=s)\right)
   \mathcal X^\varsigma(s,\pi(s)).
\]
The sum is finite, so taking limits gives $\overline X^\varsigma-\overline X^\pi=\sum_{s:d_\pi(s)>0}d_\pi(s)\mathcal X^\varsigma(s,\pi(s))$. Applying the same argument to Eq.~\eqref{eq:finite-loss-bound} gives $\overline A^\pi-\overline A^\varsigma
 =\sum_{s:d_\pi(s)>0}d_\pi(s)\mathcal A^\varsigma(s,\pi(s))$.
Finally, combining the two identities with
$J^\varsigma-J^\pi
=h(\overline X^\varsigma-\overline X^\pi)
-b(\overline A^\pi-\overline A^\varsigma)$ proves Eq.~\eqref{eq:pd}.

\subsection{Computing the single-order effects}
\label{sec:computing-oneorder}
For the calculations below, combine the first two state coordinates:
$u=(y,r,p)=(x+p_1,p_2,p_3)$, with transition
$\Phi(u,a,D)=((y-D)^++r,p,a)$.
The reduced state determines costs, future reduced states, and the
reference order, so $\Delta_B(s,a)$ depends on $s$ only through $u$.
We accordingly write $\Delta_B(u,a)$ for the same single-order effect.

Let $c(y)=\E\bigl[(y-D)^++(b/h)(D-y)^+\bigr]$, $\rho=J^B/h$. For a reduced state $u$, let $Y_t^u$ be the available stock at
time $t$ when the system starts from $u$ and follows $B$ throughout. Define $f(u):=\sum_{t=0}^{\infty}
\bigl(\E_u c(Y_t^u)-\E_{\bm0}c(Y_t^{\bm0})\bigr)$.
The common reset makes this series absolutely convergent on every
fixed finite state set, and $f(\bm0)=0$.
Under $B$, the first period costs $c(y)$ and leaves the system in
state $\Phi(u,B(u),D)$. In the definition of $f$, advancing the
empty-state comparison by one period subtracts its limiting
one-period cost $\rho$. Hence
{
\begin{equation}
 f(u)=c(y)+\E f(\Phi(u,B(u),D))-\rho,\qquad f(\bm0)=0.
 \label{eq:reference-bias-recursion}
\end{equation}
}
Changing the first order to $a$ leaves the first-period cost
unchanged, so
\begin{equation}
\begin{aligned}
\Delta_B(u,a)
=\E f(\Phi(u,a,D))-\E f(\Phi(u,B(u),D))=c(y)+\E f(\Phi(u,a,D))-f(u)-\rho.
\end{aligned}
\label{eq:calendar-residual}
\end{equation}
}

\section{Proofs for Section~\ref{sec:cds}}
\label{sec:cdsproof-new}
We use the policies $\varsigma_C,\pi^{\rm dual}_C$ from
Eq.~\eqref{eq:cdspolicies} and the paired effects from
Eq.~\eqref{eq:paired}. Recall that
\[
 s_C^*=(0,0,0,C),\qquad
 \psi_C:=\Prb(D_1\ge C,\ D_1+D_2\ge2C),
\]
where $D_1,D_2$ are independent demands with law $\mathcal D$.
Throughout this section, write
\[
 T_j=\mathbb S_{\mathcal D}(j),\qquad q=T_1,\qquad v=1-q,\qquad
 \mu=\E D,\qquad \beta=b/h.
\]
Sales are $V_t=\min\{x_t+p_{t,1},D_t\}=D_t-A_t$.
For each comparison, superscript $\varsigma$ and a tilde denote
the reference and altered paths as in Section~\ref{sec:sharedstatements}.
Each subsection uses this notation; further quantities are defined locally.

\subsection{Auxiliary bounds for the CDS comparison}
For IFR demand, the tail ratios imply
$T_n\le T_jq^{n-j}$ for $n\ge j$.
For $j\ge1$, define the binomial moment $\nu_j:=\E\binom Dj$.
The identity $\binom Dj=\sum_{n=j}^{D}\binom{n-1}{j-1}$
and the tail bound $T_n\le T_jq^{n-j}$ for $n\ge j$ give
\begin{equation}
 \nu_j=\sum_{n\ge j}\binom{n-1}{j-1}T_n
 \le T_j\sum_{k\ge0}\binom{k+j-1}{j-1}q^k
 =\frac{T_j}{v^j},
 \label{eq:binomial-tail-moment}
\end{equation}
where the last equality uses the negative-binomial series and $v=1-q$.

{
\begin{lemma}\label{lem:positive-cap-holding}
Fix an integer target $S\ge1$. Starting empty, let a CBS policy use
integer caps $C_t\ge1$ chosen before the current demand. If
$\mu=\E D<1/2$ and $\E D^2<\infty$, then, with
$Q=\E[D(D-1)]/[2(1-\mu)]$,
\begin{equation}
 \liminf_{H\to\infty}\frac1H\sum_{t<H}\E x_{t+1}
 \ge S-4\mu-Q.
 \label{eq:calendar-holding}
\end{equation}
\end{lemma}
\begin{proof}
Put $\mathcal R_t:=S-\operatorname{IP}(s_t)\in\{0,\ldots,S\}$.
Inventory conservation and $C_t\ge1$ give
\[
 \mathcal R_{t+1}=(\mathcal R_t-C_t)^++V_t
 \le\mathcal R_t-\ind_{\{\mathcal R_t>0\}}+D_t.
\]
Taking expectations and summing over $t<H$ yields
$\sum_{t<H}\Prb(\mathcal R_t>0)\le H\mu+S$.
Put $I_t:=\ind_{\{\mathcal R_t>0\}}$.
Since $I_t\mathcal R_t=\mathcal R_t$ and $D_t$ is independent of
$\mathcal R_t$, squaring the recursion bound gives
\[
 \E(\mathcal R_{t+1}^2-\mathcal R_t^2)
 \le-2(1-\mu)\E\mathcal R_t
       +(1-2\mu)\Prb(\mathcal R_t>0)+\E D^2.
\]
Summing over $t<H$, using $\mathcal R_0=S$ and
$\mathcal R_H^2\ge0$, yields
\[
 2(1-\mu)\sum_{t<H}\E\mathcal R_t
 \le(1-2\mu)\sum_{t<H}\Prb(\mathcal R_t>0)
       +H\E D^2+S^2.
\]
Together with $\sum_{t<H}\Prb(\mathcal R_t>0)\le H\mu+S$
and $\mu<1/2$, this gives
$\limsup_H H^{-1}\sum_{t<H}\E\mathcal R_t\le\mu+Q$.

Inventory-position conservation and $\operatorname{IP}(s_0)=0$ give $\sum_{t<H}\E a_t
 =\E\operatorname{IP}(s_H)+\sum_{t<H}\E V_t
 \le S+H\mu$. Because the system starts empty, $p_{t,j}=a_{t-(4-j)}$ for
$j=1,2,3$, with $a_k=0$ for $k<0$. Hence
$\sum_{t<H}\E p_{t,j}\le S+H\mu$ for each $j$.
Using $\operatorname{IP}(s_t)=S-\mathcal R_t$ then gives
\[
 \frac1H\sum_{t<H}\E x_t
 =S-\frac1H\sum_{t<H}\E\mathcal R_t
       -\frac1H\sum_{j=1}^3\sum_{t<H}\E p_{t,j}
 \ge S-\frac1H\sum_{t<H}\E\mathcal R_t-3\mu-\frac{3S}{H}.
\]
Finally, $0\le x_H\le S$ and $x_0=0$, so
\[
 \left|\frac1H\sum_{t<H}\E x_{t+1}
       -\frac1H\sum_{t<H}\E x_t\right|
 =\frac{\E x_H}{H}\le\frac SH.
\]
Taking $H\to\infty$ proves Eq.~\eqref{eq:calendar-holding}.
\end{proof}
}

\subsection{Proof of Lemma~\ref{lem:cdsmechanism}}
\label{sec:cdsmechanism-new}
Fix an integer $C\ge1$ and assume only iid nonnegative integer demand
with finite mean and $q<1$. Define
\[
 \alpha_C:=d_{\pi^{\rm dual}_C}(s_C^*),\qquad
 \Gamma_{X,C}:=\mathcal X^{\varsigma_C}(s_C^*,C-1),\qquad
 \Gamma_{A,C}:=\mathcal A^{\varsigma_C}(s_C^*,C-1).
\]
We first bound $\alpha_C$, then $\Gamma_{X,C}$ and
$\Gamma_{A,C}$, and finally derive the long-run differences. The resulting bounds and identities will also be used in the
theorem proof:
\begin{equation}
 \begin{gathered}
 v^6\psi_C\le\alpha_C\le\psi_C,\\
 |\Gamma_{X,C}-1|\le(1+10C/v^5)\psi_C,\qquad
 |\Gamma_{A,C}-\psi_C|\le(5\mu/v^5)\psi_C,\\
 \overline X^{\varsigma_C}-\overline X^{\pi^{\rm dual}_C}
       =\alpha_C\Gamma_{X,C},\qquad
 \overline A^{\pi^{\rm dual}_C}-\overline A^{\varsigma_C}
       =\alpha_C\Gamma_{A,C}.
 \end{gathered}
 \label{eq:cds-explicit-new}
\end{equation}

\medskip
\noindent\underline{\textbf{Step 1: the changed order and stationary frequency.}}
We first show that, on the common finite state set below, the two
policies differ only at $s_C^*=(0,0,0,C)$, and that its stationary
probability under $\pi_C^{\rm dual}$ satisfies
\[
 v^6\psi_C\le\alpha_C\le\psi_C.
\]
Both policies preserve the finite set
\[
 \mathsf X_C^{\rm D}:=
 \{(x,p_1,p_2,p_3)\in\mathbb Z_{\ge0}^4:
       x+p_1+p_2+p_3\le2C,\ p_1,p_2,p_3\le C\}.
\]
The additional cap $2C-1-p_3$ can reduce the CBS order only when
$p_3=C$ and that order is $C$. This requires
$\operatorname{IP}(s)\le C$, hence $s=s_C^*$; the orders there
are $C$ and $C-1$.

Six consecutive zero demands bring CDS from any state in this set
to $f_C:=(2C,0,0,0)$; hence $d_{\pi_C^{\rm dual}}(f_C)\ge v^6$.

Run CDS in stationarity and take $t\ge3$. Define
\[
 E_t:=\{D_{t-2}\ge C,\ D_{t-2}+D_{t-1}\ge2C\},\qquad
 F_t:=\{s_{t-2}=f_C\}\cap E_t.
\]
The event $F_t$ is sufficient and $E_t$ is necessary for $s_t=s_C^*$:
\[
 F_t\subseteq\{s_t=s_C^*\}\subseteq E_t.
\]
Taking probabilities gives
\[
 v^6\psi_C\le d_{\pi^{\rm dual}_C}(f_C)\psi_C
 =\Prb(F_t)\le\alpha_C\le\Prb(E_t)=\psi_C.
\]

\medskip
\noindent\underline{\textbf{Step 2: holding and loss effects of the delay.}}
Start the paired systems at $s_C^*$ at time zero.
Their first two orders are $(C,0)$ and $(C-1,1)$; both follow
$\varsigma_C$ from time one. One unit is therefore delayed from
time three to time four. The event of an extra loss at time three is
\[
 E_C:=\{D_3\ge C,\ D_2+D_3\ge2C\},\qquad \Prb(E_C)=\psi_C.
\]
On $E_C^c$, the holding saving is one in period three and future
costs agree from period four onward. On $E_C$, both ending
inventories in period three are zero and the altered system has
one extra lost sale.

To bound the remaining effects, divide demands starting at time
four into disjoint blocks of five. Let $T$ be the number of periods
through the end of the first all-zero block. Five consecutive zero
demands bring CBS from any state in $\mathsf X_C^{\rm D}$ to $f_C$,
so the full states agree by time $4+T$, and
\[
 \E T=\frac5{v^5},\qquad
 \E\sum_{t=4}^{3+T}D_t=\mu\,\E T=\frac{5\mu}{v^5}.
\]
For the second identity, write the sum as
$\sum_{n\ge0}\ind_{\{T>n\}}D_{4+n}$: the event $\{T>n\}$
is determined before $D_{4+n}$, so independence applies.
Both $T$ and this demand sum depend only on demands from time
four onward and are independent of $E_C$.

Define the remaining holding and loss differences by
\[
 H_{\rm rem}:=\sum_{t=4}^{3+T}(x_{t+1}^{\varsigma_C}
                  -\widetilde x_{t+1}),\qquad
 L_{\rm rem}:=\sum_{t=4}^{3+T}(\widetilde A_t-A_t^{\varsigma_C}).
\]
They vanish on $E_C^c$, and satisfy
$|H_{\rm rem}|\le2CT$ and
$|L_{\rm rem}|\le\sum_{t=4}^{3+T}D_t$.
Thus all paired effects are absolutely integrable, and
\[
 \Gamma_{X,C}=1-\psi_C+\E[\ind_{E_C}H_{\rm rem}],\qquad
 \Gamma_{A,C}=\psi_C+\E[\ind_{E_C}L_{\rm rem}].
\]
Bounding the two corrections by the independent future duration
and demand sum yields
\[
 |\Gamma_{X,C}-1|\le\psi_C+\frac{10C}{v^5}\psi_C,\qquad
 |\Gamma_{A,C}-\psi_C|\le\frac{5\mu}{v^5}\psi_C.
\]

\medskip
\noindent\underline{\textbf{Step 3: average and conclude the asymptotics.}}
Only $s_C^*$ has a changed order, so Lemma~\ref{lem:averaging}
gives the two stationary identities in Eq.~\eqref{eq:cds-explicit-new}.
When $q=0$, both policies eventually hold $2C$ units, and the
same identities hold with $\alpha_C=\psi_C=0$.

If $\psi_C>0$, the frequency bounds give
$|1-\alpha_C/\psi_C|\le1-v^6\le6\mu$.
Also $\psi_C\le\mu/C$ and $5/v^5\le5/(1-\mu)^5$ for $\mu<1$.
Consequently,
\[
 \begin{aligned}
 \left|\frac{\overline X^{\varsigma_C}
              -\overline X^{\pi^{\rm dual}_C}}{\psi_C}-1\right|
 \le6\mu+(1+10C/v^5)\psi_C
 \le(6+1/C+10/v^5)\mu,\qquad
 \left|\frac{\overline A^{\pi^{\rm dual}_C}
              -\overline A^{\varsigma_C}}{\psi_C^2}-1\right|
 &\le6\mu+5\mu/v^5.
 \end{aligned}
\]
For fixed $C$, these are uniform $O_C(\mu)$ bounds. If $\psi_C=0$,
then $\alpha_C=0$ and both stationary differences vanish.
Thus in all cases
\[
 \overline X^{\varsigma_C}-\overline X^{\pi^{\rm dual}_C}
 =\psi_C(1+O_C(\mu)),\qquad
 \overline A^{\pi^{\rm dual}_C}-\overline A^{\varsigma_C}
 =\psi_C^2(1+O_C(\mu)),
\]
which proves Lemma~\ref{lem:cdsmechanism} without an IFR assumption.

\subsection{Proof of Theorem~\ref{thm:cds}}
\label{sec:cdstheorem-new}
In this section, we first state four screening lemmas, two for fixed caps and two
for time-varying caps. We then derive the IFR demand bounds, establish the cost gap against $\varsigma_1$ using Lemma~\ref{lem:cdsmechanism}, and apply the screening lemmas to
cover every adjacent-integer CBS policy. 

Assume Eq.~\eqref{eq:cdswindow}, including IFR demand, and recall
$\theta=\beta\mu^2\in(1/3,1/2)$.
\begin{lemma}\label{lem:cdstargetscreen-new}
Under the assumptions of Theorem~\ref{thm:cds}, every pair of
nonnegative integers $S,C$ with $S\ne2$ or $C=0$ satisfies
\[
 J^{\pi_{\rm CBS}(\cdot;C,S)}
 >J^{\pi^{\rm dual}_1}+\frac h4(1-2\theta)\mu^2.
\]
\end{lemma}
Its proof is given in Section~\ref{sec:cdstargetscreen-new}.
\begin{lemma}\label{lem:cdscapscreen-new}
Under the same assumptions, every integer $C\ge2$ satisfies
\[
 J^{\pi_{\rm CBS}(\cdot;C,2)}
 =J^{\pi_{\rm CBS}(\cdot;2,2)}\ge J^{\varsigma_1}.
\]
\end{lemma}
Its proof is given in Section~\ref{sec:cdscapscreen-new}.

{
\begin{lemma}\label{lem:cds-positive-calendar}
Under the assumptions of Theorem~\ref{thm:cds}, let $S\ge1$ be an
integer and let $C_t\ge1$ be integer caps fixed independently of demand.
Then
\[
 J^{\pi_{\rm CBS}(\cdot;\bm C,S)}\ge J^{\varsigma_1}\quad(S=2),
 \qquad
 J^{\pi_{\rm CBS}(\cdot;\bm C,S)}>J^{\pi^{\rm dual}_1}
          +\frac h4(1-2\theta)\mu^2\quad(S\ne2).
\]
\end{lemma}
Its proof is given in Section~\ref{sec:cds-positive-calendar}.
\begin{lemma}\label{lem:cds-binary-calendar}
Under the assumptions of Theorem~\ref{thm:cds}, every cap sequence
$C_t\in\{0,1\}$ fixed independently of demand satisfies
\[
 J^{\pi_{\rm CBS}(\cdot;\bm C,S)}\ge J^{\varsigma_1}\quad(S\ge2),\qquad
 J^{\pi_{\rm CBS}(\cdot;\bm C,1)}\ge J^{\pi_{\rm CBS}(\cdot;1,1)}.
\]
\end{lemma}
Its proof is given in Section~\ref{sec:cds-binary-calendar}.
}

\medskip
\noindent\underline{\textbf{Step 1: demand bounds.}}
IFR and $T_0=1$ imply $T_{i+j}\le T_iT_j$:
the $j$ successive tail ratios starting at $i$ are no greater
than those starting at zero (zero tails are immediate).
Thus $T_j\le q^j$, and
\begin{equation}
 \frac{\mu}{1+\mu}\le q\le\mu,\qquad
 \nu_j:=\E\binom Dj
 =\sum_{i_1,\ldots,i_j\ge1}T_{i_1+\cdots+i_j}
 \le\mu^j\quad(j\ge1).
 \label{eq:cds-ifr-new}
\end{equation}
Indeed, $\mu=\sum_{j\ge1}T_j\le q/(1-q)$ gives the lower bound
on $q$. The moment identity follows from
Eq.~\eqref{eq:binomial-tail-moment}, since there are
$\binom{n-1}{j-1}$ positive integer $j$-tuples summing to $n$.

\medskip
\noindent\underline{\textbf{Step 2: compare with $\varsigma_1$.}}
Splitting on whether the first demand equals one or is at least
two gives
\[
 \psi_1=(q-T_2)q+T_2=q^2+vT_2,\qquad
 q^2\le\psi_1\le2\mu^2.
\]
Since $\mu\le1/100$ and $q\le\mu$, we have $5/v^5\le16/3$.
The holding lower bound and loss upper bound in
Eq.~\eqref{eq:cds-explicit-new}, with $C=1$, imply
\[
 \begin{aligned}
 \Gamma_{X,1}-\beta\Gamma_{A,1}
 &\ge1-\frac{35}{3}\psi_1-\beta\psi_1(1+16\mu/3)\ge(1-2\theta)-\frac{32\theta}{3}\mu-\frac{70}{3}\mu^2\ge(1-2\theta)-\frac{167}{30}\mu
 \ge\frac{1-2\theta}{2}.
 \end{aligned}
\]
The last two inequalities use $\theta<1/2$, $\mu\le1/100$ and
$\mu\le(1-2\theta)/12$, respectively.
The frequency lower bound in Eq.~\eqref{eq:cds-explicit-new} and
Eq.~\eqref{eq:cds-ifr-new} also give
\[
 \alpha_1\ge v^6\psi_1\ge
 \frac{(1-6\mu)\mu^2}{(1+\mu)^2}\ge\frac{\mu^2}{2},
\]
where $\mu\le1/100$ ensures the last inequality.
The exact stationary identities therefore yield
\[
 J^{\varsigma_1}-J^{\pi^{\rm dual}_1}
 =h\alpha_1(\Gamma_{X,1}-\beta\Gamma_{A,1})
 \ge\frac h4(1-2\theta)\mu^2.
\]

\medskip
\noindent\underline{\textbf{Step 3: cover all integer CBS parameters.}}
Lemma~\ref{lem:cdstargetscreen-new} treats $S\ne2$ and $C=0$.
At target two, $C=1$ is the benchmark just compared, and
Lemma~\ref{lem:cdscapscreen-new} treats every $C\ge2$.
Taking the infimum over all nonnegative integer parameters gives
\[
 J^{\pi^{\rm dual}_1}
 \le {\inf_{S,C\in\mathbb Z_{\ge0}}
          J^{\pi_{\rm CBS}(\cdot;C,S)}}-\frac h4(1-2\theta)\mu^2,
\]
{for the fixed-cap benchmark.}

{
\medskip
\noindent\underline{\textbf{Step 4: include every adjacent-integer cap sequence.}}
Lemma~\ref{lem:cds-positive-calendar} treats positive lower caps.
Lemma~\ref{lem:cds-binary-calendar} treats zero/one caps at targets
$S\ge2$ and reduces target one to its fixed-cap benchmark, already
excluded by Lemma~\ref{lem:cdstargetscreen-new}. Target zero never
orders. These cases exhaust Eq.~\eqref{eq:cbsbenchmark} and
retain the gap $h(1-2\theta)\mu^2/4$.
}

\subsection{Proof of Lemma~\ref{lem:cdstargetscreen-new}}
\label{sec:cdstargetscreen-new}
Assume Eq.~\eqref{eq:cdswindow} and use the notation and IFR bounds
in Section~\ref{sec:cdstheorem-new}. We first upper-bound the cost
of $\pi^{\rm dual}_1$. We then compare this bound with target-one
CBS, CBS with targets at least three, and policies that never order.
These are exactly the cases covered by the lemma.

\medskip
\noindent\underline{\textbf{Step 1: candidate loss and cost.}}
For every state of $\pi^{\rm dual}_1$,
\[
 \operatorname{IP}(s_k)+a_k+a_{k+1}+a_{k+2}\ge2.
\]
If the inventory position is two, no further argument is needed.
If it is zero, the next three orders are $(1,0,1)$ before any
arrival. If it is one, either $p_{k,3}=0$ and $a_k=1$, or
$p_{k,3}=1$ and $(a_k,a_{k+1})=(0,1)$.
These cases prove the inequality for every demand path.

All these units are available by $t=k+5$. Conservation, $V_i\le D_i$,
and nonnegative stock give
\begin{equation}
 x_t+p_{t,1}
 =\operatorname{IP}(s_{t-5})+\sum_{i=t-5}^{t-3}a_i
       -\sum_{i=t-5}^{t-1}V_i
 \ge(2-R_5)^+,\qquad R_5:=\sum_{i=t-5}^{t-1}D_i.
 \label{eq:cds-available-stock}
\end{equation}
For nonnegative integers $d,r$,
\begin{equation}
 (d-(2-r)^+)^+\le\binom d3+r\binom d2+d\binom r2.
 \label{eq:cds-positive-part}
\end{equation}
For $r=0,1,\ge2$, respectively, the first, second, or third term
bounds the loss. Eq.~\eqref{eq:cds-available-stock} gives
\begin{equation}
 A_t=(D_t-x_t-p_{t,1})^+
 \le\bigl(D_t-(2-R_5)^+\bigr)^+.
 \label{eq:cds-period-loss}
\end{equation}
Since demands are i.i.d., we have
\begin{equation}
 \E R_5=5\mu,\qquad
 \E\binom{R_5}{2}=5\nu_2+10\mu^2.
 \label{eq:cds-five-demand-moments}
\end{equation}
Using Eqs.~\eqref{eq:cds-positive-part}--\eqref{eq:cds-five-demand-moments},
independence of $D_t$ and $R_5$, and $\nu_j\le\mu^j$
from Eq.~\eqref{eq:cds-ifr-new}, we obtain in stationarity
\[
 \begin{aligned}
 \overline A^{\pi^{\rm dual}_1}
 &\le\E\bigl[\bigl(D_t-(2-R_5)^+\bigr)^+\bigr]\le\E\binom{D_t}{3}
       +\E\left[R_5\binom{D_t}{2}\right]
       +\E\left[D_t\binom{R_5}{2}\right]\\
 &=\nu_3+5\mu\nu_2+\mu(5\nu_2+10\mu^2)=\nu_3+10\mu\nu_2+10\mu^3\le21\mu^3.
 \end{aligned}
\]
Ending inventory is at most two, so
\begin{subequations}\label{eq:cds-candidatecost-new}
\begin{align}
 \overline X^{\pi^{\rm dual}_1}&\le2,\label{eq:cds-candidate-holding}\\
 \overline A^{\pi^{\rm dual}_1}&\le21\mu^3,\label{eq:cds-candidate-loss}\\
 J^{\pi^{\rm dual}_1}/h&\le2+21\theta\mu.\label{eq:cds-candidate-cost}
\end{align}
\end{subequations}

\medskip
\noindent\underline{\textbf{Step 2: target one.}}
Every positive integer cap at target one yields
$\pi_{\rm CBS}(\cdot;1,1)$.
After each sale, three periods pass before a replacement arrives;
the waiting time from its arrival through the next sale has mean
$1/q$. Each cycle has one sale and mean length $3+1/q$.
Thus
\begin{subequations}\label{eq:cds-targetone-moments}
\begin{align}
 \E V_t&=\frac{q}{1+3q},\label{eq:cds-targetone-sales}\\
 \overline X^{\pi_{\rm CBS}(\cdot;1,1)}
 &=\frac{1-q}{1+3q}\ge1-4\mu,\label{eq:cds-targetone-holding}\\
 \overline A^{\pi_{\rm CBS}(\cdot;1,1)}
 &=\mu-\frac{q}{1+3q}\ge\frac{3q^2}{1+3q}.\label{eq:cds-targetone-loss}
\end{align}
\end{subequations}
Using Eqs.~\eqref{eq:cds-candidate-cost},
\eqref{eq:cds-targetone-holding}, and~\eqref{eq:cds-targetone-loss},
$q\ge\mu/(1+\mu)$, $q\le\mu$, and $\beta\mu^2=\theta$ gives
\[
 \begin{aligned}
 \frac{J^{\pi_{\rm CBS}(\cdot;1,1)}-J^{\pi^{\rm dual}_1}}h
 &\ge-1-4\mu+
       \frac{3\theta}{(1+\mu)^2(1+3\mu)}-21\theta\mu\\
 &\ge3\theta-1-(4+39\theta)\mu
 \ge\frac{3\theta-1}{2}\ge24\mu
 >\frac{(1-2\theta)\mu^2}{4}.
 \end{aligned}
\]
Here $(1+\mu)^2(1+3\mu)\le1+6\mu$ and
$1/(1+6\mu)\ge1-6\mu$ for $\mu\le1/100$.
The last line also uses $4+39\theta<24$ and
$\mu\le(3\theta-1)/48$.

\medskip
\noindent\underline{\textbf{Step 3: targets at least three and zero ordering.}}
For $S,C\ge1$, {Lemma~\ref{lem:positive-cap-holding}}
requires only $\mu<1/2$ and a finite second moment.
Together with $\nu_2\le\mu^2$
from Eq.~\eqref{eq:cds-ifr-new}, it gives
\begin{equation}
 \overline X^{\pi_{\rm CBS}(\cdot;C,S)}
 \ge S-4\mu-\frac{\nu_2}{1-\mu}
 \ge S-4\mu-\frac{\mu^2}{1-\mu}.
 \label{eq:cds-cbs-holding}
\end{equation}
In particular, using $\theta<1/2$, $\mu\le1/100$ and
Eqs.~\eqref{eq:cds-candidate-cost} and~\eqref{eq:cds-cbs-holding},
every $S\ge3,C\ge1$ satisfies
\[
 \frac{J^{\pi_{\rm CBS}(\cdot;C,S)}-J^{\pi^{\rm dual}_1}}h
 \ge1-(4+21\theta)\mu-\frac{\mu^2}{1-\mu}
 \ge\frac{17}{20}>\frac{(1-2\theta)\mu^2}{4}.
\]
If $S=0$ or $C=0$, empty-start CBS never orders, so
Eq.~\eqref{eq:cds-candidate-cost} gives
\[
 \frac{J^{\pi_{\rm CBS}(\cdot;C,S)}-J^{\pi^{\rm dual}_1}}h
 \ge\frac{\theta}{\mu}-2-21\theta\mu
 >30>\frac{(1-2\theta)\mu^2}{4}.
\]
Together with Step 2, these cases prove the lemma.

\subsection{Proof of Lemma~\ref{lem:cdscapscreen-new}}
\label{sec:cdscapscreen-new}
Assume Eq.~\eqref{eq:cdswindow} and use
$T_2\le\mu^2$ and $5/v^5\le16/3$ from
Section~\ref{sec:cdstheorem-new}.
At target two, all caps $C\ge2$ give the same empty-start policy
$\pi_{\rm CBS}(\cdot;2,2)$.
Compare this reference with $\varsigma_1=\pi_{\rm CBS}(\cdot;1,2)$;
both preserve $\{s\in\mathbb Z_{\ge0}^4:\operatorname{IP}(s)\le2\}$.

\medskip
\noindent\underline{\textbf{Step 1: one delayed unit.}}
The policies differ only at $\bm0$.
In the paired comparison, orders at times zero and one are
$(2,0)$ and $(1,1)$, so a unit is delayed from time three to four.
The first extra loss event is $\{D_3\ge2\}$, of probability $T_2$.
Off this event the total holding saving is one and future costs
agree; on it both ending stocks at time three are zero.
Five zero demands reset the reference on the common set.
The future-block calculation in Step 2 of
Section~\ref{sec:cdsmechanism-new} bounds the expected remaining
duration by $5/v^5$ and demand by $5\mu/v^5$, independently of
$\{D_3\ge2\}$. On this event, subsequent holding differences
are at least $-2$ per period, and subsequent extra losses are
at most the remaining demand. Hence
\[
 \begin{aligned}
 \mathcal X^{\pi_{\rm CBS}(\cdot;2,2)}(\bm0,1)
 &\ge(1-T_2)-\frac{10}{v^5}T_2
 =1-\left(1+\frac{10}{v^5}\right)T_2
 \ge1-\frac{35}{3}T_2,\\
 \mathcal A^{\pi_{\rm CBS}(\cdot;2,2)}(\bm0,1)
 &\le T_2+\frac{5\mu}{v^5}T_2
 \le T_2(1+16\mu/3).
 \end{aligned}
\]
The last inequalities use $q\le\mu\le1/100$, which gives
$5/v^5\le16/3$ and $1+10/v^5\le35/3$.

\medskip
\noindent\underline{\textbf{Step 2: the stationary comparison.}}
Since the policies differ only at $\bm0$, Lemma~\ref{lem:averaging}
and the bounds in Step 1 give
\[
 \begin{aligned}
 \frac{J^{\pi_{\rm CBS}(\cdot;2,2)}-J^{\varsigma_1}}h
 &=d_{\varsigma_1}(\bm0)
 \left[\mathcal X^{\pi_{\rm CBS}(\cdot;2,2)}(\bm0,1)
       -\beta\mathcal A^{\pi_{\rm CBS}(\cdot;2,2)}(\bm0,1)\right]\\
 &\ge d_{\varsigma_1}(\bm0)
 \left[1-\frac{35}{3}T_2-\beta T_2(1+16\mu/3)\right]\ge d_{\varsigma_1}(\bm0)
 \left[1-\theta-\frac{16\theta}{3}\mu-\frac{35}{3}\mu^2\right],
 \end{aligned}
\]
where the last inequality uses $T_2\le\mu^2$ and
$\beta\mu^2=\theta$. Since $\mu\le1/100$ and $\theta<1/2$,
\[
 1-\theta-\frac{16\theta}{3}\mu-\frac{35}{3}\mu^2
 >\frac12-\frac{2}{75}-\frac7{6000}
 =\frac{2833}{6000}>0.
\]
The frequency $d_{\varsigma_1}(\bm0)$ is nonnegative, so
$J^{\pi_{\rm CBS}(\cdot;2,2)}\ge J^{\varsigma_1}$.
All caps $C\ge2$ are equivalent at target two, completing the proof.

{
\subsection{Proof of Lemma~\ref{lem:cds-positive-calendar}}
\label{sec:cds-positive-calendar}
Assume Eq.~\eqref{eq:cdswindow} and fix a deterministic positive-cap
sequence. We use $T_j,q,v,\mu,\beta$ from the beginning of
Section~\ref{sec:cdsproof-new} and $\nu_j$ from
Eq.~\eqref{eq:cds-ifr-new}. Put $z=T_2$ and $\lambda=\beta q^2$.
The IFR bounds in Eq.~\eqref{eq:cds-ifr-new} give
\begin{equation}
 \begin{gathered}
 q\le\mu\le q/v,\quad z\le q^2,\quad T_j\le q^j,\quad
 \nu_j=\E\binom Dj\le\mu^j,\quad
 3\lambda-1\ge45\mu\ge45q,\quad\lambda\le1/2.
 \end{gathered}
 \label{eq:cds-calendar-conditions}
\end{equation}
For the bound $3\lambda-1\ge45\mu$, indeed, the IFR bound $\mu\le q/(1-q)$ implies
$q\ge\mu/(1+\mu)$. Hence $\lambda=\theta\frac{q^2}{\mu^2}
\ge\frac{\theta}{(1+\mu)^2}
\ge\theta(1-2\mu)$. Since $3\theta-1\ge48\mu$ by Eq.~\eqref{eq:cdswindow} and
$\theta<1/2$, it follows that $3\lambda-1\ge(3\theta-1)-6\theta\mu
\ge48\mu-6\theta\mu\ge45\mu$. In the following, we consider the case with $S=1$, $S=2$, and $S\ge 3$.

\medskip
\noindent\underline{\textbf{Step 1: compare target two with $\varsigma_1$.}}
We show $J^{\pi_{\rm CBS}(\cdot;\bm C,2)}\ge J^{\varsigma_1}$.
Take $B=\varsigma_1$. With a positive cap sequence, the target-two
policy agrees with $B$ whenever the inventory position is positive.
At $u=(0,0,0)$ it orders either one unit, as $B$ does, or two units.
Thus we only need to check the effect of the latter order.
Recall $c,f,\rho$ from Section~\ref{sec:computing-oneorder}.
Applying Eq.~\eqref{eq:reference-bias-recursion} with
$B=\varsigma_1$ at the states needed below gives the following
equations, where $f_{yrp}=f(y,r,p)$:
\begin{equation}
\begin{gathered}
 f_{000}=0,\quad f_{001}=\rho-\beta\mu,\quad
 f_{011}=2f_{001},\quad f_{110}=3f_{001},\\
 vf_{010}=c(1)-\beta\mu+(1-3q)f_{001},\quad
 f_{100}=2(c(1)-\beta\mu)+(1-3q)f_{001},\\
 f_{101}=f_{010}+f_{001},\quad
 f_{020}=f_{200}-f_{001},\quad f_{002}=f_{200}-2f_{001},\\
 4f_{001}=c(1)-\beta\mu+vf_{200}+qf_{100},\quad
 f_{001}+qf_{200}=c(2)-\beta\mu+(q-z)f_{100}.
\end{gathered}
\label{eq:cds-calendar-core}
\end{equation}
We also use these equations in Section~\ref{sec:cds-binary-calendar}.
At $u=(0,0,0)$, orders two and one lead respectively to
$(0,0,2)$ and $(0,0,1)$. Hence
Eq.~\eqref{eq:calendar-residual} gives
$\Delta_B((0,0,0),2)=f_{002}-f_{001}=f_{200}-3f_{001}$.
Solving the displayed equations yields, with
$d=1+2q+z+3q^2+3q^2z-4qz\in[1,1.021)$,
\[
 \Delta_B((0,0,0),2)
 =\frac{1+2q+3q^2-4z-2qz-\beta vz}{d}
 \ge\frac{1-\theta-4q^2-2q^3}{d}>0.
\]
All other single-order effects vanish. The finite-horizon
decomposition in Lemma~\ref{lem:finite-horizon-decomposition}
therefore yields $J^{\pi_{\rm CBS}(\cdot;\bm C,2)}
\ge J^{\varsigma_1}$.

\medskip
\noindent\underline{\textbf{Step 2: exclude target one.}}
For $S=1$, every $C_t\ge1$ gives
$a_t=[1-\operatorname{IP}(s_t)]^+$, so
\[
 J^{\pi_{\rm CBS}(\cdot;\bm C,1)}
 =J^{\pi_{\rm CBS}(\cdot;1,1)}
 >J^{\pi^{\rm dual}_1}+\frac h4(1-2\theta)\mu^2
\]
by Lemma~\ref{lem:cdstargetscreen-new}.

\medskip
\noindent\underline{\textbf{Step 3: exclude $S\ge3$ by holding costs.}}
By Lemma~\ref{lem:positive-cap-holding} and
Eq.~\eqref{eq:cds-candidate-cost}, we have
\[
\begin{aligned}
\frac{J^{\pi_{\rm CBS}(\cdot;\bm C,S)}-J^{\pi^{\rm dual}_1}}{h}\ge (S-4\mu-Q)-(2+21\theta\mu)\ge 1-(4+21\theta)\mu-\frac{\mu^2}{1-\mu}\ge \frac{17}{20}
> \frac14(1-2\theta)\mu^2.
\end{aligned}
\]
Here the second inequality uses $S\ge3$ and
$Q\le\mu^2/(1-\mu)$; the third and final inequalities follow from $\theta<1/2$ and $\mu\le1/100$ in
Eq.~\eqref{eq:cdswindow}.

\subsection{Proof of Lemma~\ref{lem:cds-binary-calendar}}
\label{sec:cds-binary-calendar}
Assume Eq.~\eqref{eq:cdswindow}. Fix an integer target $S\ge1$
and a deterministic cap sequence $C_t\in\{0,1\}$.
Let $\pi=\pi_{\rm CBS}(\cdot;\bm C,S)$ start empty, with reduced
state $u_t=(y_t,r_t,p_t)=(x_t+p_{t,1},p_{t,2},p_{t,3})$
and order $a_t=\pi(s_t,t)$.
We use $q,v,\mu,\beta$ from the beginning of
Section~\ref{sec:cdsproof-new}, and $z=T_2$ and the bounds in
Eq.~\eqref{eq:cds-calendar-conditions} from
Section~\ref{sec:cds-positive-calendar}.

For $S\ge2$, take $B=\varsigma_1$ and use the single-order effects
from Section~\ref{sec:computing-oneorder}.
Step 1.1 bounds these effects. Although skipping an order can
have a negative effect at some states, the four-period comparison
in Step 1.2 shows that these negative effects are compensated
by positive effects, giving
$\sum_{t<H}\E\Delta_B(u_t,a_t)\ge-6$ for $H\ge4$.
Lemma~\ref{lem:finite-horizon-decomposition} then yields
$J^\pi\ge J^{\varsigma_1}$.
Finally, Step 2 treats $S=1$ by comparing $\pi$ with
$\pi_{\rm CBS}(\cdot;1,1)$.

\medskip
\noindent\underline{\textbf{Step 1.1: single-order effects for $S\ge2$.}}
We establish the bounds needed for the comparison:
\begin{equation}
\begin{gathered}
-2/q\le f(y,r,p)\le0 \quad (y+r+p\le2),\quad
\Delta_B((0,0,0),0)\ge1/(10q),\\
\Delta_B((1,0,0),0)\ge 32q,\quad \Delta_B((0,0,1),0),
\Delta_B((0,1,0),0)\ge-1,\\
\Delta_B((y,r,p),1)\ge3/(5q)\ge60
\qquad (y+r+p\ge2).
\end{gathered}
\label{eq:cds-binary-residuals}
\end{equation}
Thus only skipped orders at $(0,0,1)$ and $(0,1,0)$
may have a negative effect.

For $y+r+p\le2$, the first two lines follow by
substituting the finite-state equations in
Eq.~\eqref{eq:cds-calendar-core} into
Eq.~\eqref{eq:calendar-residual} and using
Eq.~\eqref{eq:cds-calendar-conditions}.

To prove the last line at higher states, define
\[
 l(x,p)=f(x,p,1)-f(x,p,0),\quad
 m(x)=f(x,1,0)-f(x,0,0),\quad
 d_j=f(j,0,0)-f(j-1,0,0).
\]
The reference recursion gives, successively,
\begin{align*}
  &  \Delta_B((y,r,p),1)=\E l((y-D)^++r,p),\quad 
 l(x,p)=\E m((x-D)^++p)\quad(x+p\ge2),\\
 & \qquad\qquad \qquad \qquad \qquad \qquad m(x)=\E d_{(x-D)^++1}\quad(x\ge2).
\end{align*}
The finite-state calculation also gives
$d_1\ge-2/q$ and $d_2=-3\Delta_B((0,0,1),0)\ge0$.

We next show $d_j\ge9/(10q)$ for all $j\ge3$.
For such $j$, Eq.~\eqref{eq:cds-calendar-conditions} gives
$T_j\le q^j$, $T_{j-1}\le q^{j-1}$, and
$\beta q^2\le1/2$. Assuming inductively that
$d_i\ge0$ for $2\le i<j$, the reference recursion yields
\[
 \begin{aligned}
 qd_j
 =1-(\beta+1)T_j+
   \sum_{k=1}^{j-1}\Prb(D=k)d_{j-k}
 \ge1-(\beta+1)T_j-\frac2q\Prb(D=j-1)
 \ge1-(5/2+q^2)q^{j-2}\ge9/10.
 \end{aligned}
\]
The last inequality uses $j\ge3$ and $q\le1/100$.
The low-state bound on $f$ also gives $m(0),m(1)\ge-2/q$
and $l(x,p)\ge-2/q$ when $x+p<2$.
In each of the three identities, zero demand retains the
positive bound, while a lower-state contribution under positive
demand is at least $-2/q$. Thus, in reverse order,
\[
 \begin{aligned}
 m(x)&\ge\frac{9v}{10q}-\frac2q\Prb(D\ge x)
       \ge\frac{9v}{10q}-\frac{2z}{q}
       \ge\frac4{5q}, &&x\ge2,\\
 l(x,p)&\ge\frac{4v}{5q}-\frac2q\Prb(D\ge1)
       =\frac{4v}{5q}-2
       \ge\frac7{10q}, &&x+p\ge2,\\
 \Delta_B((y,r,p),1)&\ge\frac{7v}{10q}-\frac2q\Prb(D\ge1)
       =\frac{7v}{10q}-2
       \ge\frac3{5q}\ge60, &&y+r+p\ge2.
 \end{aligned}
\]
The last inequalities use $z\le q^2$, $v=1-q$, and
$q\le1/100$.

\medskip
\noindent\underline{\textbf{Step 1.2: four-period compensation excludes every $S\ge2$.}}
We prove $\sum_{t<H}\E\Delta_B(u_t,a_t)\ge-6$ for $H\ge4$.
Let $I_t=y_t+r_t+p_t=\operatorname{IP}(s_t)$ and define $ V_t^*\coloneqq\mathbf1\{u_t\in\{(0,0,1),(0,1,0)\},a_t=0\}$, 
\[
 W_t=\left\{
\begin{array}{ll@{\qquad}ll}
 1/(10q),&(u_t,a_t)=((0,0,0),0),&
 32q,&(u_t,a_t)=((1,0,0),0),\\
 10,&I_t\ge2,\ a_t=1,
 &0,&\text{otherwise}.
\end{array}
\right.
\]
By Step 1.1, $\Delta_B(u_t,a_t)\ge W_t-V_t^*$.
Since the system starts empty and $a_t\in\{0,1\}$, every reachable
state satisfies $r_t,p_t\in\{0,1\}$ and $I_t\le S$.
For a fixed initial state with $r_0,p_0\in\{0,1\}$ and $I_0\le S$, we show
\begin{equation}
 \E\sum_{j=0}^3W_j+4V_0^*\ge8\E V_3^*.
 \label{eq:cds-fourperiod-calendar}
\end{equation}
If $C_3=1$ or $a_0=1$, then $V_3^*=0$ and this
inequality is immediate. Otherwise,
\[
 \{V_3^*=1\}\subseteq
 \{a_1+a_2=1,\ D_0+D_1+D_2\ge I_0\}.
\]

\noindent\underline{$I_0\le1$.}
For $I_0=0$, $W_0=1/(10q)\ge10$.
For $I_0=1$, $8\E V_3^*\le24q$, whereas
\[
 u_0=(1,0,0)\ \Longrightarrow\ W_0=32q;\qquad
 u_0\in\{(0,0,1),(0,1,0)\}\ \Longrightarrow\ 4V_0^*=4\ge24q.
\]

\noindent\underline{$2\le I_0<S$.}
Let $j=\min\{i\in\{1,2\}:C_i=1\}$.
If no such $j$ exists, $V_3^*=0$; otherwise
\[
 \{D_0=\cdots=D_{j-1}=0\}
 \subseteq\{I_j=I_0,\ a_j=1,\ W_j=10\}\quad\Rightarrow\quad
 \E\sum_{i=0}^3W_i\ge10v^j\ge10v^2\ge8.
\]

\noindent\underline{$I_0=S=2$.}
The common compensating event satisfies
\[
 \{D_0=D_1=0,D_2=1\}\subseteq\{u_3=(1,0,0),a_3=0\}\quad\Rightarrow\quad \E W_3\ge32qv^2(q-z)\ge32q^2v^3.
\]
Write $A=\{V_3^*=1\}$ and $c=(C_1,C_2)$.
For $c=(0,0)$, $A=\varnothing$.
For $c=(1,0)$,
\[
 A\subseteq\{D_0\ge2\}\cup\{D_0\ge1,D_1+D_2\ge1\}\quad\Rightarrow\quad 
 \Prb(A)\le z+2q^2\le3q^2.
\]
For $c=(1,1)$,
\[
 \begin{aligned}
 A\cap\{(a_1,a_2)=(1,0)\}&\subseteq\{D_0\ge1,D_2\ge1\},\\
 A\cap\{(a_1,a_2)=(0,1)\}&\subseteq
 \{D_1\ge2\}\cup\{D_1\ge1,D_2\ge1\},
 \end{aligned}
 \quad\Rightarrow\quad  \Prb(A)\le z+2q^2\le3q^2.
\]
In both cases, $8\Prb(A)\le24q^2\le\E W_3$.
For $c=(0,1)$ and initial $p_0=0$, $u_0\in\{(2,0,0),(1,1,0)\}$ and
\[
 \{D_0=1\}\subseteq\{u_1=(1,0,0),a_1=0\}\quad\Rightarrow\quad 
 \E W_1\ge32q(q-z)\ge32q^2v,
\]
\[
 A\subseteq\{D_0+D_1\ge1,\ D_0+D_1+D_2\ge2\}\quad\Rightarrow\quad 
 \Prb(A)\le2z+3q^2\le5q^2,
\]
\[
 \E(W_1+W_3)\ge32q^2(v+v^3)\ge40q^2\ge8\Prb(A).
\]
For $c=(0,1)$ and initial $p_0=1$,
\[
 A\subseteq\{D_0+D_1\ge1,D_2\ge1\}\quad\Rightarrow\quad 
 8\Prb(A)\le16q^2\le\E W_3.
\]

\noindent\underline{$I_0=S\ge3$.}
If $C_1=C_2=0$, then $V_3^*=0$. Otherwise let
$j:=\min\{i\in\{1,2\}:C_i=1\}$. Since
$y_0\ge I_0-2\ge1$,
\[
 \{D_0=1,D_1=\cdots=D_{j-1}=0\}
 \subseteq\{I_j=S-1\ge2,a_j=1,W_j=10\},
\]
\[
 \E\sum_{i=0}^3W_i\ge10(q-z)v\ge10qv^2
 \ge80q^3/v^3
 \ge8(3\nu_3+6\mu\nu_2+\mu^3)\ge8\E V_3^*.
\]
These cases prove Eq.~\eqref{eq:cds-fourperiod-calendar}.
Condition on the state at time $t$ and apply this four-period
bound to periods $t,\ldots,t+3$. Summing over
$t=0,\ldots,H-4$ gives $4\sum_{t<H}\E W_t+4\sum_{t<H}\E V_t^*
 \ge8\sum_{t=3}^{H-1}\E V_t^*$. Consequently, by the per-period bound from Step 1.1,
\[
 \sum_{t<H}\E\Delta_B(u_t,a_t)
 \ge\sum_{t<H}\E W_t-\sum_{t<H}\E V_t^*
 \ge-2\sum_{t=0}^2\E V_t^*\ge-6.
\]
Lemma~\ref{lem:finite-horizon-decomposition} now gives
$J^{\pi_{\rm CBS}(\cdot;\bm C,S)}\ge J^{\varsigma_1}$.

\medskip
\noindent\underline{\textbf{Step 2: comparison for $S=1$.}}
Take $B=\pi_{\rm CBS}(\cdot;1,1)$.
Its renewal cycle has mean length $3+1/q$, one sale, and mean
cumulative ending inventory $1/q-1$. Hence
\[
 \overline V^B=\frac{q}{1+3q},\qquad
 \overline X^B=\frac{1-q}{1+3q},\qquad
 \frac{J^B}h=1+\beta\mu-\frac{(\beta+4)q}{1+3q}.
\]
The only deviation is skipping at $(0,0,0)$, with
\[
 \Delta_B((0,0,0),0)=\beta\mu-\frac{J^B}h
 =\frac{(\beta+4)q}{1+3q}-1>0
 \qquad\left(\beta q\ge\frac1{3q},\ q\le\frac1{100}\right).
\]
Lemma~\ref{lem:finite-horizon-decomposition} gives $J^\pi\ge J^B$.
}

\section{Implementation details of the policy-discovery system}
\label{app:discoveryimplementation}
This section records the search-control and evaluation details used in Section~\ref{sec:methods}. A candidate is an SOF or an SOF template, as defined in Section~\ref{sec:sec:SOF-policy-class}. Candidate classes, evaluation instances, fitted gaps, search scores, and Pareto frontiers are defined in Section~\ref{sec:methods}. The settings below apply to the same search loop in all three discovery stages.

\subsection{{\color{black}Population} updates and adaptive parent selection}
\label{app:discoverysearchcontrol}
The search selects parents {\color{black}directly from the current population $\mathcal G_r^{\rm pop}$}, and uses adaptive scores to reward parents that generate useful offspring. {\color{black}The LLM proposes 20 edits per iteration. The offspring that pass $\mathsf{Validate}$ form $\mathcal G_r^{\rm new}$.}

\medskip
\noindent\underline{\textbf{Parent-set construction.}}
At iteration $r$, the parent set $\mathcal G^{\rm par}_r\subseteq{\color{black}\mathcal G_r^{\rm pop}}$ contains at most 24 distinct candidates and is constructed in two steps. First, the search reserves two places for candidates with the lowest search scores, two for behaviorally novel candidates, two for recent lineages, and six for representatives of the two Pareto frontiers, denoted by $\mathcal F^{\rm par}({\color{black}\mathcal G_r^{\rm pop}})$ and $\mathcal F^{\rm cmp}({\color{black}\mathcal G_r^{\rm pop}})$. The distinct candidates selected from these categories form the initial $\mathcal G^{\rm par}_r$.

Second, if $\mathcal G^{\rm par}_r$ has not reached the parent-count limit, adaptive sampling adds further parents. Let $\mathcal G^{\rm rem}_r={\color{black}\mathcal G_r^{\rm pop}}\setminus\mathcal G^{\rm par}_r$ contain the candidates not yet selected. At each draw, a candidate $\mathfrak c_{\rm sel}$ is sampled from $\mathcal G^{\rm rem}_r$ according to
\begin{equation}
 \mathfrak c_{\rm sel}\sim\Pr\nolimits^{\rm adp}(\,\cdot\mid\mathcal G^{\rm rem}_r),
 \qquad \Pr\nolimits^{\rm adp}(\mathfrak c\mid\mathcal G^{\rm rem}_r)
 =\frac{\exp(B_r(\mathfrak c)/\tau)}{\sum_{\mathfrak c'\in\mathcal G^{\rm rem}_r}\exp(B_r(\mathfrak c')/\tau)},
 \qquad \mathfrak c\in\mathcal G^{\rm rem}_r.
 \label{eq:methodislandparentsampling}
\end{equation}
After each draw, we update the sets by $\mathcal G^{\rm par}_r\leftarrow\mathcal G^{\rm par}_r\cup\{\mathfrak c_{\rm sel}\}$ and $\mathcal G^{\rm rem}_r\leftarrow\mathcal G^{\rm rem}_r\setminus\{\mathfrak c_{\rm sel}\}$. Here, $B_r(\mathfrak c)$ is an adaptive parent score that rewards useful offspring and encourages exploration of parents with little offspring evidence (see Eq.~\eqref{eq:methodparentscore} below for its detailed definition), and the temperature $\tau=0.5$ controls the concentration of the sampling probabilities. Sampling continues until the parent-count limit is reached or $\mathcal G^{\rm rem}_r$ is empty. {\color{black}Every five iterations, the reflection agent uses the selected parents $\mathcal G_r^{\rm par}$ and their records $\mathcal E_r$ to propose motifs.}

\medskip
\noindent\underline{\textbf{{\color{black}Population} updates.}}
{\color{black}The search first fits and evaluates the initial population $\mathcal G_0^{\rm pop}$ on $\mathcal T$.} Write $\mathcal F_r^k=\mathcal F^k({\color{black}\mathcal G_r^{\rm pop}})$ for the global frontiers, $k\in\{\mathrm{par},\mathrm{cmp}\}$.
Applying the LLM-proposed typed AST edits and validating the resulting candidates yields a set $\mathcal G^{\rm new}_r\subseteq\mathcal H$ of valid offspring at iteration $r$.
For each $\mathfrak c\in\mathcal G^{\rm new}_r$, parameters are fitted separately at every $\Theta_j\in\mathcal T$, and the resulting costs determine its score $f(\mathfrak c)$ through Eq.~\eqref{eq:methodsearchscore}.
{\color{black}After evaluating the offspring, we update the population by}
\begin{equation}
{\color{black}\mathcal G_{r+1}^{\rm pop}=\mathcal G_r^{\rm pop}\cup\mathcal G_r^{\rm new}.}
\label{eq:methodpopulationrecursion}
\end{equation}
The global frontiers are updated as $\mathcal F_{r+1}^k=\mathcal F^k({\color{black}\mathcal G_{r+1}^{\rm pop}})$. After $r_{\max}$ iterations, the final {\color{black}population $\mathcal G_{r_{\max}}^{\rm pop}$} determines the output frontiers.

\medskip
\noindent\underline{\textbf{Adaptive parent scores.}}
The adaptive rule measures a parent's demonstrated value for generating useful offspring. Let $\operatorname{pa}(\mathfrak c)$ be the nonempty set of recorded parents of a completed offspring $\mathfrak c$, and $f_{\rm pa}(\mathfrak c)=\min_{\mathfrak c'\in\operatorname{pa}(\mathfrak c)}f(\mathfrak c')$. {\color{black}For an offspring evaluated in iteration $r$,} the offspring reward is
\begin{equation}
\rho(\mathfrak c)=\min\left\{4,\frac{[f_{\rm pa}(\mathfrak c)-f(\mathfrak c)]^+}{|f_{\rm pa}(\mathfrak c)|+\epsilon}\right\} +\tfrac12\bigl(\mathbf 1\{\mathfrak c\in\mathcal F_{{\color{black}r+1}}^{\rm par}\}+\mathbf 1\{\mathfrak c\in\mathcal F_{{\color{black}r+1}}^{\rm cmp}\}\bigr).
\label{eq:methodoffspringreward}
\end{equation}
Here, $\epsilon>0$ is a fixed numerical stabilizer for the reward denominator. The reward is divided equally among all the recorded parents; candidates without recorded parents generate no parent credit.
Let $\mathcal C_r(\mathfrak c)$ be the completed offspring with recorded parent $\mathfrak c$. Define $n_r(\mathfrak c)=|\mathcal C_r(\mathfrak c)|$ and $u_r(\mathfrak c)=\sum_{\mathfrak c'\in\mathcal C_r(\mathfrak c)}\rho(\mathfrak c')/|\operatorname{pa}(\mathfrak c')|$. The mean reward is $\bar\rho_r(\mathfrak c)=u_r(\mathfrak c)/n_r(\mathfrak c)$ when $n_r(\mathfrak c)>0$, and zero otherwise; $n_r^{\rm tot}=\sum_{\mathfrak c\in{\color{black}\mathcal G_r^{\rm pop}}}n_r(\mathfrak c)$ counts completed parent--offspring attributions. The adaptive parent score is
\begin{equation}
B_r(\mathfrak c)=\bar\rho_r(\mathfrak c) +\beta\sqrt{\frac{\log(n_r^{\rm tot}+2)}{1+n_r(\mathfrak c)}} +\frac{\eta}{1+n_r(\mathfrak c)}.
\label{eq:methodparentscore}
\end{equation}
Here $\beta=0.8$ controls exploration and $\eta=0.5$ provides a bonus for candidates with little offspring evidence.
These scores determine the softmax probabilities in Eq.~\eqref{eq:methodislandparentsampling}, with weights renormalized over the remaining candidates at each draw.

\subsection{Symbolic edits and candidate construction}
\label{app:discoveryedits}
The operator $\mathsf{Edit}^{\rm LLM}$ in Eq.~\eqref{eq:methodpromptgeneration} maps the selected parents and proposal prompt to proposed candidate structures. For an individual candidate generated from $k$ parents $\mathfrak c_1,\ldots,\mathfrak c_k\in\mathcal G^{\rm par}_r$, we write
\begin{equation}
 \mathfrak c^{\rm new}=\mathsf{Edit}^{\rm LLM}(\mathfrak c_1,\ldots,\mathfrak c_k;\mathit{Prompt}_r),
 \qquad k\in\{1,2\}.
 \label{eq:methodtypedproposal}
\end{equation}
Here, $k=1$ denotes a single-parent edit and $k=2$ a crossover between two parents. Write each SOF parent as $\mathfrak c_i=(g_i,\bm\alpha_i)$ and each SOF-template parent as $\mathfrak c_i=(g_i,\Phi_i,\bm\alpha_i)$. An edit constructs a new ordering expression $g^{\rm new}$ and its required parameter declarations $\bm\alpha^{\rm new}$; for a template, it also constructs the required summary definitions $\Phi^{\rm new}$. Thus the offspring is $(g^{\rm new},\bm\alpha^{\rm new})$ for an SOF and $(g^{\rm new},\Phi^{\rm new},\bm\alpha^{\rm new})$ for a template. The component changes for each action are specified below.

\medskip
\noindent\underline{\textbf{Edit actions.}}
Each edit changes the symbolic ordering rule of a parent and constructs a new candidate. All edit types below apply to both SOFs and SOF templates except \texttt{replace\_summary}, which applies only to templates. Below, $\oplus$ concatenates parameter vectors, listing each shared parameter only once.
\begin{enumerate}
\item \texttt{replace\_subtree} selects a subtree of $g_1$ and replaces it with an admissible expression $u$, giving $g^{\rm new}$; all other branches remain unchanged. Every parameter or summary newly used inside $u$ must be declared independently in the offspring. \emph{SOF:} $\mathfrak c^{\rm new}=(g^{\rm new},\bm\alpha^{\rm new})$, where $\bm\alpha^{\rm new}=\bm\alpha^{\rm keep}\oplus\bm\alpha^{\rm add}$: $\bm\alpha^{\rm keep}$ retains the parameters of $\bm\alpha_1$ still used by $g^{\rm new}$, and $\bm\alpha^{\rm add}$ lists the new parameters introduced by $u$. \emph{SOF-template:} $\mathfrak c^{\rm new}=(g^{\rm new},\Phi^{\rm new},\bm\alpha^{\rm new})$, where $\Phi^{\rm new}$ contains the existing and newly supplied summary definitions required by $g^{\rm new}$. The same parameter rule applies, counting parameters in both $g^{\rm new}$ and its summary definitions.
\item \texttt{cap} bounds the complete parent expression using a parameter $C$: $g^{\rm new}=\min\{g_1,C\}$ for an upper bound and $g^{\rm new}=\max\{g_1,C\}$ for a lower bound. \emph{SOF:} $\mathfrak c^{\rm new}=(g^{\rm new},\bm\alpha^{\rm new})$, where $\bm\alpha^{\rm new}=\bm\alpha_1$ if $C$ is already declared, and $\bm\alpha^{\rm new}=(\bm\alpha_1^\top,C)^\top$ otherwise; a new $C$ is supplied with its parameter domain. \emph{SOF-template:} $\mathfrak c^{\rm new}=(g^{\rm new},\Phi_1,\bm\alpha^{\rm new})$, where $\bm\alpha^{\rm new}$ follows the same rule.
\item \texttt{crossover} combines two complete parent expressions using $g^{\rm new}=\min\{g_1,g_2\}$ or $g^{\rm new}=\max\{g_1,g_2\}$. For both SOFs and templates, $\bm\alpha^{\rm new}=\bm\alpha_1\oplus\bm\alpha_2$; for templates, $\Phi^{\rm new}=\Phi_1\cup\Phi_2$, with shared summary names having identical definitions.
\item \texttt{replace\_summary} selects exactly one summary-reference leaf in $g_1$ and replaces only that leaf with a reference to another summary, giving $g^{\rm new}$. Every other summary reference, operator, and summary definition remains unchanged. \emph{SOF-template:} $\mathfrak c^{\rm new}=(g^{\rm new},\Phi^{\rm new},\bm\alpha^{\rm new})$, where $\Phi^{\rm new}=\Phi^{\rm keep}\cup\{\varphi^{\rm replacement}\}$. Here $\Phi^{\rm keep}$ contains all parent definitions still required by the revised expression, and $\varphi^{\rm replacement}$ is the existing or newly supplied replacement summary. Also, $\bm\alpha^{\rm new}=\bm\alpha^{\rm keep}\oplus\bm\alpha^{\rm add}$, where $\bm\alpha^{\rm keep}$ retains parent parameters still used by $g^{\rm new}$ and its summary definitions, and $\bm\alpha^{\rm add}$ lists any new parameters in the replacement summary.
\item \texttt{change\_cap\_type} selects an existing upper-cap node, or a minimum node with exactly one declared cap parameter $C$, and changes its type to \texttt{integer}, \texttt{continuous}, or \texttt{xi}. The cap limit becomes $\lfloor C \rceil$, $C$, or $\xi(C,t)$, respectively; the capped expression and all other branches remain unchanged. \emph{SOF:} $\mathfrak c^{\rm new}=(g^{\rm new},\bm\alpha_1)$. \emph{SOF-template:} $\mathfrak c^{\rm new}=(g^{\rm new},\Phi_1,\bm\alpha_1)$. The parameter $C$ is fitted on an integer domain for the \texttt{integer} type and a continuous domain for the other two types, within its declared bounds.
\end{enumerate}

\medskip
\noindent\underline{\textbf{Candidate validation.}}
After an edit produces a candidate, $\mathsf{Validate}$ checks that the returned AST matches the declared edit, all parameter, summary, and motif references have valid definitions, and the candidate satisfies the type, grammar, parameter-count, and expanded-complexity limits. We allow at most eight fitted parameters, AST depth ten, 64 AST nodes, and expanded complexity 150. Parameter domains must be declared, and references must contain no recursive cycles. For templates, summary definitions must satisfy Eq.~\eqref{eq:state-summary} and produce valid expressions at every target configuration. Candidates passing these checks proceed to parameter fitting and cost evaluation.

\subsection{Screening, strict evaluation, and fitting budgets}
\label{app:discoveryevaluation}
Following Section~\ref{sec:methodfitting}, parameters are fitted separately for each candidate $\mathfrak c$ and instance $\Theta$ to minimize the numerical cost $J^{\rm eval}$ over horizon $H=1000$. Let $\mathcal A_{\mathfrak c,\Theta}$ be the admissible parameter domain. The fitting objective is
\begin{equation}
\bm{\alpha}^\star(\mathfrak c,\Theta) \in\arg\min_{\bm{\alpha}\in\mathcal A_{\mathfrak c,\Theta}} J^{\rm eval}\bigl(\pi_{g_{\bm\sigma}}(\,\cdot\,;\bm{\alpha});\Theta\bigr).
\label{eq:methodinneroptimization}
\end{equation}

Directly solving Eq.~\eqref{eq:methodinneroptimization} is difficult because $g_{\bm\sigma}$ can contain operators such as minimum, maximum, clipping, and indicators, and its output is rounded to an integer order quantity through Eq.~\eqref{eq:definition-sof}. These operations make the cost objective generally nonsmooth as a function of $\bm\alpha$. Moreover, the parameter domain may contain integer, continuous, or grid-valued coordinates.

To reduce evaluation cost, trial parameter vectors are first screened using a cheaper cost estimate. Promising vectors are then shortlisted for strict evaluation, which supplies $J^{\rm eval}$ and determines the fitted vector. Under state propagation, the two tiers differ in whether the remaining horizon is approximated; under Monte Carlo simulation, they differ in path count.

When the reachable state space permits exact propagation, screening may use a stationary approximation after the state distribution and expected one-period cost have stabilized. This acceleration applies only to policies whose actions depend on the inventory state and fixed instance parameters, with no explicit dependence on time or the remaining horizon; policies using $\xi(C,t)$ instead use full-horizon Monte Carlo simulation. For a trial vector $\bm\alpha\in\mathcal A_{\mathfrak c,\Theta}$, abbreviate the induced policy as $\pi=\pi_{g_{\bm\sigma}}(\,\cdot\,;\bm\alpha)$. Using the system state $s_t\in\mathcal S_{\bm\sigma}$ from Eq.~\eqref{eq:system-state}, define $d_t(s)=\Pr_{\pi,\Theta}(s_t=s\mid s_0=\bm0)$ for a state realization $s\in\mathcal S_{\bm\sigma}$. Starting from $d_0(s)=1\{s=\bm0\}$, the state distribution is propagated by
\begin{equation}
 d_{t+1}(s')=\sum_{s\in\mathcal S_{\bm\sigma}}d_t(s)\Pr\nolimits_\Theta(s_{t+1}=s'\mid s_t=s,\ a_t=\pi(s)).
 \label{eq:methodscreenpropagation}
\end{equation}
Given $d_t$, the expected one-period cost is
\begin{equation}
 c_t=\sum_{s\in\mathcal S_{\bm\sigma}}d_t(s)\mathbb E_\Theta[G_t\mid s_t=s,\ a_t=\pi(s)].
 \label{eq:methodscreenperiodcost}
\end{equation}
If the state distribution approaches a stationary regime, changes in both the distribution and the expected one-period cost diminish, motivating an approximation of the remaining horizon cost. We monitor these changes and stop propagation early once $\|d_{t+1}-d_t\|_\infty\le\delta_{\rm dist}$ and $|c_t-c_{t-1}|\le\delta_{\rm cost}$ hold for three consecutive periods; otherwise, propagation continues through the full horizon. We use $\delta_{\rm dist}=\delta_{\rm cost}=10^{-10}$ for these absolute checks; evaluators using cost-scaled checks instead use coefficient $10^{-12}$ for the cost criterion. Let $T\le H$ be the number of computed periods, with costs $c_0,\ldots,c_{T-1}$, and $\|d_{t+1}-d_t\|_\infty=\sup_{s\in\mathcal S_{\bm\sigma}}|d_{t+1}(s)-d_t(s)|$. The screening estimate is then
\begin{equation}
 J^{\rm scr}(\pi;\Theta)=\frac1H\left(\sum_{t=0}^{T-1}c_t+(H-T)c_{T-1}\right).
\label{eq:methodscreenextrapolation}
\end{equation}
The first term contains the $T$ computed transient costs, and the second extrapolates the last computed one-period cost over the remaining $H-T$ periods. If $T=H$, the extrapolation vanishes and the screening cost equals $J_H$. Strict evaluation propagates shortlisted vectors through the full horizon without this stationary approximation, giving $J^{\rm eval}=J_H$.

When the reachable state space is too large for exact propagation, screening uses a smaller fixed Monte Carlo path bank, and shortlisted vectors are reevaluated using a larger strict path bank. Both tiers simulate the full horizon. For $R$ paths, the cost estimator is
\begin{equation}
 \widehat J_{H,R}(\pi;\Theta)
 =\frac{1}{R}\sum_{n=1}^{R}\frac{1}{H}
   \sum_{t=0}^{H-1}G_t^{(n)}(\pi;\Theta).
\label{eq:methodmc-eval}
\end{equation}
Here, $G_t^{(n)}(\pi;\Theta)$ is the period-$t$ cost on path $n$. Strict evaluation supplies $J^{\rm eval}(\pi;\Theta)=\widehat J_{H,R}(\pi;\Theta)$ using the larger path bank.

For either evaluation method, let $\mathcal A^{\rm strict}_{\mathfrak c,\Theta}\subseteq\mathcal A_{\mathfrak c,\Theta}$ be the nonempty set of parameter vectors that have completed strict evaluation. The fitted vector is
\begin{equation}
 \widehat{\bm\alpha}(\mathfrak c,\Theta)
 \in\arg\min_{\bm\alpha\in\mathcal A^{\rm strict}_{\mathfrak c,\Theta}}
 J^{\rm eval}\bigl(\pi_{g_{\bm\sigma}}(\,\cdot\,;\bm\alpha);\Theta\bigr).
 \label{eq:methodtwotierfit}
\end{equation}

Small finite parameter domains are enumerated. For larger domains, Sobol points provide broad coverage, differential evolution searches without gradients, and deterministic coordinate refinement improves promising vectors.

The fitting budget limits the number of parameter vectors evaluated for each candidate--instance pair. Standard, used during discovery, permits at most 5,000 screening evaluations and 96 strict evaluations. After discovery, selected structures are independently refitted under Certify for benchmark evaluation; its base allocation permits 20,000 screening evaluations and 256 strict evaluations. Above three fitted parameters, each additional parameter adds 8,000 screening evaluations and 48 strict evaluations, with the increment capped at ten parameters. This extension is unnecessary for small discrete domains that are exhaustively enumerated. The fitted vector is the best found under the applicable budget; global optimality is not guaranteed.
For Monte Carlo fitting, the path count separately controls the fidelity of each evaluation. Standard uses 128 paths for screening and 1,000 for strict evaluation, whereas Certify uses 512 and 8,000, respectively.
For final benchmark reporting, fitted parameters are evaluated on demand paths generated independently of those used for fitting. Policies evaluated on the same instance share the final evaluation paths to reduce simulation noise in their comparisons.

\subsection{Diagnostic-bank construction and record selection}
\label{app:discoverydiagnostics}
Given an instance $\Theta\in\mathcal T$ with a precomputed dynamic-programming (DP) reference, we fix a diagnostic bank $\mathcal B_\Theta=\{(t_i,s_i,\omega_i)\}_{i=1}^{N_\Theta}$, where $N_\Theta$ is the number of diagnostic records, $t_i\in\{0,\ldots,H-1\}$ is the period of record $i$, $s_i\in\mathcal S_{\bm\sigma}$ is its state, and $\omega_i$ is its probability weight. The records are obtained by propagating the DP-optimal reference policy from $s_0=\bm0$, collecting reachable states at designated periods, and pruning states with negligible probability. The weights are normalized after pruning so that $\sum_i\omega_i=1$, and the bank is held fixed across candidates. When no DP reference is available, the decision regret is not computed or supplied to the LLM.

The bank-generation horizon is $H=1000$. We inspect 16 periods evenly spaced over the middle 60\% of the horizon. Starting with unit probability at $s_0=\bm0$, we propagate the state-occupancy distribution under the DP-optimal action. After each transition, a state is discarded when its occupancy is below $\varepsilon_{\rm occ}=10^{-13}$; the remaining occupancies are renormalized before the next transition. Each retained state at an inspected period contributes one record. After pooling the records, their weights are renormalized so that $\sum_i\omega_i=1$. The same bank is used for every candidate evaluated on $\Theta$, so diagnostic rankings compare candidates on identical states and weights.

Let $q_\Theta^{\rm diag}$ be the maximum integer order quantity permitted in the DP reference computation. Define $V_{\Theta,t}^*(s)$ as the minimum expected cumulative cost from state $s$ at period $t$ through period $H-1$, subject to this action cap, with terminal value $V_{\Theta,H}^*(s)=0$. The optimal action value is
\begin{equation}
 Q_{\Theta,t}^*(s,a)=\mathbb E_\Theta\!\left[G_t+V_{\Theta,t+1}^*(s_{t+1})\mid s_t=s,\ a_t=a\right],
 \qquad a\in\{0,\ldots,q_\Theta^{\rm diag}\}.
 \label{eq:methodreferenceactionvalue}
\end{equation}
For $t=H-1,\ldots,0$, the Bellman recursion is $V_{\Theta,t}^*(s)=\min_{a\in\{0,\ldots,q_\Theta^{\rm diag}\}}Q_{\Theta,t}^*(s,a)$. An action is DP-optimal at $(t,s)$ when it attains this minimum. Thus $V_{\Theta,t}^*(s)$ is the optimal remaining total cost, whereas $Q_{\Theta,t}^*(s,a)$ is the remaining cost of taking action $a$ now and acting optimally thereafter. 
The first diagnostic asks how much a candidate's action increases the remaining expected cost relative to a DP-optimal action. Using the fitted policy $\pi_{g_{\bm\sigma}}(\,\cdot\,;\widehat{\bm\alpha}(\mathfrak c,\Theta))$ from Section~\ref{sec:methodfitting}, the probability-weighted Bellman regret for a candidate respecting the DP action cap is
\begin{equation}
 \operatorname{Reg}_{\Theta,i}(\mathfrak c)=\omega_i
 \left[
 Q_{\Theta,t_i}^*(s_i,\pi_{g_{\bm\sigma}}(s_i,t_i;\widehat{\bm\alpha}(\mathfrak c,\Theta)))-
 \min_{a\in\{0,\ldots,q_\Theta^{\rm diag}\}}Q_{\Theta,t_i}^*(s_i,a)
 \right].
 \label{eq:methodstatewiseregret}
\end{equation}

This value is nonnegative and is zero for any DP-optimal action, including ties. For each parent $\mathfrak c\in\mathcal G_r^{\rm par}$, we rank all records in $\mathcal B_\Theta$ by $\operatorname{Reg}_{\Theta,i}(\mathfrak c)$ and retain the first $K_{\rm reg}=10$ records. Each retained record contains the instance, state, period, probability weight, parent action, set of DP-optimal actions, and Bellman cost difference. These records show where a structural change could correct a consequential ordering decision.

The second diagnostic identifies states where the parents $\mathcal G^{\rm par}_r$ selected at iteration $r$ prescribe different orders. The probability-weighted parent-action spread is
\begin{equation}
 \operatorname{Div}_{\Theta,i}(\mathcal G^{\rm par}_r)=\omega_i
 \left[\max_{\mathfrak c\in\mathcal G^{\rm par}_r}\pi_{g_{\bm\sigma}}(s_i,t_i;\widehat{\bm\alpha}(\mathfrak c,\Theta))-
 \min_{\mathfrak c\in\mathcal G^{\rm par}_r}\pi_{g_{\bm\sigma}}(s_i,t_i;\widehat{\bm\alpha}(\mathfrak c,\Theta))\right].
 \label{eq:methodactionspread}
\end{equation}

The maximum and minimum in Eq.~\eqref{eq:methodactionspread} are the largest and smallest order quantities prescribed by the selected parents at state $s_i$ in period $t_i$. Their difference measures the range of ordering decisions; multiplying by $\omega_i$ gives greater emphasis to records with higher probability weights. We rank the records in $\mathcal B_\Theta$ by $\operatorname{Div}_{\Theta,i}(\mathcal G_r^{\rm par})$ and retain the first $K_{\rm div}=10$ records. {\color{black}The parents’ evaluation records, the selected diagnostic records, and the corresponding parent actions together form $\mathcal E_r$.} This selection identifies records with the greatest probability-weighted differences between parent decisions. 

\subsection{Structural complexity}
\label{app:discoverycomplexity}
The complexity Pareto frontier in Eq.~\eqref{eq:methodbilevel} uses the expanded AST complexity as its structural-simplicity measure. Each motif reference is first replaced by its registered parameter-free AST. Nested motif references are expanded recursively until only operators and terminals from the admissible grammar remain. Type and grammar validation also use this fully expanded representation.

For the resulting AST, let $n_{\rm node}$ be the number of nodes, $d_{\rm AST}$ its maximum depth, $n_{\rm par}$ the number of distinct fitted parameter symbols, and $n_{\rm nl}$ the number of nonlinear operators (minimum, maximum, absolute value, indicator, guarded ratio, or clipping). The reported expanded complexity is
\begin{equation}
\zeta_{\rm cmp}(\mathfrak c)=n_{\rm node}+d_{\rm AST}+2n_{\rm par}+1.5n_{\rm nl}.
\label{eq:methodexpandedcomplexity}
\end{equation}
This weighted measure penalizes expression size and nesting through the node and depth counts, while assigning additional weight to fitted parameters and nonlinear operators, which add tuning burden and structural complexity. The unexpanded expression uses the same definition before motif replacement; parameter count is reported separately as $\zeta_{\rm par}(\mathfrak c)=n_{\rm par}$.
For complexity accounting, $\xi(C,t)$ is treated in the same way as an ordinary fitted constant: a single fitted-parameter leaf, whose standalone complexity is $1+1+2=4$.

\subsection{Specification of LLM inputs and outputs}
\label{app:discoveryprompts}
The four LLM calls below propose candidate edits, extract reusable motifs, or initialize templates. Each prompt supplies the relevant system configurations and state-coordinate meanings, available demand and cost features, and the admissible grammar and structural limits. The grammar follows Section~\ref{sec:sec:SOF-policy-class}; edit actions and validation follow Section~\ref{app:discoveryedits}, and diagnostic records follow Section~\ref{app:discoverydiagnostics}. All responses contain JSON only. The following specification records the supplied information and required outputs without repeating the grammar or individual evaluation records.

\medskip
\noindent\textbf{SOF proposal.}
At a fixed configuration, the input contains the {\color{black}selected SOF parents $\mathcal G_r^{\rm par}$}, represented by their ordering expressions and ASTs with parameter declarations and domains, together with {\color{black}their evaluation and diagnostic records $\mathcal E_r$ and the current motif library $\mathcal M_r$}. The task is to identify structural changes that reduce search scores while keeping the ordering expressions compact. The output contains 20 edit proposals under the actions in Section~\ref{app:discoveryedits}. Direct state-coordinate leaves are permitted.

\medskip
\noindent\textbf{SOF-template proposal.}
The input contains {\color{black}the selected parent templates $\mathcal G_r^{\rm par}$}, represented by their ordering expressions and ASTs, parameter declarations and domains, and summary definitions. It also supplies {\color{black}their records $\mathcal E_r$, including overall and configuration-specific scores and available decision-regret and parent-action-spread records. The prompt specifies the target configurations, summary registry and definition schema, and motif library $\mathcal M_r$.} The output contains 20 edit proposals to ordering expressions or summary references under the actions in Section~\ref{app:discoveryedits}. State access in the ordering AST and motifs is through declared summaries; direct state-coordinate leaves are forbidden. Any new summary must be defined according to Eq.~\eqref{eq:state-summary}.

\medskip
\noindent\textbf{Reflection.}
{\color{black}The reflection agent receives the same parents $\mathcal G_r^{\rm par}$ and records $\mathcal E_r$ as the proposal agent. The prompt includes the parents’ ordering expressions, ASTs and parameter declarations, the target configurations, and the current motif library $\mathcal M_r$. For template parents, it also includes their summary definitions and the summary registry.} The task is to identify useful symbolic building blocks not already in the library. The output contains at most two new motifs, each with a name, complete formula and AST, and a brief rationale. {\color{black}Each motif must be a parameter-free subexpression of a selected parent and satisfy the grammar for its candidate class. For templates, a motif accesses the state only through existing summaries. Validated motifs enter $\mathcal M_{r+1}$.}

\medskip
\noindent\textbf{Template initialization.}
The input contains {\color{black}the AST of a seed SOF from the Stage-II frontiers $\mathcal F_{\rm II}^{\rm par}\cup\mathcal F_{\rm II}^{\rm cmp}$}, parameter declarations and domains, existing summary and motif definitions, together with the target configurations, domain vocabulary, and summary registry and schema. The task is to propose five or six distinct SOF templates by generalizing state summaries while preserving the seed's ordering logic. The output contains a \texttt{families} array; each entry identifies its parent and supplies a complete AST, parameter declarations, and any required summary or motif definitions. Every template must be valid across the target configurations; branches on individual configurations are forbidden. {\color{black}The validated template triples $(g,\Phi,\bm\alpha)$ form the initial Stage-III population $\mathcal G_{{\rm III},0}^{\rm pop}$.}

\medskip
\noindent\textbf{Proposal response fields.}
Both proposal calls return a \texttt{proposals} array. Each entry contains \texttt{parent}, \texttt{action}, the complete offspring \texttt{formula} and \texttt{ast}, and \texttt{rationale}; crossover additionally requires \texttt{other\_parent} and specifies whether the combination uses minimum or maximum. The response also supplies parameter declarations and domains for any new parameters and, for templates, definitions of any new summaries.

\section{Codex baseline formulas}
\label{app:codexbaseline}
Table~\ref{tab:codexbaseline} lists the SOF templates generated using Codex (\texttt{gpt-6-astra}), where we define the summaries $\varphi_{\rm stock}=\sum_{m=1}^{\sigma^{\rm life,\natural}}x_m$, $\varphi_{\rm pipe}=\sum_{r=1}^{\sigma^{\rm lead}}p_r$, and $\varphi_{\rm back}=\sum_{\kappa=1}^{\sigma^{\rm back}}b_\kappa$. Here $p_1$ is the receipt due in the current period and $\mu$ is the benchmark demand mean. For lead time--perishability with $1\le\sigma^{\rm life}\le10$, define the linear summaries $\psi_j=\sum_{m=1}^{\sigma^{\rm life}}\mathbf1\{m=j\}x_m$, $j=1,\ldots,10$, and $E_{10}=\left[\psi_1-\mu+\left[\psi_2-\mu+\cdots+\left[\psi_{10}-\mu\right]^+\cdots\right]^+\right]^+$.
The fixed ten-level expression $E_{10}$ measures stock exceeding mean demand before successive expiry deadlines, with $\psi_j=0$ for $j>\sigma^{\rm life}$. Each expression is converted to an order using the same rounding and action bounds as the other policies.
\begin{table}[!htbp]
\TABLE{Codex baseline expressions before order rounding and clipping.\label{tab:codexbaseline}}{%
\small
\setlength{\tabcolsep}{4pt}
\begin{tabular}{@{}l l r r@{}}
\hline
Setting & SOF template generated by Codex & Mean complexity & \# of param.'s\\
\hline
Lead time & $\min\{c,a-\varphi_{\rm stock}-(1-v)\varphi_{\rm pipe}\}$ & $36.00$ & 3\\
Perishable & $a-\varphi_{\rm stock}+v[\mu-\varphi_{\rm stock}]^+$ & $41.50$ & 2\\
Lead time--perishability & $\min\{c,a-\varphi_{\rm stock}-\varphi_{\rm pipe}+vE_{10}\}$ & $124.17$ & 3\\
Lead time--backlog & $a-\varphi_{\rm stock}-\varphi_{\rm pipe}+\varphi_{\rm back}+v[\mu+\varphi_{\rm back}-\varphi_{\rm stock}-p_1]^+$ & $45.50$ & 2\\
\hline
\end{tabular}}{}
\end{table}

\section{Pareto-frontier policies}
\label{app:paretofrontiers}

\subsection{Lead-time setting}
\label{app:leadtimefrontiers}
Table~\ref{tab:leadtimefrontiers} lists the union of the complexity and parameter-count evaluation Pareto frontiers in Figure~\ref{fig:resultspareto}. Superscripts C and P indicate membership in the complexity and parameter-count frontiers, respectively. The benchmark configurations are those in Section~\ref{sec:resultsleadtime}.
The table uses $\varphi_{\rm IP}=x_1+\sum_{r=1}^{\sigma^{\rm lead}}p_r$, the weighted tail summary $\varphi_{\rm tail}$ defined with Eq.~\eqref{eq:resultsimbalance}, and
$\varphi_{\rm head}=\sum_{r=1}^{\lceil\sigma^{\rm lead}/2\rceil}p_r$, $
 \varphi_{\rm tail5}=\sum_{r=\sigma^{\rm lead}-\lceil\sigma^{\rm lead}/5\rceil+1}^{\sigma^{\rm lead}}p_r$.
F02 is Constant Order (CO); F08 is Capped Base Stock (CBS) after nonnegative clipping. F12 is the Base Stock benchmark comparator.

\begin{table}[!htbp]
\TABLE{Lead-time policies on the evaluation Pareto frontiers.\label{tab:leadtimefrontiers}}{%
\small
\setlength{\tabcolsep}{3pt}
\renewcommand{\arraystretch}{1.10}
\begin{tabular}{@{}l l r r r@{}}
\hline
Policy & SOF template & Complexity & \# of param.'s & $f^{\rm bench}(\%)$\\
\hline
F01$^{\mathrm{C,P}}$ & $\Pr(D=0)$ & $2.00$ & $0$ & $609.5837$\\[1pt]
F02$^{\mathrm C}$ (CO) & $\xi(C,t)$ & $4.00$ & $1$ & $16.7696$\\[1pt]
F03$^{\mathrm C}$ & $\min\{r,\xi(C,t)\}$ & $10.50$ & $2$ & $16.7675$\\[1pt]
F04$^{\mathrm C}$ & $\min\{S-x_1,\xi(C,t)\}$ & $15.50$ & $2$ & $7.5993$\\[1pt]
F05$^{\mathrm C}$ & $\min\{S-x_1-p_1,\xi(C,t)\}$ & $18.50$ & $2$ & $5.7666$\\[1pt]
F06$^{\mathrm C}$ & $\min\{S-x_1-p_1-p_2,\xi(C,t)\}$ & $19.50$ & $2$ & $4.3499$\\[1pt]
F07$^{\mathrm C}$ & $\min\{S-x_1-\varphi_{\rm head},\xi(C,t)\}$ & $22.00$ & $2$ & $4.2904$\\[1pt]
F08$^{\mathrm C}$ (CBS) & $\min\{S-\varphi_{\rm IP},\xi(C,t)\}$ & $26.00$ & $2$ & $0.7191$\\[1pt]
F09$^{\mathrm C}$ (SCDS) & $\min\{K\min\{S-\varphi_{\rm IP},F-p_{\sigma^{\rm lead}}\},\xi(C,t)\}$ & $40.50$ & $4$ & $0.4367$\\[1pt]
F10$^{\mathrm C}$ & $\min\{K(S-\varphi_{\rm IP}),S_1-\varphi_{\rm tail5},\xi(C,t)\}$ & $42.29$ & $4$ & $0.4289$\\[1pt]
F11$^{\mathrm{C,P}}$ (IA-SCDS) & $\min\{\xi(C,t),K\min\{S-\varphi_{\rm IP},F-\varphi_{\rm tail}\}\}$ & $62.43$ & $4$ & $0.3086$\\[1pt]
F12$^{\mathrm P}$ (Base Stock) & $S-\varphi_{\rm IP}$ & $19.50$ & $1$ & $16.2352$\\[1pt]
F13$^{\mathrm P}$ & $\min\{K(S-\varphi_{\rm tail}),S+\mu-\varphi_{\rm IP}\}$ & $45.29$ & $2$ & $0.6465$\\[1pt]
\hline
\end{tabular}}{}
\end{table}

\subsection{Perishable-inventory setting}
\label{app:perishablefrontiers}
\begin{figure}[!htbp]
\centering
\includegraphics[width=\textwidth]{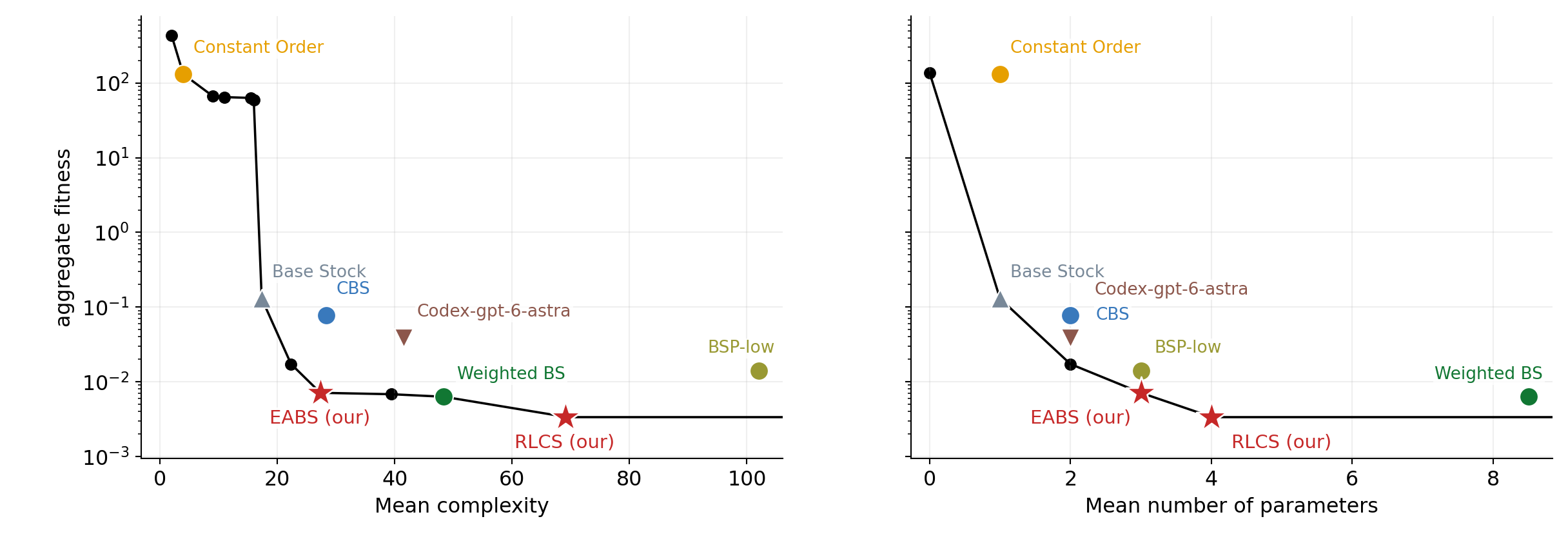}
\caption{Evaluation Pareto frontiers of the framework's final output candidates, SOF literature baselines, and Codex in the perishable-inventory setting, with EABS and RLCS highlighted. The horizontal axes report mean expanded AST complexity (left) and mean parameter count (right); both panels report aggregate benchmark fitness across $\sigma^{\rm life}=2,\ldots,15$.}
\label{fig:resultsperishablepareto}
\end{figure}

Figure~\ref{fig:resultsperishablepareto} shows the evaluation Pareto frontiers obtained by rescoring the framework's final output candidates, SOF literature baselines, and Codex on benchmark instances with $\sigma^{\rm life}=2,\ldots,15$. RLCS attains the lowest aggregate fitness on both frontiers, while EABS provides a simpler three-parameter alternative on both frontiers. Base Stock is also on both frontiers. CBS and BSP-low are dominated by discovered SOF templates in both panels. Weighted Base Stock belongs to the complexity frontier, offering lower fitness than EABS with a larger expression, but is dominated by RLCS in parameter count and fitness. A simpler two-parameter frontier template improves on CBS.

Table~\ref{tab:perishablecomplexityfrontier} lists the union of the complexity and parameter-count evaluation Pareto frontiers. Superscripts C and P indicate membership in the respective frontiers. The notation follows Section~\ref{sec:resultsperishable}.

\begin{table}[!htbp]
\TABLE{Perishable-inventory policies on the evaluation Pareto frontiers.\label{tab:perishablecomplexityfrontier}}{%
\small
\setlength{\tabcolsep}{3pt}
\renewcommand{\arraystretch}{1.10}
\begin{tabular}{@{}l l r r r@{}}
\hline
Policy & SOF template & Complexity & \# of param.'s & $f^{\rm bench}(\%)$\\
\hline
F14$^{\mathrm C}$ & $\mu$ & $2.00$ & $0$ & $435.4519$\\[1pt]
F15$^{\mathrm C}$ (CO) & $\xi(C,t)$ & $4.00$ & $1$ & $131.3412$\\[1pt]
F16$^{\mathrm C}$ & $S-x_{\sigma^{\rm life}-1}$ & $9.00$ & $1$ & $66.8744$\\[1pt]
F17$^{\mathrm C}$ & $S-x_1-x_{\sigma^{\rm life}-1}$ & $11.00$ & $1$ & $64.6904$\\[1pt]
F18$^{\mathrm C}$ & $S-x_1-\min\{x_{\sigma^{\rm life}-1},\mu\}$ & $15.50$ & $1$ & $63.0717$\\[1pt]
F19$^{\mathrm C}$ & $S-x_1-\omega_2 x_{\sigma^{\rm life}-1}$ & $16.00$ & $2$ & $59.7611$\\[1pt]
F20$^{\mathrm{C,P}}$ (Base Stock) & $S-\varphi_{\rm IP}$ & $17.36$ & $1$ & $0.1299$\\[1pt]
F21$^{\mathrm{C,P}}$ & $S-\omega_1 \varphi_{\rm IP}$ & $22.36$ & $2$ & $0.0172$\\[1pt]
F22$^{\mathrm{C,P}}$ (EABS) & $S-\omega_1 \varphi_{\rm IP}+\omega_2 x_1$ & $27.36$ & $3$ & $0.0071$\\[1pt]
F23$^{\mathrm C}$ & $\min\{\xi(C,t),S-\omega_1 \varphi_{\rm IP}-\omega_2\min\{x_{\sigma^{\rm life}-1},\mu\}\}$ & $39.43$ & $4$ & $0.0068$\\[1pt]
WBS$^{\mathrm C}$ & $S-\sum_{m=1}^{\sigma^{\rm life}-1}\omega_m x_m$ & $48.36$ & $8.50$ & $0.0063$\\[1pt]
F24$^{\mathrm{C,P}}$ (RLCS) & $\min\{\xi(C,t),S-\omega_1 \varphi_{\rm IP}-\gamma \varphi_{\rm life}/(1+\varphi_{\rm IP})\}$ & $69.14$ & $4$ & $0.0034$\\[1pt]
F25$^{\mathrm P}$ & $\mu-\varphi_{\rm IP}+\sigma$ & $16.36$ & $0$ & $137.7972$\\[1pt]
\hline
\end{tabular}}{}
\end{table}

\section{Additional inventory settings}
\label{app:additionalsettings}
We apply the framework in Section~\ref{sec:methods} separately to the two joint settings, using the training and benchmark instances in Sections~\ref{app:searchinstances} and~\ref{app:demandcatalog}.

\subsection{Lead-time--perishability setting}
\label{sec:resultsltperishable}
We consider perishable inventory with a positive replenishment lead time and lost sales, so that $\bm\sigma=(\sigma^{\rm lead},\sigma^{\rm life},0)$. We compare the discovered templates against 15 baselines: eight SOF policies and six non-SOF policies from the literature, together with the SOF template found by Codex. The SOF baselines are Constant Order, Base Stock, COP$(s,q)$, Capped Base Stock, and BSP-low \citep{haijema2016stock}, together with Weighted Base Stock, BSP-2D, and BSP-OIR \citep{haijema2019improved}. The non-SOF baselines are EWA \citep{broekmeulen2009heuristic}; BSP-low-EW, SQmax-EW, and BSP-WS3 \citep{haijema2019improved}; and PIL and APIL \citep{bu2026perishable}.

\medskip
\noindent\underline{\textbf{Representative policy templates.}} We select two representative templates: \emph{scaled capped base stock} (SCBS) and \emph{expiry- and tail-adjusted inventory shortfall} (ETIS). Define $\varphi_{\rm IP}=\sum_{m=1}^{\sigma^{\rm life}}x_m+\sum_{r=1}^{\sigma^{\rm lead}}p_r$, and let $\varphi_{\rm tail}$ be the weighted tail summary in Section~\ref{sec:resultsleadtime}. Their expressions are
\begin{align}
 q_{\rm SCBS}=\min\!\left\{\xi(C,t),K(S-\varphi_{\rm IP})\right\},\qquad 
 q_{\rm ETIS}=\min\!\left\{\xi(C,t),S-\omega_1\varphi_{\rm IP}+\omega_2x_1-\omega_3\varphi_{\rm tail}\right\}.
 \label{eq:resultsltpetis}
\end{align}
Both templates cap an inventory-based replenishment signal at $\xi(C,t)$. SCBS scales the inventory-position shortfall by $K$ and has three parameters $(S,K,C)$. ETIS adds separate weights for total inventory position, expiring stock, and late pipeline supply, giving five parameters $(S,\omega_1,\omega_2,\omega_3,C)$, with $\omega_1,\omega_2,\omega_3\ge0$. In both templates, $S$ is a fitted target and $C$ controls the cap.

\begin{figure}[!htbp]
\centering
\begin{minipage}[b]{0.51\linewidth}
\centering
\includegraphics[width=\linewidth]{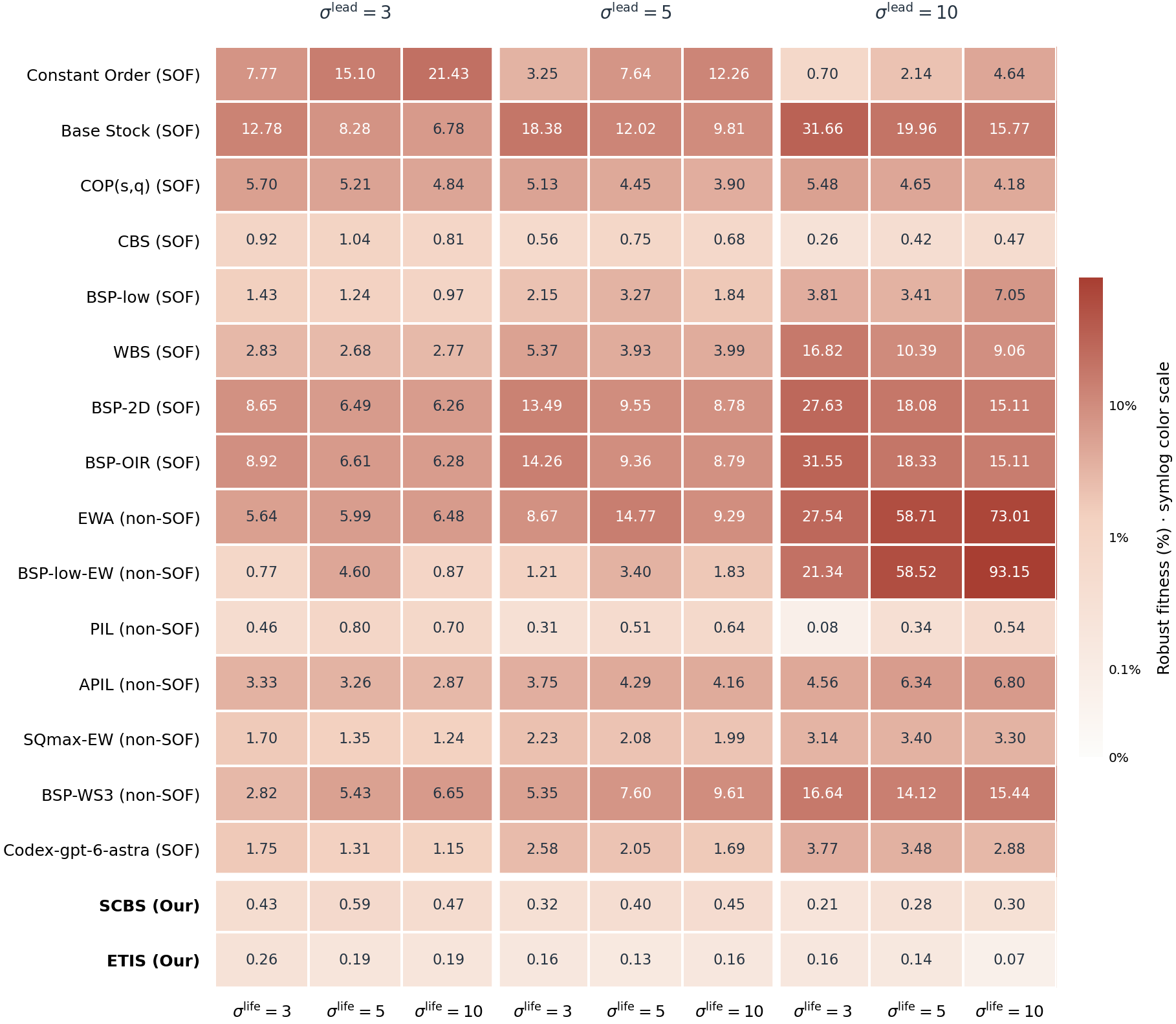}\par
\vspace{2pt}
{\small (a) Configuration-level fitness.\par}
\end{minipage}\hfill%
\begin{minipage}[b]{0.48\linewidth}
\centering
\includegraphics[width=\linewidth]{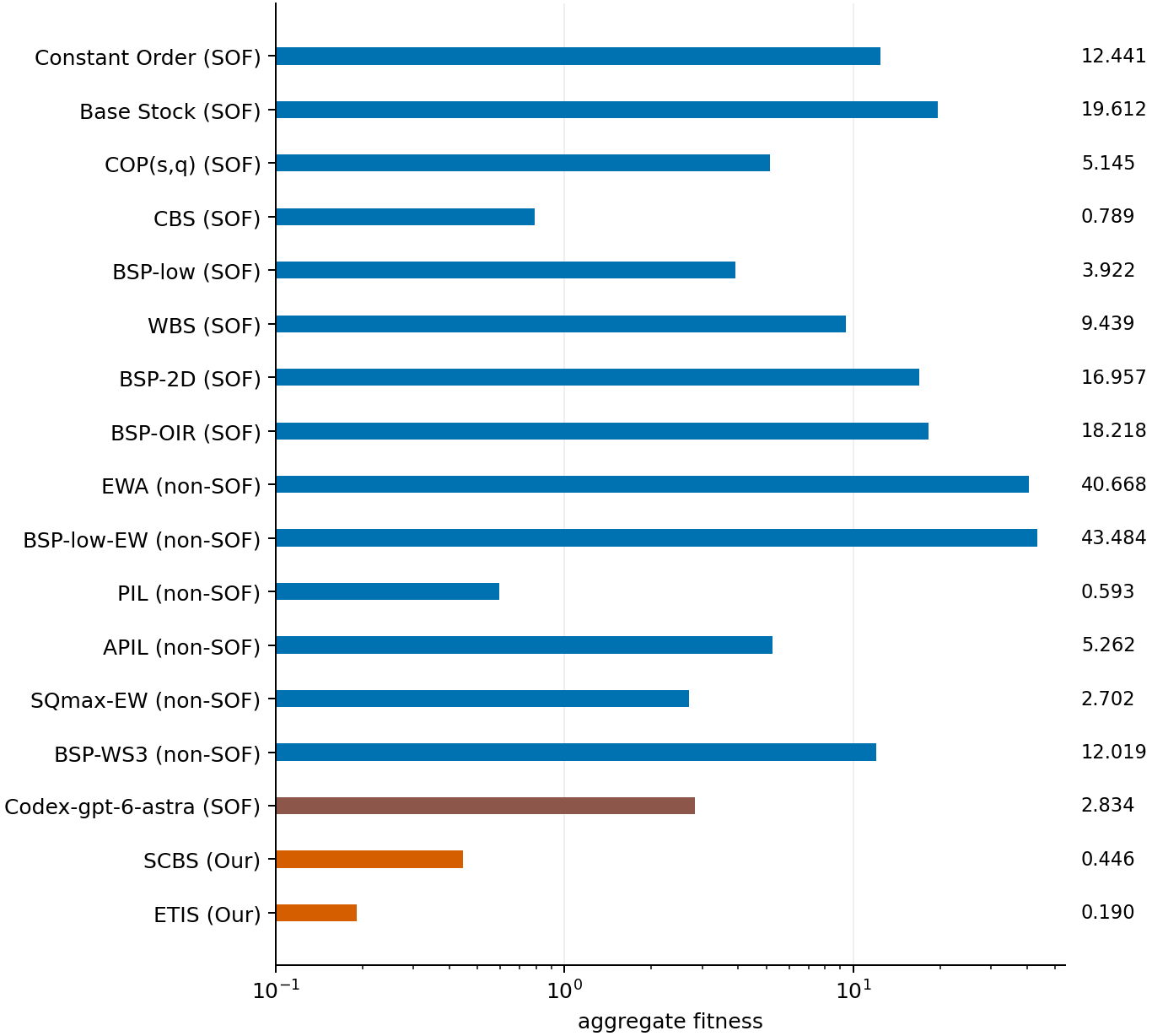}\par
\vspace{2pt}
{\small (b) Aggregate fitness.\par}
\end{minipage}
\caption{Benchmark fitness for lead time and perishability. Panel (a) reports configuration-level fitness; panel (b) reports aggregate fitness. Fitness values are percentages; lower is better.}
\label{fig:resultsltpmc}
\end{figure}

\medskip
\noindent\underline{\textbf{Evaluation of the representative templates.}} Figure~\ref{fig:resultsltpmc}(a) compares SCBS and ETIS with the 15 baselines across nine configurations. SCBS outperforms all SOF baselines in every configuration and all 15 baselines in seven configurations. ETIS improves on SCBS throughout the grid and achieves the lowest fitness in eight configurations, with only PIL outperforming it at $(\sigma^{\rm lead},\sigma^{\rm life})=(10,3)$. In Figure~\ref{fig:resultsltpmc}(b), SCBS and ETIS achieve aggregate fitness of $0.4456\%$ and $0.1903\%$, respectively, reducing fitness by $24.8\%$ and $67.9\%$ relative to PIL, the strongest baseline in aggregate.

\medskip
\noindent\underline{\textbf{Operational insights behind the representative templates.}}
SCBS allows the replenishment response to the inventory-position shortfall to be tuned through $K$, rather than filling the shortfall up to the cap as CBS does. Jointly adjusting the scale, target, and cap can advance replenishment when the system is empty and delay it when inventory is already available or scheduled to arrive. Notably, adding just one scaling parameter to CBS enables SCBS to outperform even the non-SOF PIL baseline in aggregate fitness.

ETIS combines the expiry correction in EABS with the weighted pipeline summary in IA-SCDS (Eq.~\ref{eq:resultsimbalance}), accounting for both expiring on-hand inventory and the arrival times of outstanding orders. At a fixed inventory position, the expiry correction $\omega_2x_1$ increases the replenishment signal when more inventory is approaching expiration, reflecting its reduced ability to protect against future shortages. Meanwhile, the tail adjustment $-\omega_3\varphi_{\rm tail}$ restrains replenishment when more outstanding supply is scheduled to arrive closer to the newly placed order, as these deliveries have less time to be consumed by intervening demand. By accounting for both expiration and arrival timing, ETIS can distinguish inventory states with the same aggregate inventory position that SCBS treats identically.

\subsection{Lead-time--backlog setting}
\label{sec:resultsltbacklog}
We consider nonperishable inventory with a positive replenishment lead time and deterministic demand patience, so that $\bm\sigma=(\sigma^{\rm lead},\infty,\sigma^{\rm back})$. We compare the discovered templates against eight baselines: six SOF policies and one non-SOF policy from the literature, together with the Codex baseline. The SOF baselines are Constant Order \citep{goldberg2016constant}, net Base Stock \citep{bu2026patience}, net Capped Base Stock adapted from \citet{xin2021capped}, the stockout-compensation $(s,q)$ rule (SQ) of \citet{lim2024partial}, the lot-size--reorder-point $(Q,r)$ rule (QR) of \citet{kim1985inventory}, and Generalized Base Stock (GBS) \citep{stolyar2022exploiting}. The non-SOF baseline is Myopic from \citet{bu2026patience}.

\medskip
\noindent\underline{\textbf{Representative policy templates.}} We select two representative templates: \emph{backlog-adjusted target shortfall} (BATS) and its capped variant (CBATS).
Let $\varphi_{\rm IP}=x_1+\sum_{r=1}^{\sigma^{\rm lead}}p_r-\sum_{\kappa=1}^{\sigma^{\rm back}}b_\kappa$ denote the net inventory position, and $\varphi_{\rm back}=\sum_{\kappa=1}^{\sigma^{\rm back}}b_\kappa$ the total backlog.
The following two templates share the scaled shortfall $q^{\rm short}=K\bigl(b_1+[S_1-\varphi_{\rm IP}-\varphi_{\rm back}]^+\bigr)$ and have expressions
\begin{align}
 q_{\rm BATS}&=\min\{q^{\rm short},S_2-\varphi_{\rm IP}\},\qquad 
 q_{\rm CBATS}=\min\{\xi(C,t),q^{\rm short},S_2-\varphi_{\rm IP}\},\label{eq:resultsltbbatsd}
\end{align}
BATS bounds the net-position shortfall $S_2-\varphi_{\rm IP}$ by $q^{\rm short}$, and CBATS adds the cap $\xi(C,t)$. The fitted parameters are $(S_1,S_2,K)$ for BATS, and $(S_1,S_2,K,C)$ for CBATS.

\begin{figure}[!htbp]
\centering
\begin{minipage}[b]{0.51\linewidth}
\centering
\includegraphics[width=\linewidth]{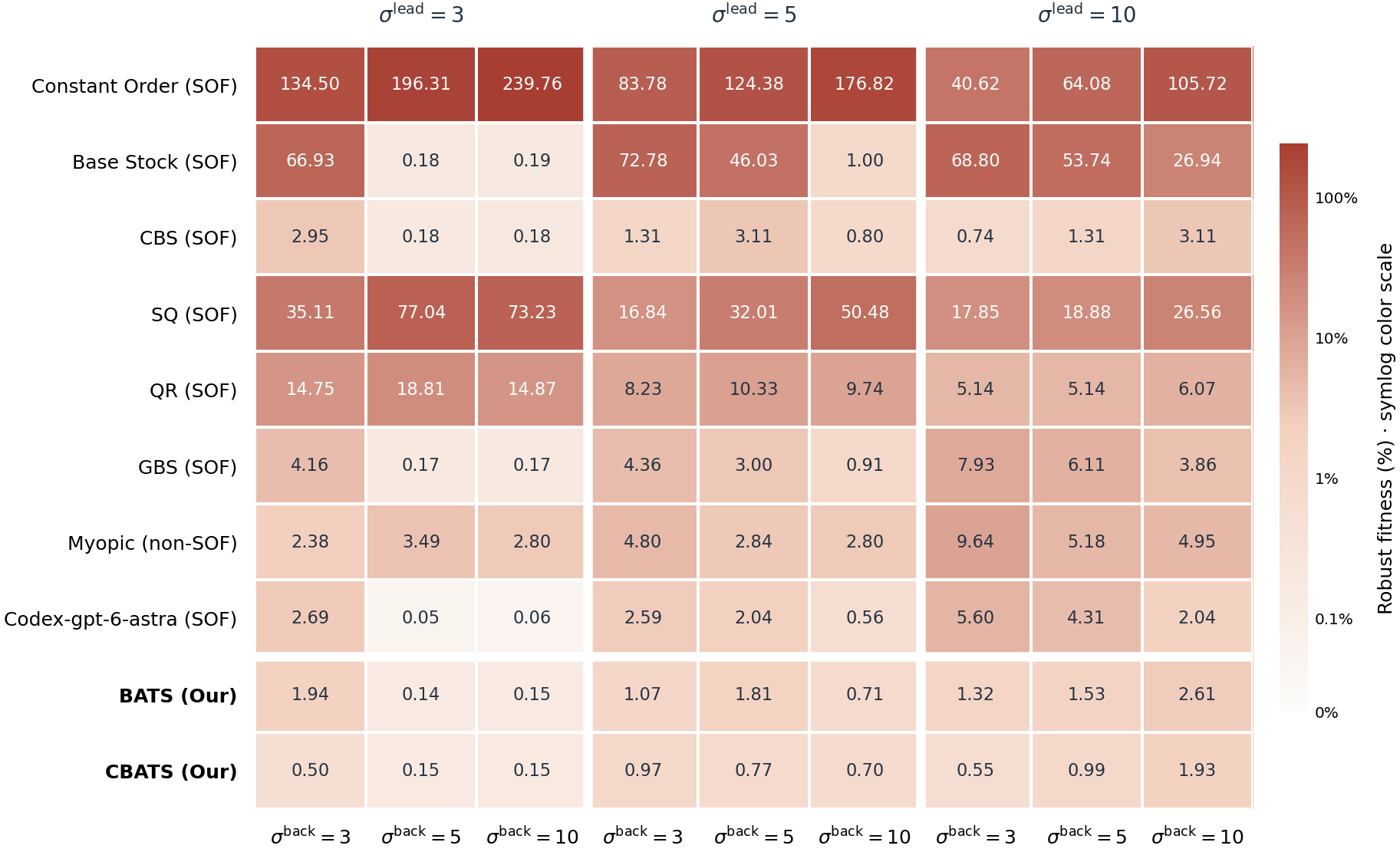}\par
{\small (a) Configuration-level fitness.\par}
\end{minipage}\hfill%
\begin{minipage}[b]{0.48\linewidth}
\centering
\includegraphics[width=\linewidth]{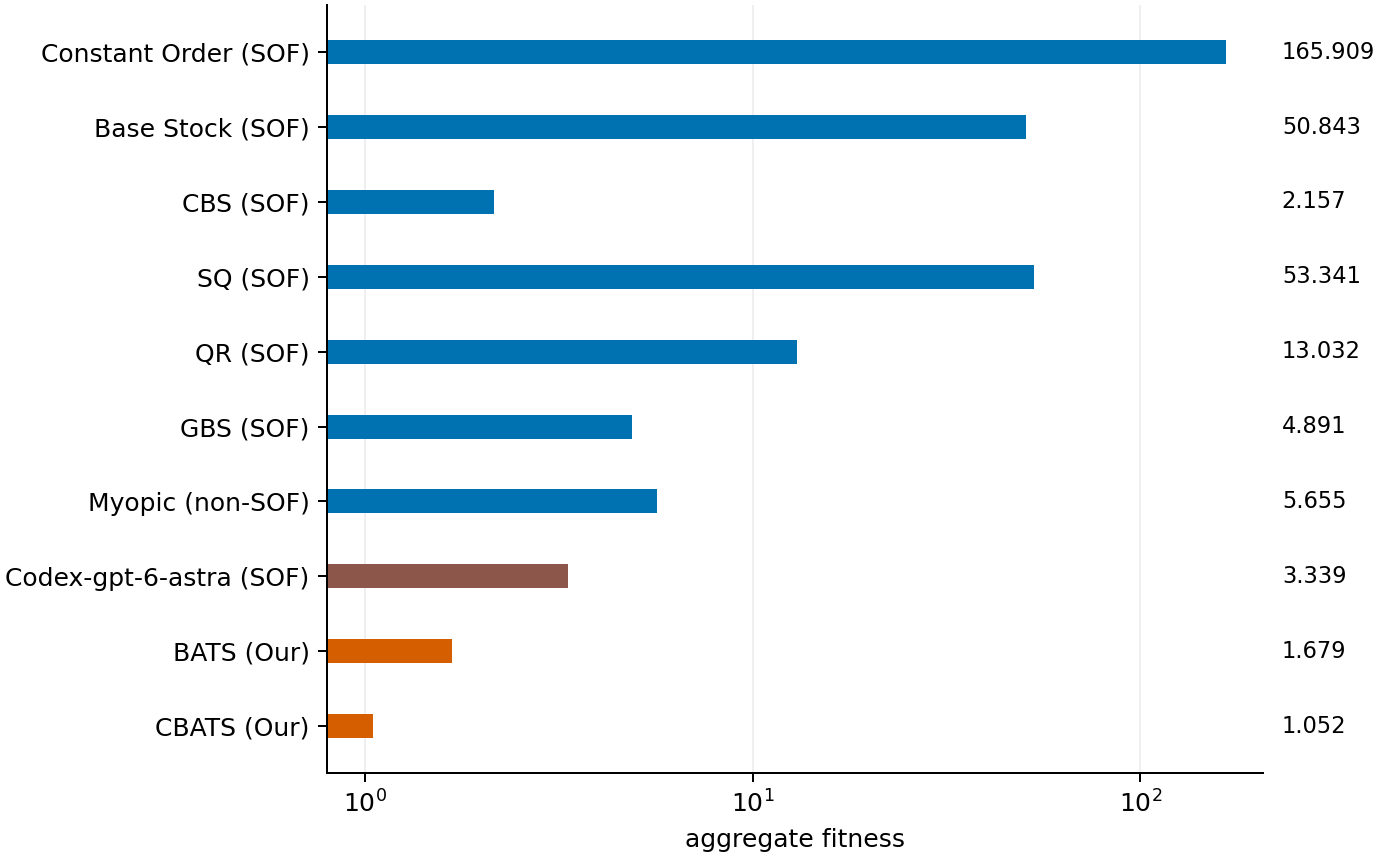}\par
{\small (b) Aggregate fitness.\par}
\end{minipage}
\caption{Benchmark fitness for lead time and backlog. Panel (a) reports configuration-level fitness by lead time and patience horizon; panel (b) reports aggregate fitness. Fitness values are percentages; lower is better.}
\label{fig:resultsltbmc}
\end{figure}

\medskip
\noindent\underline{\textbf{Evaluation of the representative templates.}} Figure~\ref{fig:resultsltbmc}(a) compares the two discovered policies with the eight baselines across nine configurations. CBATS improves on BATS in seven configurations, with slightly higher fitness at $(\sigma^{\rm lead},\sigma^{\rm back})=(3,5)$ and $(3,10)$, and outperforms all baselines in six configurations, with only Codex achieving lower fitness at $(3,5)$, $(3,10)$, and $(5,10)$. In Figure~\ref{fig:resultsltbmc}(b), BATS and CBATS achieve aggregate fitness of $1.6792\%$ and $1.0517\%$, respectively, reducing fitness by $22.1\%$ and $51.2\%$ relative to Capped Base Stock, the strongest baseline in aggregate.

\medskip
\noindent\underline{\textbf{Operational insights behind the representative templates.}}
BATS addresses the risk of excessive replenishment when the net inventory position treats all outstanding backlog as persistent demand, even though some customers may abandon before newly ordered inventory arrives. To limit this response, BATS combines the conventional net-position shortfall $S_2-\varphi_{\rm IP}$ with the adaptive cap $q^{\rm short}$. Since $\varphi_{\rm IP}+\varphi_{\rm back}=x_1+\sum_{r=1}^{\sigma^{\rm lead}}p_r$, this cap responds to the gross inventory-position ($x_1 + \sum_{r=1}^{\sigma^{\rm lead}}p_r$) shortfall and recently accumulated unmet demand $b_1$, rather than total backlog. Older backlog can therefore increase the net-position shortfall without directly increasing the adaptive cap. This structure limits replenishment driven by older backlog while preserving responsiveness to current supply conditions and recent stockouts. The scale $K$ controls the strength of this adaptive cap, while CBATS introduces the additional cap $\xi(C,t)$. Furthermore, under bounded demand, both the BATS and CBATS families contain base-stock policies as special cases. Consequently, the best achievable average cost within each family inherits the corresponding optimality guarantees established in \cite{bu2026patience}, under the assumptions of those results.

\section{Training and evaluation instances}
\label{app:experimentinstances}
This section specifies the training configurations and demand--cost profiles used for symbolic discovery, followed by the complete benchmark instance sets used for final evaluation.

\subsection{Training instances}
\label{app:searchinstances}
Each setting uses the ten demand--cost profiles in Table~\ref{tab:searchinstances}, combining standard distributions with intermittent and separated finite-support demand. Stages I--II fix the active configuration coordinate(s) at 3: Stage I searches separately on each profile, while Stage II shares an SOF across all ten. Stage III reuses the profiles at $\sigma^{\rm lead}\in\{3,5,8,10\}$ for lead time, at $\sigma^{\rm life}\in\{3,5,8,10\}$ for pure perishability, and at the six pairs $\{(2,2),(3,3),(1,4),(4,1),(5,3),(3,5)\}$ for each joint setting. These pairs denote $(\sigma^{\rm lead},\sigma^{\rm life})$ or $(\sigma^{\rm lead},\sigma^{\rm back})$, respectively, giving 40 Stage-III instances per single-feature setting and 60 per joint setting.
In Table~\ref{tab:searchinstances}, $P_\lambda$ denotes Poisson demand, $G_p$ denotes geometric demand with $\Pr(D=d)=p(1-p)^d$, and $U_k$ denotes uniform demand on $\{0,\ldots,k\}$. $Z$ denotes zero-inflated Poisson demand with rate 4 and extra zero mass 0.75. We write $D_k(\bm d;\bm p)$ for a $k$-point distribution assigning probabilities $\bm p$ to support points $\bm d$, for $k\in\{2,3,4\}$; when $\bm p$ is omitted, the probabilities are equal.

\begin{table}[H]
\centering
\caption{Ten training profiles per setting, written as $(\text{distribution},\text{cost tuple})$.}
\label{tab:searchinstances}
\label{tab:ltsearchinstances}
\label{tab:perishablesearchinstances}
\footnotesize
\setlength{\tabcolsep}{3pt}
\renewcommand{\arraystretch}{1.15}
\begin{tabular}{@{}>{\raggedright\arraybackslash}p{0.97\textwidth}@{}}
\hline
\textbf{Lead time:} $(c^{\rm hold},c^{\rm lost})$\\
$(D_2((0,4)),(1,3))$, $(P_{0.75},(1,3))$, $(P_4,(2,6))$, $(G_{1/3},(1,9))$, $(U_4,(1,4))$, $(D_2((1,9)),(2,9))$, $(Z,(1,9))$, $(D_3((0,4,12);(0.55,0.30,0.15)),(0.5,12))$, $(D_2((0,14);(0.95,0.05)),(1,20))$, $(D_3((0,2,6);(1/3,1/2,1/6)),(1,3))$.\\
\hline
\textbf{Perishability:} $(c^{\rm hold},c^{\rm lost},c^{\rm expiry})$\\
$(G_{0.2},(1,12,16))$, $(G_{0.2},(1,12,8))$, $(G_{0.2},(1,12,4))$, $(G_{0.2},(1,12,2))$, $(D_4((0,1,13,14);(0.42,0.33,0.17,0.08)),(1,9,8))$, $(U_8,(1,6,8))$, $(G_{0.25},(0.5,6,2))$, $(G_{0.25},(0.5,6,4))$, $(P_4,(2,6,0.5))$, $(D_3((0,2,8);(0.45,0.35,0.20)),(1,6,2))$.\\
\hline
\textbf{Lead time--perishability:} $(c^{\rm hold},c^{\rm lost},c^{\rm expiry})$\\
$(P_1,(1,30,3))$, $(P_1,(1,7,7))$, $(P_4,(1,30,120))$, $(P_4,(1,3,15))$, $(G_{0.5},(1,30,3))$, $(G_{0.5},(1,15,3.75))$, $(U_8,(1,30,120))$, $(U_8,(1,3,3))$, $(D_3((0,4,12)),(1,15,15))$, $(D_3((0,4,12)),(1,7,7))$.\\
\hline
\textbf{Lead time--backlog:} $(c^{\rm hold},c^{\rm backlog},c^{\rm lost})$\\
$(P_1,(1,0.75,3))$, $(G_{0.5},(1,0.75,3))$, $(D_3((0,4,12)),(1,0.3,3))$, $(D_3((0,4,12)),(1,7.5,30))$, $(P_4,(1,0.75,3))$, $(P_4,(1,3.75,15))$, $(G_{0.5},(1,3,0.3))$, $(G_{0.5},(1,0.03,12))$, $(U_8,(1,1.75,7))$, $(D_3((0,4,12)),(1,3.75,15))$.\\
\hline
\end{tabular}
\end{table}

\subsection{Evaluation instances}
\label{app:demandcatalog}
All four settings use the same 50-distribution catalog $\mathfrak D^{\rm bench}$. Unbounded demand distributions are represented by $\min\{D,50\}$, assigning all probability mass above 50 to demand 50; the distribution parameters listed below describe the distributions before this tail folding. Table~\ref{tab:evaluationinstances} gives the configuration and nonzero cost grids combined with that catalog; cost coefficients not listed are zero. Each row therefore defines the complete benchmark set for one setting.

\begin{table}[H]
\centering
\caption{Evaluation configurations and cost profiles. The configurations have the forms $(\sigma^{\rm lead},\infty,0)$, $(0,\sigma^{\rm life},0)$, $(\sigma^{\rm lead},\sigma^{\rm life},0)$, and $(\sigma^{\rm lead},\infty,\sigma^{\rm back})$, respectively.}
\label{tab:evaluationinstances}
\footnotesize
\setlength{\tabcolsep}{4pt}
\renewcommand{\arraystretch}{1.2}
\begin{tabular}{@{}>{\raggedright\arraybackslash}p{0.21\textwidth}>{\raggedright\arraybackslash}p{0.21\textwidth}>{\raggedright\arraybackslash}p{0.42\textwidth}>{\raggedright\arraybackslash}p{0.09\textwidth}@{}}
\hline
Setting & Configurations & Nonzero cost coefficients & Instances \\
\hline
Lead time & $\sigma^{\rm lead}\in\{2,\ldots,15\}$ & $c^{\rm hold}=1$, $c^{\rm lost}\in\{1,3,9,20\}$ & 2,800 \\
Perishability & $\sigma^{\rm life}\in\{2,\ldots,15\}$ & $c^{\rm hold}=1$; $(c^{\rm lost},c^{\rm expiry})\in$ $\{(1,0.25),\allowbreak(1,16),\allowbreak(3,1),\allowbreak(3,4),\allowbreak(9,0.25),\allowbreak(9,16),\allowbreak(20,1),\allowbreak(20,4)\}$ & 5,600 \\
Lead time--perishability & $(\sigma^{\rm lead},\sigma^{\rm life})\in\{3,5,10\}^2$ & $c^{\rm hold}=1$; same eight $(c^{\rm lost},c^{\rm expiry})$ pairs as above & 3,600 \\
Lead time--backlog & $(\sigma^{\rm lead},\sigma^{\rm back})\in\{3,5,10\}^2$ & $c^{\rm hold}=1$; $(c^{\rm backlog},c^{\rm lost})\in$ $\{(0.05,0.25),\allowbreak(0.05,1),\allowbreak(0.25,0.25),\allowbreak(0.25,1),\allowbreak(1,5),\allowbreak(1,20),\allowbreak(5,5),\allowbreak(5,20)\}$ & 3,600 \\
\hline
\end{tabular}
\end{table}

\medskip
\noindent\textbf{Demand catalog.} Each row below specifies five distributions, giving $|\mathfrak D^{\rm bench}|=50$. Finite-support probabilities follow the listed support order; omitted probabilities are equal. Geometric and negative-binomial supports start at zero.
\par\smallskip
{\footnotesize
\setlength{\tabcolsep}{3pt}
\renewcommand{\arraystretch}{1.15}
\noindent\begin{tabular}{@{}>{\raggedright\arraybackslash}p{0.23\textwidth}>{\raggedright\arraybackslash}p{\dimexpr0.77\textwidth-6pt\relax}@{}}
\hline
\textbf{Family} & \textbf{Five parameter choices or supports}\\
\hline
Poisson & $\lambda\in\{0.5,1,2,4,8\}$.\\
Geometric & $\mu\in\{0.5,1,2,4,8\}$; $p=1/(1+\mu)$.\\
Negative binomial & Mean--shape pairs $(\mu,r)\in\{(1,1),(2,1),(2,2),(4,2),(8,4)\}$; $p=r/(r+\mu)$.\\
Binomial & Trial count and success probability $(n,p)\in\{(2,0.5),(4,0.5),(8,0.25),(8,0.5),(16,0.5)\}$.\\
Discrete uniform & Integer supports $[0,2]$, $[0,4]$, $[0,8]$, $[1,5]$, $[2,10]$.\\
Two-point & $D_2((0,4))$, $D_2((0,8))$, $D_2((1,9))$, $D_2((0,10);(0.8,0.2))$, $D_2((2,12);(0.2,0.8))$.\\
Three-point & $D_3((0,2,6))$, $D_3((0,4,12))$, $D_3((1,3,9);(0.2,0.6,0.2))$, $D_3((1,5,15);(0.5,0.3,0.2))$, $D_3((2,6,12);(0.2,0.3,0.5))$.\\
Intermittent shifted Poisson & $\Pr(D=0)=z\in\{0.5,0.7,0.8,0.9,0.95\}$; otherwise $D\sim1+\operatorname{Poisson}(z/(1-z))$.\\
Poisson mixture & Component rates and first-component weight $(\lambda_1,\lambda_2,w)\in\{(0.5,4,0.5),(0.5,8,0.5),(1,6,0.3),(1,10,0.7),(2,10,0.5)\}$.\\
Random finite PMF & One fixed symmetric Dirichlet draw on each support $\{0,\ldots,m\}$, with maximum--concentration pairs $(m,a)\in\{(5,0.2),(8,0.3),(10,0.5),(12,1),(15,3)\}$.\\
\hline
\end{tabular}\par}

\end{document}